\documentclass{article}

\usepackage[utf8]{inputenc}
\usepackage{amsmath,amsfonts,amssymb}
\usepackage{amsthm,upref}
\usepackage[a4paper]{geometry}
\usepackage{color} 
\usepackage[dvipsnames]{xcolor}
\usepackage{graphicx}
\usepackage{subfigure}
\usepackage{tabularx,float}
\usepackage[hidelinks]{hyperref}
\usepackage{comment} 
\usepackage{tabularx} 
\usepackage[pagewise]{lineno}
\usepackage{titlefoot} 
\usepackage{fancyhdr}

\newtheorem{thm}{Theorem}[section]
\newtheorem{definition}[thm]{Definition}

\newtheorem{proposition}[thm]{Proposition}
\newtheorem{corollary}[thm]{Corollary}
\newtheorem{lemma}[thm]{Lemma}
\newtheorem{remark}[thm]{Remark}
\newtheorem{assumption}{Assumption}

\newcommand{\eps}{\varepsilon}

\newcommand{\imag}{\textrm{Im}\,}

\newcommand\norm[1]{\left\lVert#1\right\rVert}

\title{Discrete Geometric Singular Perturbation Theory}
\author{S.~Jelbart\thanks{Department of Mathematics, The Technical University of Munich. Garching Bavaria, 85748, Germany.} \& C.~Kuehn$^\ast$}
\date{\today}

\begin{document}

	\maketitle
	

	\begin{abstract}
		We propose a mathematical formalism for discrete multi-scale dynamical systems induced by maps which parallels the established \textit{geometric singular perturbation theory} for continuous-time fast-slow systems. We identify limiting maps corresponding to both `fast' and `slow' iteration under the map. A notion of normal hyperbolicity is defined by a spectral gap requirement for the multipliers of the fast limiting map along a critical fixed-point manifold $S$. We provide a set of Fenichel-like perturbation theorems by reformulating pre-existing results so that they apply near compact, normally hyperbolic submanifolds of $S$. The persistence of the critical manifold $S$, local stable/unstable manifolds $W^{s/u}_{loc}(S)$ and foliations of $W^{s/u}_{loc}(S)$ by stable/unstable fibers is described in detail. The practical utility of the resulting \textit{discrete geometric singular perturbation theory (DGSPT)} is demonstrated in applications. First, we use DGSPT to identify singular geometry corresponding to excitability, relaxation, chaotic and non-chaotic bursting in a map-based neural model. Second, we derive results which relate the geometry and dynamics of fast-slow ODEs with non-trivial time-scale separation and their Euler-discretized counterpart. Finally, we show that fast-slow ODE systems with fast rotation give rise to fast-slow Poincar\'e maps, the geometry and dynamics of which can be described in detail using DGSPT.
	\end{abstract}
	
	
	\unmarkedfntext{\textbf{Keywords:} geometric singular perturbation theory, discrete dynamical systems, multi-scale dynamical systems, singularly perturbed maps, invariant manifolds.}
	
	\unmarkedfntext{\noindent \textbf{MSC2020:} 37C05, 37D10, 37C86, 34D15, 37C15.}

	\section{Introduction}
	\label{sec:introcuction}

	The primary aim of this manuscript is to provide a mathematical framework for the geometric analysis of multi-scale discrete dynamical systems induced by maps. In essence, we aim to provide a \textit{discrete geometric singular perturbation theory}, or simply \textit{(DGSPT)}, which parallels the established \textit{geometric singular perturbation theory (GSPT)} for continuous-time fast-slow systems \cite{Fenichel1979,Jones1995,Kaper1999,Kuehn2015,Wechselberger2019}. Our starting point is motivated by the recently developed formulation in \cite{Wechselberger2019} of GSPT for continuous-time, $C^r-$smooth fast-slow ODE systems in the general form
	\begin{equation}
		\label{eq:nonstnd_form}
		z' = N(z) f(z) + \eps G(z,\eps) ,
	\end{equation}
	where $z \in \mathbb R^n$, $z' = dz / dt$, $N$ is an $n \times (n-k)$ matrix, $f$ and $G$ are vector-valued functions of dimensions $(n-k) \times 1$ and $n \times 1$ respectively, $0 < \eps \ll 1$ is a perturbation parameter and the existence of a $(k < n)-$dimensional \textit{critical manifold} $S = \{z \in \mathbb R^n : f(z) = 0\}$ for $\eps = 0$ is assumed. The class of systems defined by \eqref{eq:nonstnd_form} includes the (perhaps better known) class of fast-slow systems in the so-called \textit{standard form}
	\begin{equation}
		\label{eq:stnd_form_ode}
		\begin{split}
			x' &= \eps \tilde g(x,y,\eps) , \\
			y' &= \tilde f(x,y,\eps) ,
		\end{split}
	\end{equation}
	with $(x,y) \in \mathbb R^k \times \mathbb R^{n-k}$ as a special case, since \eqref{eq:stnd_form_ode} can be written in the form \eqref{eq:nonstnd_form} after setting $z = (x,y)^\top$, $f(z) = \tilde f(x,y,0)$, $N(z) = (O_{k,n-k}, I_{n-k})^\textnormal{T}$ and $G(z,\eps) = (\tilde g(x,y,\eps), \tilde f(x,y,\eps) - \tilde f(x,y,0))^\textnormal{T}$. Conceptually, the formulation presented in \cite{Wechselberger2019} can be considered as a coordinate-independent extension of earlier formulations for standard form systems \eqref{eq:stnd_form_ode} in e.g.~\cite{Jones1995,Kaper1999,Kuehn2015}, motivated by a large number of applications for which fast-slow dynamics arises in systems that are `beyond the standard form', see e.g.~\cite{deMaesschalck2021b,Huber2005,Jelbart2022,Jelbart2020a,Kaleda2011,Gucwa2009,Kosiuk2011,Kosiuk2016,KuehnSzmolyan2015,Schecter2009}.
	
	\
	
	For the development of DGSPT, we consider fast-slow maps in the general form
	\begin{equation}
		\label{eq:nonstandard_maps}
		z \mapsto \bar z = z + N(z) f(z) + \eps G(z,\eps) ,
	\end{equation}
	where $z \in \mathbb R^n$, $N$, $f$ and $G$ are defined as in the continuous-time system \eqref{eq:nonstnd_form}, $0 < \eps \ll 1$ is a perturbation parameter and the existence of a $(k < n)-$dimensional critical (fixed point) manifold $S = \{z \in \mathbb R^n : f(z) = 0\}$ for $\eps = 0$ is assumed. Similarly to continuous-time setting, the class of maps defined by \eqref{eq:nonstandard_maps} includes the (perhaps better known) class of fast-slow maps in standard form
	\begin{equation}
		\label{eq:stnd_form_maps}
		\begin{split}
			x &\mapsto \bar x = x + \eps \tilde g(x,y,\eps) , \\
			y &\mapsto \bar y = y + \tilde f(x,y,\eps) ,
		\end{split}
	\end{equation}
	where $(x,y) \in \mathbb R^k \times \mathbb R^{n-k}$ as a special case, since \eqref{eq:stnd_form_maps} can be written in the form \eqref{eq:nonstandard_maps} using with the same choices for $N$, $f$ and $G$ which put system \eqref{eq:stnd_form_ode} into the form \eqref{eq:nonstnd_form}.
	
	\
	
	The utility of GSPT in the continuous-time setting 
	depends upon the availability of both
	\begin{itemize}
		\item[(I)] \textit{Singular theory}: A mathematical framework or formalism for the geometric analysis of non-equivalent limiting problems for each `time-scale' as $\eps \to 0$, and
		\item[(II)] \textit{Perturbation/invariant manifold theorems}: A collection of results on the perturbation of dynamical and geometric structure identified using the singular theory in (I) for $0 < \eps \ll 1$,
	\end{itemize}
	and this work, we shall consider (I) and (II) together as the basic requirements for a prospective DGSPT.
	
	In the continuous-time setting, a singular theory as required by (I) was already present in \cite{Fenichel1979}, later refined and clarified for fast-slow systems in standard form \eqref{eq:stnd_form_ode} in e.g.~\cite{Guckenheimer1996,Jones1995,Kaper1999,Kuehn2015,Szmolyan1991}, and finally for the more general class \eqref{eq:nonstnd_form} in \cite{Wechselberger2019} and related works, e.g.~\cite{deMaesschalck2021,Goeke2014,Kruff2019,Lizarraga2020c,Lizarraga2020b,Lizarraga2020}. This theory provides a mathematical framework for identifying and analysing the geometry and dynamics of the so-called \textit{layer (or fast subsystem)} and \textit{reduced (or slow subsystem)} problems, obtained after taking $\eps \to 0$ in \eqref{eq:nonstnd_form} on the fast and slow time-scales $t$ and $\tau = \eps t$ respectively. A typical analysis proceeds by constructing `singular orbits', sometimes also called candidate orbits, by a geometric concatenation of fast trajectory segments of the layer problem on $\mathbb R^n \setminus S$ and slow trajectory segments of the reduced problem on $S$.
	
	The fundamental perturbation theorems required by (II) are provided by \textit{Fenichel theory}, which is constituted by a collection of perturbation and invariant manifold theorems derived in \cite{Fenichel1974,Fenichel1977,Fenichel1971} and culminating in \cite{Fenichel1979}
	(though it is also important to mention the earlier works of \cite{Hirsch1970,Tikhonov1952}, which provided foundational understanding for the field). Fenichel theory 
	ensures that suitably constructed singular orbits 
	perturb in a regular fashion for $0 < \eps \ll 1$ in regions of phase space bounded away from certain singularities on $S$. Fenichel theory applies locally near \textit{normally hyperbolic} submanifolds of $S$ for which the linearization along the ($k-$dimensional) critical manifold $S$ has $n-k$ non-trivial eigenvalues bounded away from the imaginary axis. For fast-slow systems \eqref{eq:nonstnd_form}, this is equivalent to the requirement that the eigenvalues of the $(n-k)\times(n-k)$ square matrix $DfN|_S$ are bounded away from the imaginary axis \cite{Wechselberger2019}. Under normally hyperbolic conditions, the geometry and dynamics for $0 < \eps \ll 1$ are described up to $O(\eps)$ accuracy by the dynamics of layer and reduced problems, which are almost invariably far more tractable analytically. 
	Additional perturbation theorems are required in order to describe the geometry and dynamics near non-normally hyperbolic submanifolds of $S$. Many authors have demonstrated the utility of a method of geometric desingularization known as the \textit{blow-up method} for such purposes; here we simply cite the seminal works in 
	\cite{Dumortier1996,Krupa2001a} and refer to the recent survey~\cite{Jardon2019b}. 
	This combination of GSPT and blow-up techniques has been applied by many authors \cite{Carter2018,Wechselberger2015,Hayes2016,Huber2005,Jardon2019,Jelbart2022,Gucwa2009,Kosiuk2011,Kosiuk2016,KuehnSzmolyan2015,Szmolyan2004}.

	
	\
	
	In the discrete setting, the picture is less complete. To the best of our knowledge, a singular theory in the sense of (I) does not yet exist for maps. This is perhaps because of difficulties relating to the fact that there is no direct analogue for the time rescaling $\tau = \eps t$ which leads to an equivalent `slow formulation' of the map \eqref{eq:nonstandard_maps}. Thus it is not immediately clear how a to obtain a `reduced map' which describes slow iteration on or close to $S$. On the other hand, the perturbation and invariant manifold theory required by (II) for the (discrete analogue of the) normally hyperbolic regime is in principle quite established, dating back at least to the work of Hirsch, Pugh \& Shub \cite{Hirsch1970}, and results on the existence of invariant manifolds and the foliation of the adjacent space in 2-dimensional maps arising in the analysis of fast-slow ODE systems in particular are given in \cite{Szmolyan2004}. We also mention Pötzsche \cite{Potzsche2003}, Nipp \& Stoffer \cite{Nipp1992,Nipp2013} and 
	Shub \cite{Shub2013} (many more references can be found in the books \cite{Nipp2013,Shub2013}). In many practical situations, however, a direct application of the pre-existing results of the results in \cite{Hirsch1970} can be difficult due to their generality and the relatively abstract formulation of necessary and sufficient conditions of their applicability. On the other hand, concrete and more applicable formulations such as those in Nipp \& Stoffer \cite{Nipp2013} depend upon the identification of suitable coordinates, and often require the rather extensive use of nonlinear coordinate transformations and cutoff techniques in order to `prepare' the equations. Thus in many situations, there remains a practical barrier to the application of these results. With regard to perturbation results on the dynamics in the non-normally hyperbolic setting we mention \cite{Arcidiacono2019,Baesens1991,Do2013,Engel2020b,Engel2019,Engel2019b,Fruchard2009,Nipp2013,Nipp2009}.
	
	\
	
	In order to obtain a satisfactory DGSPT, our first task is to develop a singular theory in the sense of (I) for fast-slow maps \eqref{eq:nonstandard_maps}. We begin by defining a \textit{layer map} by setting $\eps = 0$ in \eqref{eq:nonstandard_maps}, which allows us introduce a notion of normal hyperbolicity of $S$ in terms of a spectral gap requirement. Algebraic formulae for the non-trivial multipliers needed to determine the normal hyperbolicity of $S$ are given solely in terms of the initial data $N$ and $f$. Specifically, 
	it suffices to verify the existence of an annular spectral gap about the unit circle for the matrix $I_{n-k} + Df N|_S$.
	
	Next, we show that a \textit{reduced map} which approximates slow iteration near $S$ to an accuracy of $O(\eps^2)$ can be derived under normally hyperbolic conditions. The reduced map is conceptually distinct from the reduced problem in the continuous-time setting since it reduces to the trivial map $z \mapsto z$ as $\eps \to 0$. This is necessarily so, since the $\eps-$dependence in the leading order cannot in general be `divided out' (as is achieved by moving to the slow time-scale $\tau = \eps t$ in the continuous-time setting). Fortunately, however, this is no obstacle in practice, where one is primarily concerned with understanding the leading order dynamics near $S$. We also show that a reduced $m$'th iterate map induced by repeated iteration of \eqref{eq:nonstandard_maps} can be derived locally near $S$ using the asymptotic self-similarity properties of \eqref{eq:nonstandard_maps}. Interestingly, this $m$'th iterate map can be related to a suitable discretization of the continuous-time reduced problem associated to fast-slow ODEs \eqref{eq:nonstnd_form} if the number of iterates $m$ is comparable to $\eps^{-1}$. 
	From an applied point of view, it is significant that both reduced and $m$'th iterate maps are given by closed form algebraic formulae defined purely in terms of $N$, $f$ and $G$.
	
	Given the prevalence of fast-slow maps in standard form \eqref{eq:stnd_form_maps} in applications, we also consider the form of the corresponding singular theory as a special case of the theory developed for general fast-slow maps \eqref{eq:nonstandard_maps}. In addition to describing general features of the layer and reduced maps, we prove local equivalence of the maps \eqref{eq:nonstandard_maps} and \eqref{eq:stnd_form_maps} near an arbitrary point on $S$. Just as in the continuous-time setting, however, it is worthy to emphasise such an equivalence is strictly local, and typically only useful for theoretical purposes.
	
	\
	
	Having developed a singular theory, we turn our attention to the coupling of this singular theory to suitable perturbation and invariant manifold theorems as required by (II). This is achieved via the adaptation of pre-existing perturbation and invariant manifold theorems in the formulation of \cite{Nipp2013} for fast-slow maps \eqref{eq:nonstandard_maps}. We provide persistence theorems which parallel Fenichel's invariant manifold theorem's for flows in continuous-time fast-slow systems, which apply under normally hyperbolic conditions in the discrete sense described above. These results characterise the perturbation of compact normally hyperbolic submanifolds of the critical manifold $S$, as well as its corresponding local stable and unstable manifolds $W_{loc}^s(S)$ and $W_{loc}^u(S)$ respectively for $0 < \eps \ll 1$. Local invariance and smoothness properties of the perturbed counterparts are described, as well as the asymptotic rate foliation of perturbed stable and unstable manifolds by stable and unstable fast fibers, respectively.
	
	Although most of the invariant manifold theorems presented herein have direct and in most cases more general analogues 
	in \cite{Nipp2013} and other pre-existing work dating back to \cite{Hirsch1970}, our main contribution is to extend and reformulate these results in a manner well-suited to applications. We emphasise in particular the following:
	\begin{itemize}
		\item As they are formulated herein, perturbation theorems do not depend on a special choice of coordinates, so the equations do not need to be `prepared'.
		\item In applications it suffices to check normal hyperbolicity of $S$, which amounts to calculating the eigenvalues of the matrix $DfN|_S$.
		\item Results apply to compact, normally hyperbolic submanifolds $S_n \subseteq S$. Perturbed counterparts of $S$, $W^{s/u}_{loc}(S)$ and foliations for $W^{s/u}_{loc}(S)$ are typically locally (as opposed to globally) invariant objects.
	\end{itemize}
	In essence, our formulation leads to perturbation theorems which parallel Fenichel's theorems in the continuous-time setting. Our main perturbation theorems are obtained by formulating a number of necessary and sufficient conditions in \cite{Nipp2013} in terms of spectral bounds for the multipliers of the layer map along $S$, which do not depend on a special choice of coordinates. This allows us to derive results for the map \eqref{eq:nonstandard_maps} via the application of results in \cite{Nipp2013} to a suitable `normal form'. The `price' of coordinate-independence in our approach, is that it requires a sufficient degree of smoothness in the map \eqref{eq:nonstandard_maps} (it must be at least $C^1$ in order for the spectrum to be well-defined), while many results in \cite{Nipp2013} apply minimally under Lipschitz continuous conditions. The extension to the locally (as opposed to globally) invariant case is delicate but standard, and achieved by the use of cutoff techniques.
	
	
	In addition to the reformulation of results in \cite{Nipp2013} as coordinate-independent results 
	in the locally invariant setting, we also provide concrete results for the map \eqref{eq:nonstandard_maps}, including explicit formulae for the perturbed slow manifold $S_\eps$ in local coordinates up to $O(\eps^2)$, and quantitative estimates for contraction and repulsion along stable and unstable fibers respectively in terms of the size of the spectral gap associated to the matrix $I_{n-k} + DfN|_S$.
	
	\
	
	The utility of DGSPT 
	is demonstrated in the context of three different applications. The first of these is a 2-dimensional map-based model for neuronal bursting known as the \textit{Chialvo map}, introduced in \cite{Chialvo1995} and considered further in e.g.~\cite{Jing2006,Muni2022,Trujillo2021,Wang2018}. This map takes the standard form \eqref{eq:stnd_form_maps} in a suitable parameter regime, 
	with an S-shaped critical manifold having two regular fold points and a flip-type (period-doubling) bifurcation in the layer map appearing as non-normally hyperbolic singularities on $S$. We show how DGSPT can be used to identify four open parameter sets corresponding to excitable, relaxation-type, non-chaotic bursting and (potentially chaotic) bursting dynamics. This analysis extends (in some directions) the work in \cite{Chialvo1995,Jing2006,Trujillo2021}, though it remains only partial since the dynamics near the non-normally hyperbolic fold points is not yet understood in detail.
	
	The last two applications demonstrate the utility of DGSPT in a more theoretical setting. In the first of these, we consider the geometry and dynamics of maps arising by Euler discretization of fast-slow systems in the general non-standard form \eqref{eq:nonstnd_form}. These results, which relate the geometry and dynamics of the ODE and discretized systems, parallel pre-existing results on discretized fast-slow systems for larger classes of discretizations in e.g.~\cite{Nipp1995,Nipp1996}. We restrict to the simplest case of Euler discretized systems in order to show clearly how DGSPT applies in such contexts. To the best of our knowledge, the extension of these results to case of fast-slow systems in the more general form \eqref{eq:nonstandard_maps} is also novel.
	
	Finally, we use DGSPT in order to analyse fast-slow Poincar\'e maps associated to fast-slow systems in standard form with a single slow variable, a situation which arises often in applications if one allows a parameter $\alpha$ to evolve slowly in time; see e.g.~\cite{Benoit2006,Fruchard2009} and the references therein for examples in the context of dynamic bifurcation theory. On the assumption that the (continuous-time) layer problem has a hyperbolic limit cycle for $\alpha = \alpha_\ast$, there exists a 2-dimensional manifold of limit cycles $\mathcal M$ in $\mathbb R^n \times V_\alpha$, where $V_\alpha$ is a sufficiently small neighbourhood about $\alpha_\ast$ (see e.g.~\cite{Stiefenhofer1998} for an example in 3 dimensions). After showing that the Poincar\'e map on a transverse section $\Delta$ is a fast-slow map with normally hyperbolic critical manifold $S = \Delta \cap \mathcal M$, we are able to characterise the geometry and dynamics of the Poincar\'e map in detail using DGSPT. Using information about the Poincar\'e map we can then infer geometric and dynamical properties for the higher-dimensional ODE system. In particular, we extend previous results due to Anosova \cite{Anosova1999,Anosova2002} which characterise the persistence of $\mathcal M$ as a nearby locally invariant manifold $\mathcal M_\eps$.
	
	\
	
	The manuscript is structured as follows: In Section \ref{sec:a_coordinate-independent_framework_for_fast-slow_maps} we develop the singular theory. Layer and reduced maps are introduced and characterised in Sections \ref{sub:the_layer_map} and \ref{sub:reduced_map} respectively, and the relationship to the special subclass of fast-slow maps in standard form \eqref{eq:stnd_form_maps} is considered in Section \ref{sub:fast-slow_maps_in_standard_form}. The main invariant manifold theorems are stated and described in Section \ref{sec:slow_manifold_theorems}, and proved in Section \ref{sec:proof_of_theorem_fenichel}. The applications are treated in Section \ref{sec:examples}; Section \ref{sub:Chialvo_model} contains the geometric analysis of the map-based neural model, Section \ref{sub:Euler_discretization} contains the analysis of Euler discretized fast-slow ODEs, and the application to fast-slow Poincar\'e maps induced by ODEs with slowly varying parameters and persistence results for limit cycle manifolds are given in Section \ref{sub:Poincare_maps}. Finally in Section \ref{sec:outlook}, we summarise and conclude the manuscript.

	\section{Coordinate-independent GSPT for fast-slow maps}
	\label{sec:a_coordinate-independent_framework_for_fast-slow_maps}

	In this section we extend the singular GSPT framework developed for fast-slow ODEs in \cite{Wechselberger2019} to fast-slow maps. Specifically, we consider maps of the form
	\begin{equation}
		\label{eq:general_map_1}
		z \mapsto \bar z = H(z,\eps) ,
	\end{equation}
	with variables $z \in \mathbb R^n$ and perturbation parameter $\eps$. The map $H : \mathbb R^n \times [0,\eps_0) \to \mathbb R^n$ is assumed to be $C^r-$smooth in both $z$ and $\eps$. For simplicity, we shall assume that $r$ is sufficiently large for the validity of certain calculations. This will frequently lead to simplified statements, however most results can be derived and stated minimally for $r \geq 1$.
	
	Since \eqref{eq:general_map_1} is (at least) $C^1-$smooth in $\eps$, the map can be written as
	\begin{equation}
		\label{eq:general_map_2}
		z \mapsto \bar z = z + h(z) + \eps G(z,\eps) ,
	\end{equation}
	where the functions $h:\mathbb R^n \to \mathbb R^n$ and $G:\mathbb R^n \times [0,\eps_0) \to \mathbb R^n$ are $C^r-$smooth in $z$ and (for the latter) $C^{r-1}-$smooth in $\eps$.

	\begin{definition}
		\label{def:singular_pert_map}
		\textup{(Regularly perturbed maps/fast-slow maps)} The map \eqref{eq:general_map_2} is called a fast-slow map if the level set
		\begin{equation}
			\label{eq:critical_set}
			C:= \left\{z \in \mathbb R^n : h(z) = O_n \right\} ,
		\end{equation}
		where $O_n$ denotes the $n \times 1$ zero vector, contains a $(k<n)-$dimensional regularly embedded submanifold $S \subset \mathbb R^n$. We also say that the map \eqref{eq:general_map_2} is singularly perturbed. If the above condition is not satisfied, we say that the map \eqref{eq:general_map_2} is regularly perturbed.
	\end{definition}
	
	Definition \ref{def:singular_pert_map} is directly analogous to the geometric definition of singularly perturbed ODEs in \cite{Fenichel1979}, see also \cite[Definitions 3.1-3.2]{Wechselberger2019}.
	In this work we are interested in the dynamics of fast-slow maps, so we restrict our focus accordingly.
	
	\begin{assumption}
		\label{ass:1}
		\textup{(Restriction to fast-slow maps)} 
		The map \eqref{eq:general_map_2} is a fast-slow map in the sense of Definition \ref{def:singular_pert_map}. We also assume for simplicity that the submanifold $S \subseteq C \subset \mathbb R^n$ in Definition \ref{def:singular_pert_map} is connected. 
	\end{assumption}
	
	As stated, the assumption that $S$ is connected is made for simplicity. Generalisations and adaptions of all results to cases in which $S$ is a disjoint union of connected submanifolds are straightforward.
	

	\subsection{The layer map}
	\label{sub:the_layer_map}
	
	Since \eqref{eq:general_map_1} is smooth (and thus continuous) in $\eps$, the limit $\eps \to 0$ is well-defined.
	
	\begin{definition}
		\textup{(Layer map)} The map
		\begin{equation}
			\label{eq:layer_map}
			z \mapsto \bar z = z + h(z) 
		\end{equation}
		obtained from \eqref{eq:gen_maps} in the limit $\eps \to 0$ is called the layer map.
	\end{definition}
	

	It follows from Assumption \ref{ass:1} that the set of fixed points $C$ contains a $k-$dimensional submanifold $S$. In order that this submanifold $S$ itself has a level set representation, we impose an algebraic assumption on the existence of a suitable factorisation of $h(z)$. Such assumptions are also made in the ODE setting in \cite{Wechselberger2019}.
	
	\begin{assumption}
		\label{ass:factorisation}
		\textup{(Factorisation of the layer map)} 
		The function $h(z)$ can be factorised as follows:
		\begin{equation}
			\label{eq:factorisation}
			h(z) = N(z) f(z) ,
		\end{equation}
		where $N(z)$ and $f(z)$ are matrices of size $n \times (n-k)$ and $(n-k) \times 1$ respectively. We assume the matrix $N(z)$ has full column rank for all $z \in S$, and that fixed points $z_\ast \notin S$ such that $N(z_\ast)f(z_\ast) = O_n$, if they exist, are isolated.
	\end{assumption}
	
	Using Assumption \ref{ass:factorisation}, we shall hereafter write the map \eqref{eq:general_map_2} in the form
	\begin{equation}
		\label{eq:gen_maps}
		z \mapsto \bar z = z + N(z) f(z) + \eps G(z,\eps) .
	\end{equation}
	All subsequent results will be stated for fast-slow maps in the form \eqref{eq:gen_maps} satisfying Assumptions \ref{ass:1}-\ref{ass:factorisation}.
	
	\begin{remark}
		\label{rem:factorisation}
		The factorisation \eqref{eq:factorisation} can always be locally attained. In fact, algebraic algorithms for obtaining such a local factorisation already exist \cite{Goeke2014,Lizarraga2020}. Globally, however, it must be assumed. See \cite[Rem.~3.6]{Wechselberger2019} for discussion and \cite[Sec.~1.3]{deMaesschalck2021} for counterexamples in the continuous-time setting. 
	\end{remark}
	
	It follows from Assumption \ref{ass:factorisation} that the $k-$dimensional manifold of fixed points $S \subseteq C$ is given by the level set
	\begin{equation}
		\label{eq:S}
		S := \left\{ z\in\mathbb R^n : f(z) = O_{n-k} \right\} .
	\end{equation}
	In the theory for fast-slow ODEs, $S$ is known as the \textit{critical manifold}. We adopt the same terminology order to emphasise the similar role played by the fixed point manifold \eqref{eq:S} in fast-slow maps.
	
	\begin{definition}
		\textup{(Critical manifold)} The fixed point manifold $S$ is called the critical manifold of \eqref{eq:gen_maps}.
	\end{definition}
	
	Since $S$ is a regularly embedded submanifold of $\mathbb R^n$ by Assumption \ref{ass:1}, the $(n-k) \times n$ matrix $Df|_S$ is regular, i.e.~has full row rank.
	
	\begin{remark}
		It is common in applications that the set $S$ is not everywhere a regularly embedded submanifold in $\mathbb R^n$, since it contains, e.g.~self-intersections. We do not consider such cases in this work; they are ruled out by Assumption \ref{ass:1}.
	\end{remark}
	
	In order to describe stability properties of the critical manifold, we calculate the Jacobian
	\begin{equation}
		\label{eq:Jacobian}
		DH(z,0) = I_n + Dh(z) = I_n + N(z) Df(z) , \qquad z \in S ,
	\end{equation}
	where $I_n$ denotes the $n \times n$ identity matrix. Since $S$ is $k-$dimensional, $DH(z,0)|_S$ has $k$ multipliers equal to $1$ whose corresponding eigenvectors span the tangent space $T_zS$. The remaining $n-k$ multipliers $\mu_j(z), \ j = 1, \ldots , n-k,$ may or may not lie on the unit circle. In the following we refer to these $n-k$ multipliers as the \textit{non-trivial multipliers} 
	of $DH(z,0)|_S$.
	
	An important notion in the theory of fast-slow ODEs is the that of \textit{normal hyperbolicity}, which refers to the situation in which the non-trivial eigenvalues of the linearized layer problem 
	are all bounded away from the imaginary axis. It is straightforward to define an analogous notion of normal hyperbolicity for fast-slow maps \eqref{eq:gen_maps} in terms of the non-trivial multipliers $\mu_j(z)$.
	
	\begin{definition}
		\label{def:nh}
		\textup{(Normal hyperbolicity)} A point $z \in S$ is normally hyperbolic if the non-trivial multipliers $\mu_j(z)$ do not lie on the unit circle, i.e~
		\[
		|\mu_j(z)| \neq 1, \qquad  j=1,\ldots,n-k .
		\]
		Additionally, a normally hyperbolic point $p$ is called:
		\begin{enumerate}
			\item Attracting if $|\mu_j(z)| < 1$ for all $j=1,\ldots,n-k$;
			\item Repelling if $|\mu_j(z)| > 1$ for all $j=1,\ldots,n-k$;
			\item Saddle-type if $|\mu_j(z)| < 1$ for $n_a < n-k$ multipliers $\mu_j$, and $|\mu_j(z)| > 1$ for $n_r < n-k$ multipliers $\mu_j$, where $n_a + n_r = n-k$ and $n_a, n_r \geq 1$.
		\end{enumerate}
		These definitions are extended to sets, i.e.~a subset $S_n \subseteq S$ is called normally hyperbolic if every point in $S_n$ is normally hyperbolic, and a normally hyperbolic set $S_n$ is called attracting, repelling or saddle-type if every point in $S_n$ is attracting, repelling or saddle-type respectively.
	\end{definition}
	
	\begin{figure}[t!]
		\centering
		\subfigure[Attracting, $0 < \mu < 1$.]{\includegraphics[width=.4\textwidth]{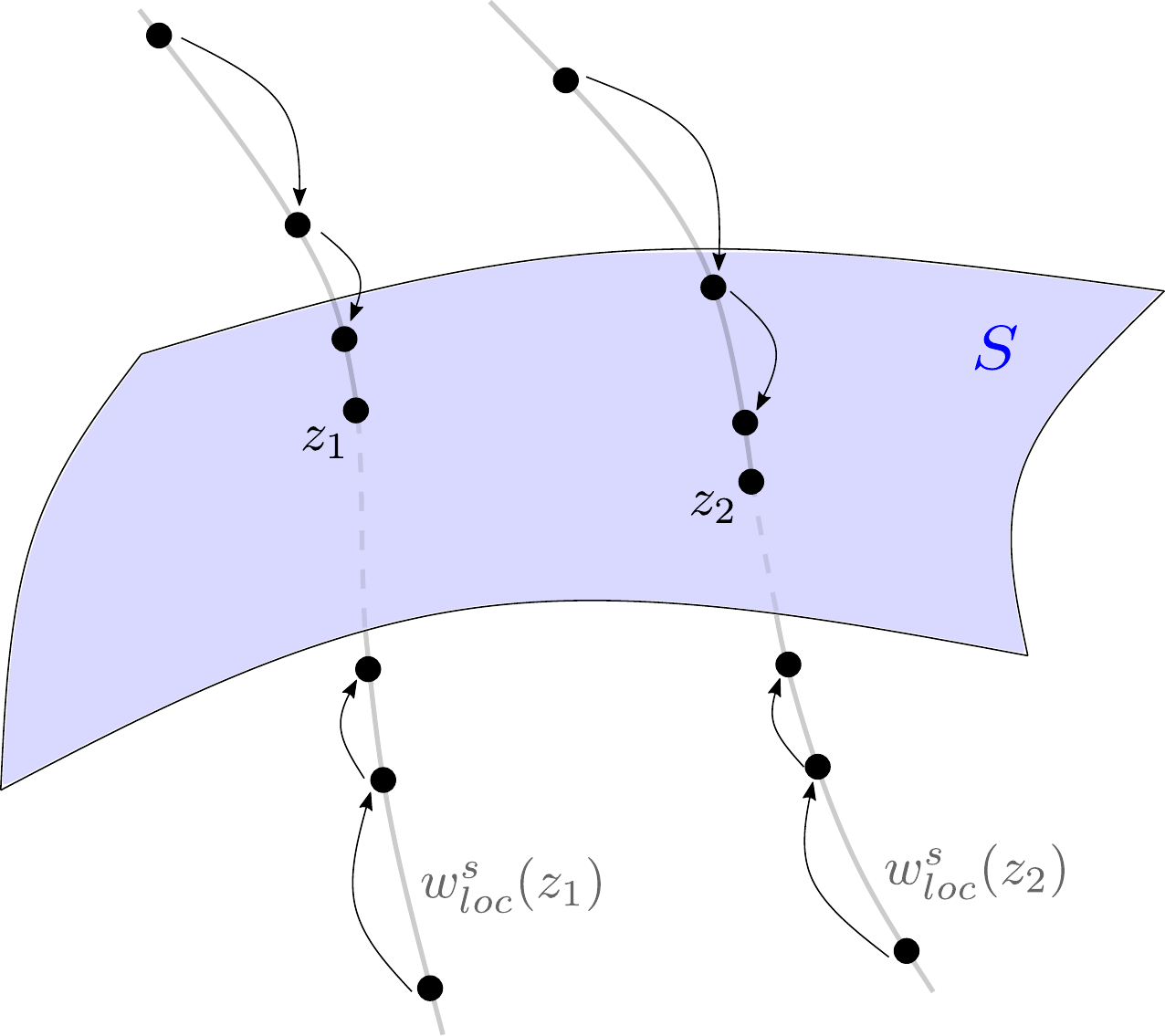}}
		\qquad
		\subfigure[Attracting, $-1 < \mu < 0$.]{\includegraphics[width=.4\textwidth]{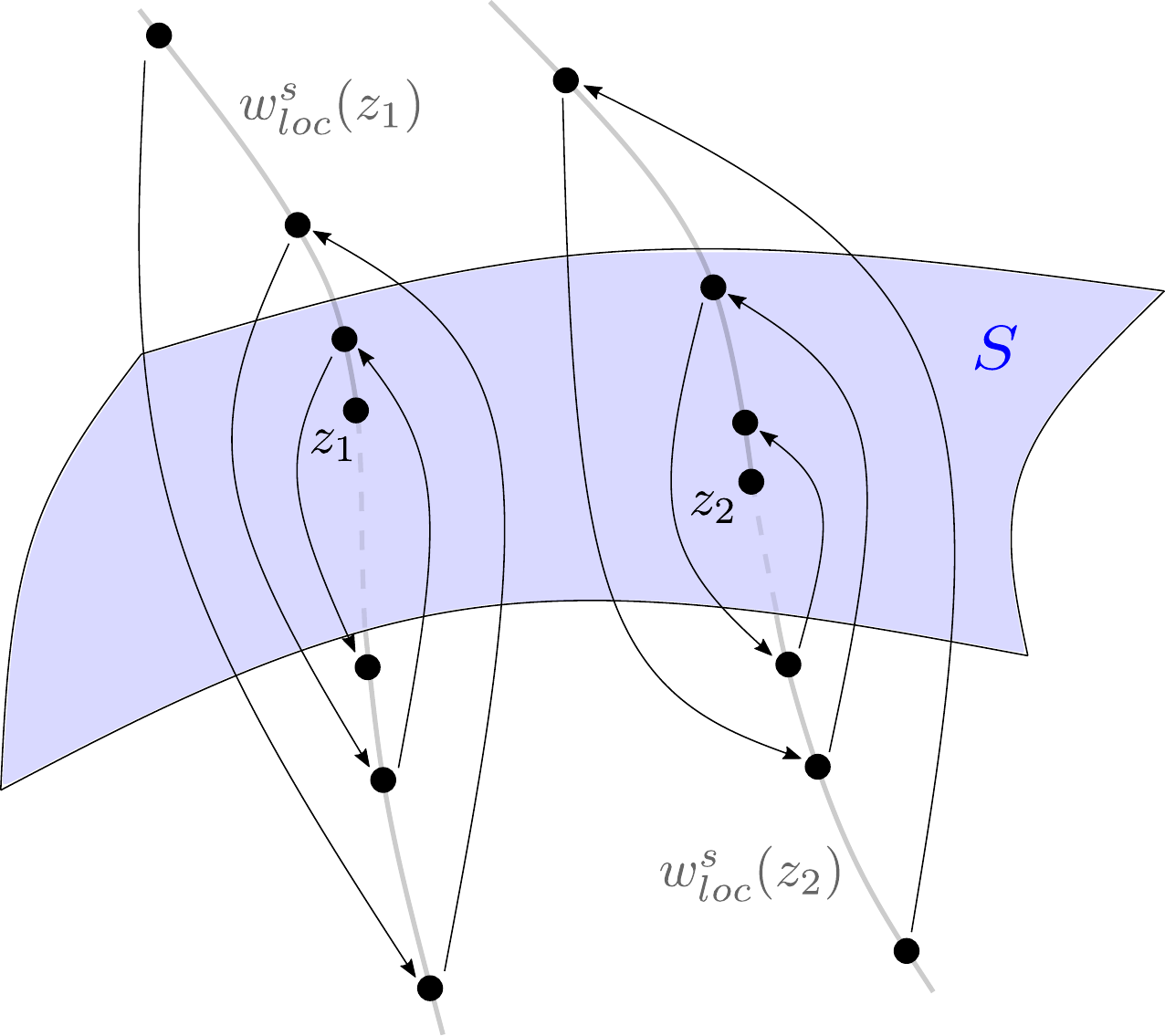}}
		\\
		\subfigure[Repelling, $ \mu > 1$.]{\includegraphics[width=.4\textwidth]{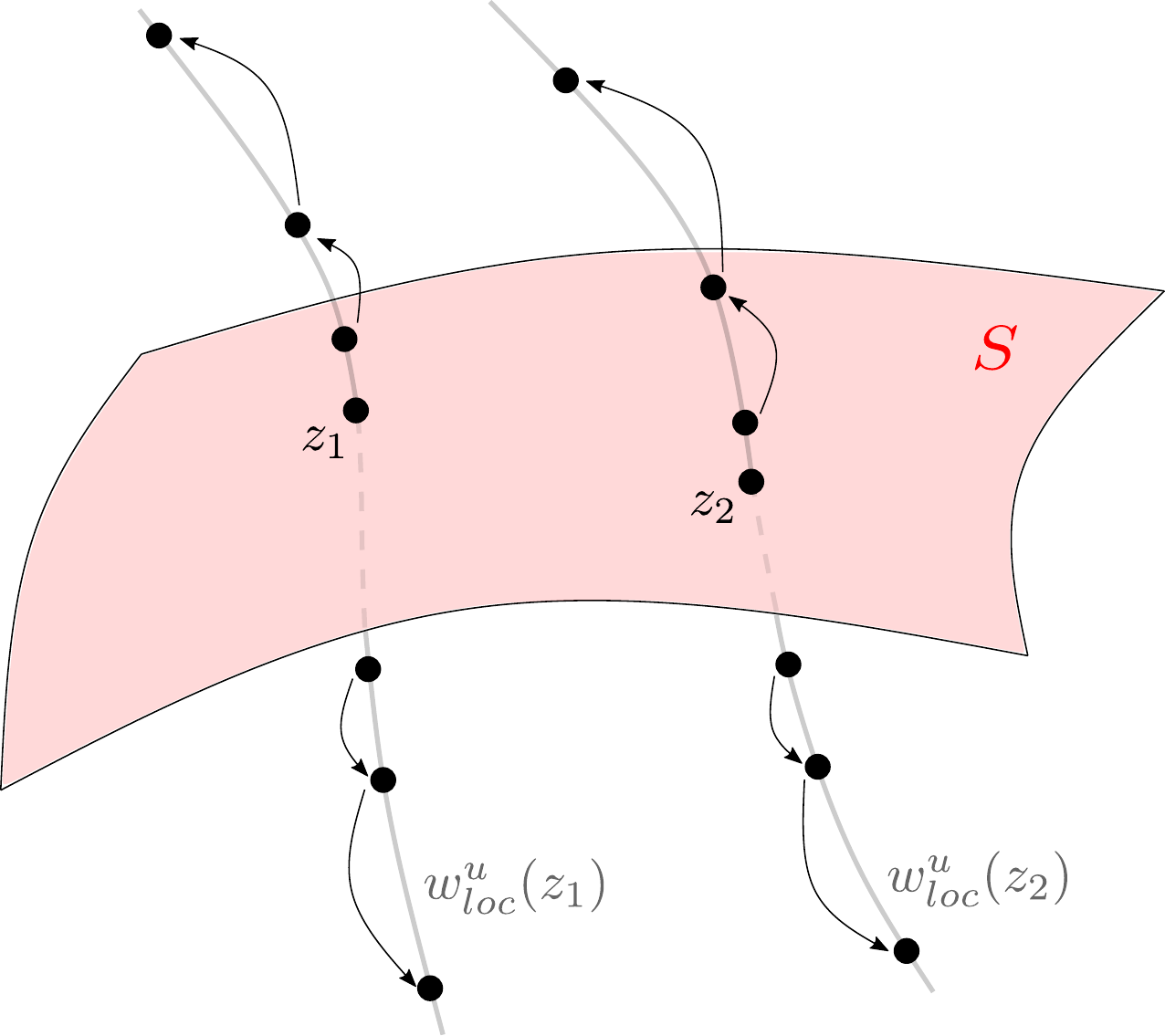}}
		\qquad
		\subfigure[Repelling, $ \mu < -1$.]{\includegraphics[width=.4\textwidth]{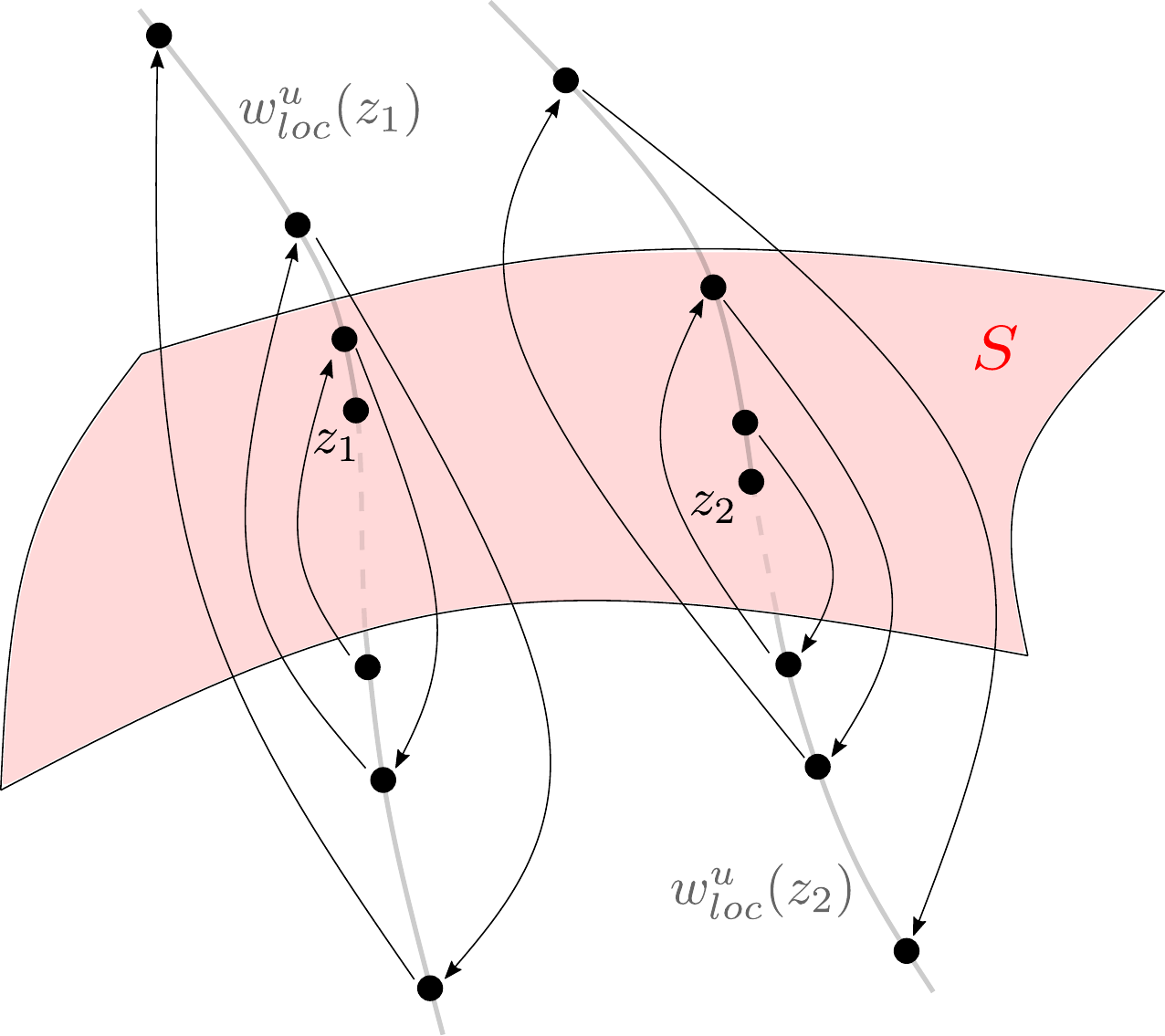}}
		\qquad
		\\
		\subfigure[Saddle-type, $0 < |\mu_a| < 1$, $| \mu_r | > 1$.]{\includegraphics[width=.45\textwidth]{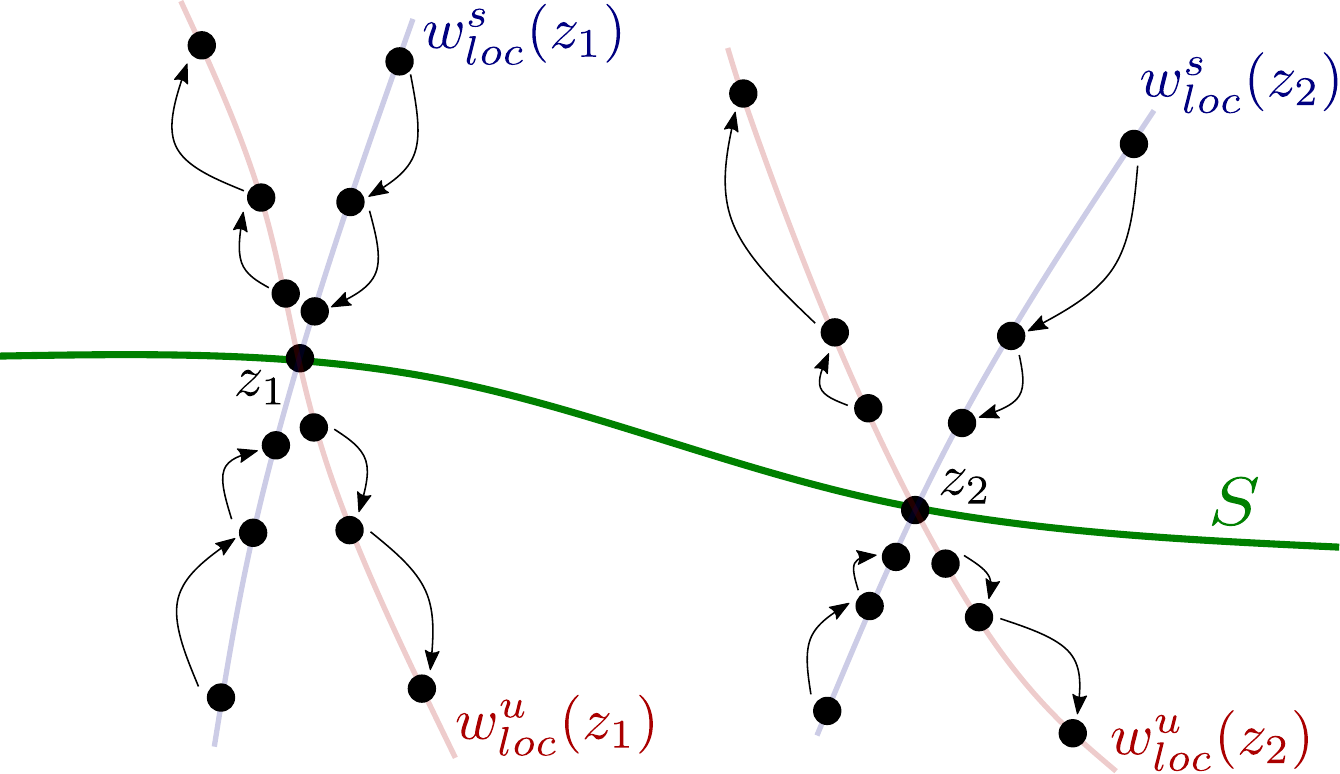}}
		\quad
		\subfigure[Saddle-type, $0 < | \mu_a | < 1$, $| \mu_r | > 1$.]{\includegraphics[width=.45\textwidth]{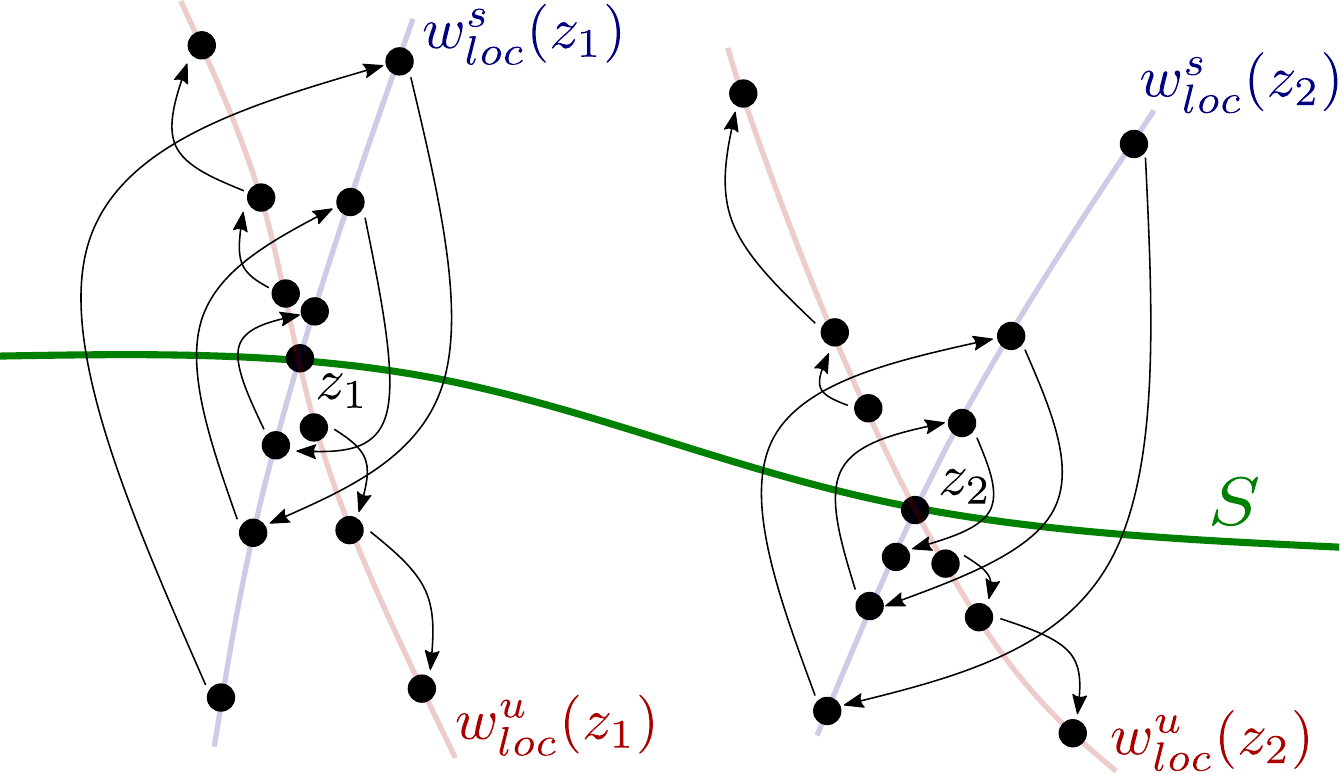}}
		\caption{Examples of normally hyperbolic critical manifolds $S$ in dimension $n=3$. In (a)-(d) $\dim S = 2$ and there exists a single non-trivial multiplier $|\mu| \neq 1$. There are four generic possibilities depending on $\mu$: (a) attracting and orientation preserving; 
			(b) attracting and orientation reversing; 
			(c) repelling and orientation preserving; 
			(d) repelling and orientation reversing. 
			In (e)-(f) $S$ is saddle-type with $\dim S = 1$ and two non-trivial multipliers $|\mu_a| < 1$ and $|\mu_r| > 1$. Two of four generic possibilities are shown: (e) orientation preserving along stable/unstable fibers $w^{s/u}_{loc}(z)$; (f) orientation reversing along $w^s_{loc}(z)$.}
		\label{fig:normal_hyperbolicity}
	\end{figure}

	Figure \ref{fig:normal_hyperbolicity} shows a number of representative scenarios in which the critical manifold $S$ is normally hyperbolic and attracting, repelling or saddle-type. Note, that the notion of normal hyperbolicity provided above applies equally for orientation-preserving and orientation-reversing maps.

	\begin{remark}
		\label{rem:singularities}
		Although we do not consider it further in this work, it is worthy to note that there are three generic codimension-1 possibilities for a loss of normal hyperbolicity along $S$ under additional parameter variation:
		\begin{enumerate}
			\item Fold-type singularities: Crossing of a real non-trivial multiplier over $1$.
			\item Flip-type singularities: Crossing of a real non-trivial multiplier over $-1$.
			\item Neimark-Sacker-type singularities: Crossing of a pair of complex conjugate non-trivial multipliers through $S^1 \setminus \{\pm 1\}$.
		\end{enumerate}
		With regards to the loss of normal hyperbolicity in fast-slow maps, the analogy to the theory for fast-slow ODEs is not a direct one. For example, the loss of normal hyperbolicity at a flip-type singularity has no direct analogue in the corresponding ODE theory.
	\end{remark}
	
	%

	The following result provides a useful method for determining the non-trivial multipliers $\mu_j$ in practice, and does not depend on a special choice of coordinates.
	
	
	\begin{proposition}
		\label{prop:EVs}
		We have the following equivalences:
		\begin{enumerate}
			\item The subset of non-trivial eigenvalues of the matrix $N Df|_S$ coincides with the set of eigenvalues of the $(n-k) \times (n-k)$ matrix $Df N|_S$.
			\item The set of non-trivial multipliers $\mu_j$ of the Jacobian matrix \eqref{eq:Jacobian} coincides with the set of eigenvalues of the $(n-k) \times (n-k)$ matrix $I_{n-k} + DfN|_S$.
		\end{enumerate}
	\end{proposition}
	
	\begin{proof}
		The proof proceeds by analogy with the corresponding statement for ODEs in \cite[Lemma 3.3]{Wechselberger2019}, with minor adaptations.
		
		We work in a neighbourhood about a point $p \in S$. Assume without loss of generality that $z = (x,y)^\textnormal{T}$ where $x \in \mathbb R^k$ and $y \in \mathbb R^{n-k}$ are chosen such that locally, $D_yf$ is an $(n-k) \times (n-k)$ regular (invertible) matrix; this is possible since $Df|_S$ is regular by Assumption \ref{ass:1}. In $(x,y)$ coordinates, the map \eqref{eq:gen_maps} is given by
		\begin{equation}
			\label{eq:gen_map_split}
			\begin{pmatrix}
				x \\
				y
			\end{pmatrix}
			\mapsto 
			\begin{pmatrix}
				\bar x \\
				\bar y
			\end{pmatrix}
			=
			\begin{pmatrix}
				x \\
				y
			\end{pmatrix}
			+
			\begin{pmatrix}
				N^x(x,y) \\
				N^y(x,y)
			\end{pmatrix}
			f(x,y) + \eps
			\begin{pmatrix}
				G^x(x,y,\eps) \\
				G^y(x,y,\eps)
			\end{pmatrix} ,
		\end{equation}
		where $N^x(x,y)$, $N^y(x,y)$ are matrices of dimensions $k \times (n-k)$ respectively $(n-k) \times (n-k)$, and $G^x(x,y,\eps)$, $G^y(x,y,\eps)$ are column vectors of length $k$ respectively $n-k$.
		
		Applying the coordinate transformation
		\begin{equation}
			\label{eq:coord_trans_straightening}
			v = f(x,y)
		\end{equation}
		with local (smooth) inverse $y = K(x,v)$ guaranteed to exist by the implicit function theorem, locally rectifies the critical manifold $S$. Explicitly, we obtain the following after local expansion about $v=0$, $\eps = 0$:
		\[
		\begin{split}
			v \mapsto \bar v & = f(\bar x, \bar y) \\
			& = f\left(x + N^x(x,y) v + \eps G^x(x,y,\eps) , y + N^y(x,y) v + \eps G^y(x,y,\eps) \right) \\
			& = f\left(x,y\right) + Df N(x,y) v + \eps Df G(x,y,\eps) + O(v^2,v\eps,\eps^2)  \\
			& = v + Df N(x,K(x,v)) v + \eps Df G(x,K(x,v),\eps) +  O(v^2,v\eps,\eps^2) ,
		\end{split}
		\]
		and so the map becomes
		\begin{equation}
			\label{eq:straight_S_map}
			\begin{split}
				x & \mapsto \bar x = x + \tilde N^x(x,v) v + \eps \tilde G^x(x,v,\eps) , \\
				v & \mapsto \bar v = v + Df N(x,K(x,v)) v + \eps Df G(x,K(x,v),\eps) +  O(v^2,v\eps,\eps^2) ,
			\end{split}
		\end{equation} 
		where we denote $\tilde N^x(x,v) := N^x(x,K(x,v))$ and $\tilde G^x(x,v,\eps) := G(x,K(x,v),\eps)$. Notice that \eqref{eq:straight_S_map} has critical manifold $S = \{v = 0\}$. The Jacobian evaluated along $S$ when $\eps = 0$ is
		\[
		\begin{pmatrix}
			I_k & \tilde N^x(x,0) \\
			O_{n-k,k} & I_{n-k} + Df N |_S
		\end{pmatrix}
		,
		\]
		which is block-diagonal with $k$ (trivial) multipliers equal to $1$, and $n-k$ non-trivial multipliers $\mu_j$ given by the eigenvalues of $I_{n-k} + DfN|_S$.
		
		Since both the coordinate transformation \eqref{eq:coord_trans_straightening} and its inverse (which is obtained via the implicit function theorem) are $C^r-$smooth with $r \geq 1$, the eigenvalues (multipliers) along $S = \{v = f(x,y) = 0\}$ are invariant under the composed transformation, i.e.~they are the same for map \eqref{eq:gen_maps}. The result follows since the point $p \in S$ about which we applied our arguments was arbitrary. 
	\end{proof}
	
	By Proposition \ref{prop:EVs}, the problem of calculating the $n-k$ non-trivial multipliers of the $n \times n$ matrix $DH(z,0)|_S$ reduces to the problem of calculating the multipliers of the $(n-k) \times (n-k)$ matrix $I_{n-k} + Df N|_S$. Of course, this further reduces to the problem of calculating the eigenvalues of $Df N|_S$.
	
	\begin{remark}
		In the proof of Proposition \ref{prop:EVs} we utilised a special choice of coordinates, in which the critical manifold $S$ is locally rectified along the $x-$axes. However, Proposition \ref{prop:EVs} itself applies independently of the choice of coordinates. The map \eqref{eq:straight_S_map} will often appear in proofs in later sections, and we shall frequently make use of the existence of special coordinates in proofs, but efforts will be made to ensure that `final results' do not depend on a choice of coordinates.
	\end{remark}
	
	\begin{remark}
		Proposition \ref{prop:EVs} applies whether or not $S$ is normally hyperbolic. Hence assertions (1)-(2) can be also used for identifying and classifying a loss of normal hyperbolicity.
	\end{remark}
	
	As an immediate consequence of Proposition \ref{prop:EVs}, we have that for each non-trivial multiplier $\mu_j$, there exists an eigenvalue $\lambda_j$ of the matrix $Df N|_S$ such that
	\begin{equation}
		\label{eq:multipliers}
		\mu_j(z) = 1 + \lambda_j(z) , \qquad z \in S.
	\end{equation}
	This leads to an alternative characterisation of normal hyperbolicity in terms of the eigenvalues $\lambda_j$.
	
	\begin{corollary}
		\label{cor:nh}
		A point $z \in S$ is normally hyperbolic if and only if eigenvalues $\lambda_j$ of the matrix $DFN|_S$ satisfy
		\[
		2 \textup{Re}\, \lambda_j(z) \neq - |\lambda_j(z)|^2 , \qquad j = 1, \ldots , n-k.
		\]
		Moreover, a normally hyperbolic point $z$ is attracting if
		\[
		2 \textup{Re}\, \lambda_j(z) < - |\lambda_j(z)|^2
		\] 
		for all $j  = 1, \ldots , n-k$, repelling if
		\[
		2 \textup{Re}\, \lambda_j(z) > - |\lambda_j(z)|^2
		\]
		for all $j = 1, \ldots , n-k$, and saddle-type otherwise. 
	\end{corollary}
	
	\begin{proof}
		This follows directly from Proposition \ref{prop:EVs} and equation \eqref{eq:multipliers}.
	\end{proof}
	
	Finally, we note that on the linear level, a normally hyperbolic critical manifold $S_n \subseteq S$ induces a natural (pointwise) splitting
	\begin{equation}
		\label{eq:splitting_pointwise}
		T_z \mathbb R^n = T_z S_n \oplus \mathcal N_z , \qquad z \in S_n,
	\end{equation}
	with 
	$\mathcal N_z := \text{span} \{N^i(z)\}_{j = 1, \ldots , n-k} = E^s(z) \cup E^u(z)$, 
	where $N^i(z)$ denotes the $i$'th column of $N(z)$ and $E^{s/u}(z)$ denote stable/unstable eigenspaces at $z$ respectively. 
	The pointwise splitting \eqref{eq:splitting_pointwise} leads to the global splitting
	\begin{equation}
		\label{eq:splitting_bundles}
		\mathbb R^n \cong T \mathbb R^n = TS_n \oplus \mathcal N ,
	\end{equation}
	where $TS_n = \cup_{z \in S_n} T_zS$ and $\mathcal N = \cup_{z \in S_n} \mathcal N_z$ denote the tangent bundle and transverse linear fiber bundle associated to $S_n$, respectively.
	
	\begin{figure}[t!]
		\centering
		\includegraphics[scale=0.45]{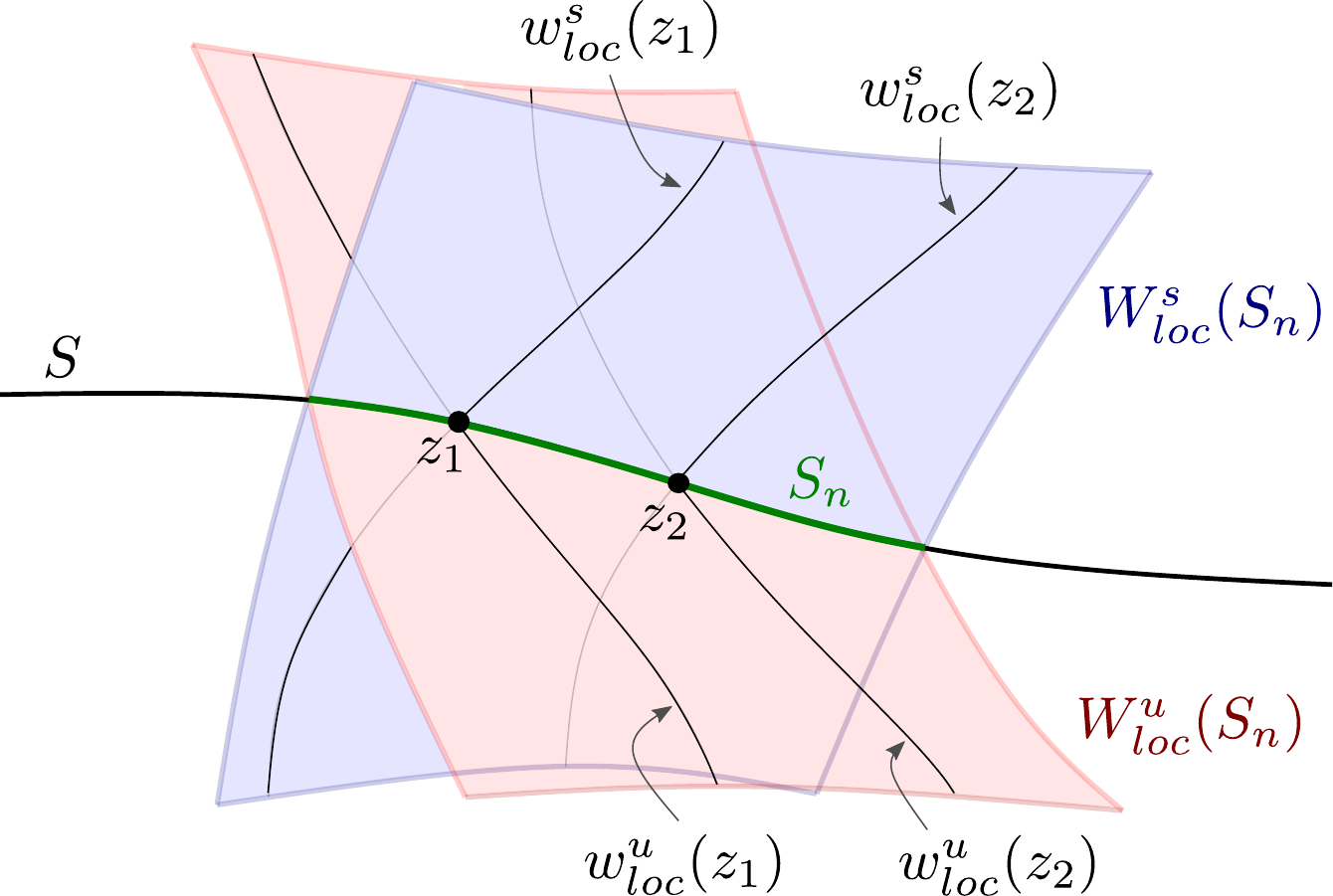}
		\caption{Stable/unstable manifolds $W^{s/u}_{loc}(S_n)$ associated to a normally hyperbolic critical manifold $S_n$ are foliated by lower dimensional stable/unstable manifolds $w^{s/u}_{loc}(z)$ associated to base points $z \in S_n$. The case of a $1-$dimensional saddle-type critical manifold in dimension $n=3$ is shown here ($S_n$ is shown in green). In this case the manifolds $W^{s/u}_{loc}(S_n)$ (shown in blue/red) are $2-$dimensional and foliated by the $1-$dimensional stable/unstable manifolds $w^{s/u}_{loc}(z)$ associated to base points $z \in S_n$. Representative stable/unstable fibers $w^{s/u}(z_i)$ with base points $z_i \in S_n$ for $i=1,2$ are also shown.}
		\label{fig:foliation}
	\end{figure}
	
	On the nonlinear level, a normally hyperbolic critical manifold $S_n \subseteq S$ induces a local foliation of the adjacent space by its local stable and unstable manifolds $w_{loc}^s(z)$ and $w_{loc}^u(z)$. We denote induced stable and unstable foliations by
	\begin{equation}
		\label{eq:layer_foliation}
		W_{loc}^s(S_n) := \bigcup_{z \in S_n} w_{loc}^s(z) , \qquad 
		W_{loc}^u(S_n) := \bigcup_{z \in S_n} w_{loc}^u(z) ,
	\end{equation}
	respectively. Note that if $|\mu_j| < 1$ for $n_a$ non-trivial multipliers and $|\mu_j| > 1$ for $n_r$ non-trivial multipliers ($n_a + n_r = n-k$), then $W_{loc}^s(S_n)$ is $(k + n_a)-$dimensional and $W_{loc}^u(S_n)$ is $(k + n_r)-$dimensional. By definition, the fibers $w_{loc}^{s/u}(z)$ and hence the induced foliations in \eqref{eq:layer_foliation} are locally invariant under the layer map \eqref{eq:layer_map}. Together, $W^{s/u}_{loc}(S_n)$ form a local nonlinear fiber bundle, whose linear part coincides with the transverse linear fiber bundle $\mathcal N$. The case of a saddle-type normally hyperbolic critical manifold in $\mathbb R^3$ with $n_s = n_u = 1$ is sketched in Figure \ref{fig:foliation}.
	
	\begin{remark}
		In the case that $S$ is normally hyperbolic, the transverse linear fiber bundle $\mathcal N$ is related to the normal bundle associated to the tangent bundle $TS$ by the smooth coordinate transformation described in Remark \ref{rem:Fenichel_normal_form} below. Thus it is common to find the terms ``transverse bundle" and ``normal bundle" used interchangeably in the literature. The former terminology is preferred in this work in order to emphasise the coordinate-independence of the formalism.
	\end{remark}

	\subsection{The reduced and $m$'th iterate maps}
	\label{sub:reduced_map}
	
	In the context of fast-slow ODEs one obtains a second, non-equivalent limiting problem -- the \textit{reduced problem} -- by considering the singular limit $\eps \to 0$ taken with respect to the system of fast-slow ODEs posed on the so-called \textit{slow time} $\tau = \eps t$. This reduced problem induces a flow on normally hyperbolic submanifolds of $S$ on the slow time-scale $\tau$. By analogy, 
	one might expect a similar equivalence between maps of the form \eqref{eq:gen_maps} and `slow' maps of the form
	\[
	z \mapsto z + \frac{1}{\eps} N(z) f(z) + G(z,\eps) .
	\]
	Unfortunately, such an equivalence is not available for maps, since there is in general no analogue of the time rescaling which would render the maps
	\[
	\bar z - z \qquad \text{and} \qquad \eps^{-1}(\bar z - z) 
	\]
	equivalent for each $0 < \eps \ll 1$. It is natural to ask, then, whether it is still possible to derive a `reduced map' which describes the limiting dynamics on (normally hyperbolic submanifolds of) the critical manifold $S$, thereby providing information on the limiting dynamics which is not present in the layer map \eqref{eq:layer_map}.
	
	In this section we show that such a map can be derived, though it must be formulated as a map governing the leading order dynamics on a locally invariant $C^r-$smooth \textit{slow manifold} $S_\eps$ which converges to $S$ as $\eps \to 0$. Under certain nondegeneracy conditions, sequential iterates of this map on the slow manifold $S_\eps$ are $O(\eps)-$close to each other. 
	We also derive the reduced $m$'th iterate map on $S_\eps$ induced by repeated iteration of \eqref{eq:gen_maps}. Interestingly, this map can be related to a discretization of the ODE reduced problem associated to fast-slow ODEs in \cite{Wechselberger2019} if the number of iterates $m$ is comparable to $\eps^{-1}$. 
	
	\begin{remark}
		As stated above, in the following we need to assume the existence of an invariant slow manifold in order to derive both the reduced and $m$'th iterate maps. The existence of slow manifolds under suitable (normally hyperbolic) conditions will be treated in Section \ref{sec:slow_manifold_theorems} (see Theorems \ref{thm:slow_manifolds} and \ref{thm:graph_slow_manifolds}), and does not depend on the existence of a reduced map.
	\end{remark}


	\subsubsection{The reduced map}
	
\begin{figure}[t!]
		\centering
		\includegraphics[scale=0.3]{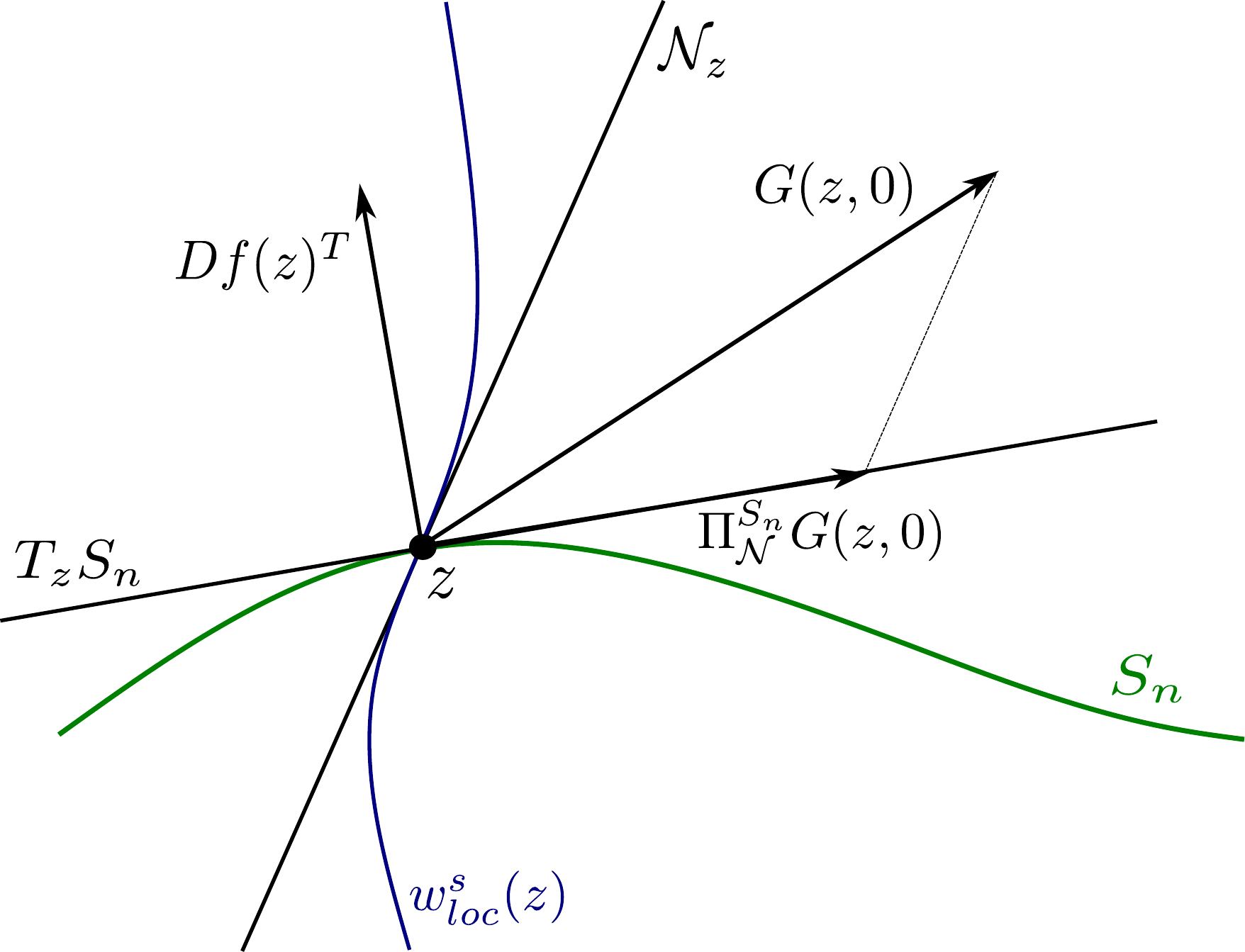}
		\caption{Action of the oblique projection operator $\Pi^{S_n}_{\mathcal N}$ along $\mathcal N$ onto $TS_n$ in the case of a $1-$dimensional attracting critical manifold $S_n$ (green). The splitting \eqref{eq:splitting_bundles} implies that the leading order perturbation $G(z,0)$ at each point $z \in S_n$ can be decomposed into components in the tangent bundle $T_z S_n$, and the (transverse) linear fast fiber $\mathcal N_z$ tangent to the nonlinear fast fiber $w^s_{loc}(z)$ shown in blue in the figure. The component in $T_zS_n$ is given by $\Pi_{\mathcal N}^{S_n}G(z,0) \in T_zS_n$, and sufficient to determine the reduced map \eqref{eq:red_map}, see also Proposition \ref{prop:reduced_map}.}
		\label{fig:projection}
	\end{figure}

	Consider again a connected, normally hyperbolic submanifold $S_n \subseteq S$. Due to the splitting \eqref{eq:splitting_bundles}, there exists a unique projection operator
	\begin{equation}
		\label{eq:proj_def}
		\Pi^{S_n}_{\mathcal N} : T\mathbb R^n = TS_n \oplus \mathcal N \to TS_n ,
	\end{equation}
	which projects vectors in $T \mathbb R^n$ onto their component in $TS_n$ along the direction of $\mathcal N$, see Figure \ref{fig:projection}. This allows one to isolate components of the map contributing to dynamics in $TS_n$. The following result should be compared with the characterisation of the reduced problem in fast-slow ODEs presented in \cite[Def.~3.8]{Wechselberger2019}, see also \cite{Fenichel1979,Goeke2014,Kruff2019,Lizarraga2020}.
	
	\begin{proposition}
		\label{prop:reduced_map}
		Fix $\eps \in [0,\eps_0)$ and denote by $S_{n,\eps}$ a locally invariant slow manifold perturbing from $S_n$, as described by Theorem \ref{thm:slow_manifolds}. Then for $\eps_0$ sufficiently small we have
		\begin{equation}
			\label{eq:red_map_expanded}
			\bar z |_{S_{n,\eps}} = z + \eps \Pi^{S_n}_{\mathcal N} G(z,0) + O(\eps^2) , \qquad z \in S_n .
		\end{equation}
		Moreover, the unique projection operator $\Pi^{S_n}_{\mathcal N}$ has matrix representation
		\begin{equation}
			\label{eq:proj}
			\Pi^{S_n}_{\mathcal N} = I_n - N (Df N)^{-1} Df |_{S_n} .
		\end{equation}
	\end{proposition}
	
	\begin{proof}
		Without loss of generality, we make the same preliminary assumptions as in the proof of Proposition \ref{prop:EVs} and consider again the map \eqref{eq:gen_map_split} in a neighbourhood about some $p \in S_n \subseteq S$. In particular, we assume that the matrix $D_yf$ is locally regular so that the critical manifold $S$ has a graph representation $y = \varphi_0(x)$.
		
		To obtain the map \eqref{eq:red_map_expanded}, we first obtain a local graph form for the invariant slow manifold by substituting the ansatz
		\[
		S_{n,\eps} : y = \varphi_\eps(x) = \varphi_0(x) + \varphi_1(x) \eps + O(\eps^2) ,
		\]
		where $f(x,\varphi_0(x)) = 0$, into the invariance equation
		\[
		\bar y = \varphi_\eps(\bar x) = \varphi_\eps(x) + N^y (x,\varphi_\eps(x)) f(x,\varphi_\eps(x)) + \eps G^y (x,\varphi_\eps(x),\eps) .
		\]
		By redefining $N(x,y)$ if necessary, we may assume that
		\[
		f(x,y) = y - \varphi_0(x) ,
		\]
		without loss of generality. By expanding in $\eps$ and requiring that $(x,y) \in S_\eps \implies (\bar x, \bar y) \in S_\eps$, we obtain the following invariance equation after matching terms at $O(\eps)$:
		\[
		\begin{split}
			\varphi_1(x) &+ N^y(x,\varphi_0(x)) \varphi_1(x) + G^y(x,\varphi_0(x),0)  \\
			&= \varphi_1(x) + D_x\varphi_0(x) N^x(x,\varphi_0(x)) \varphi_1(x) + D_x\varphi_0(x) G^x(x,\varphi_0(x),0) .
		\end{split}
		\]
		Rearranging for $\varphi_1(x)$ and using
		\[
		\begin{split}
			DfN|_S &= - D_x \varphi_0(x) N^x(x,\varphi_0(x)) + N^y(x,\varphi_0(x)) , \\
			DfG|_S &= - D_x \varphi_0(x) G^x(x,\varphi_0(x),0) + G^y(x,\varphi_0(x),0) ,
		\end{split}
		\]
		yields
		\begin{equation}
			\label{eq:slow_manifold}
			S_\eps : y = \varphi_\eps(x) = \varphi_0(x) - \eps (Df N)^{-1} Df G |_S + O(\eps^2) .
		\end{equation}
		Substituting \eqref{eq:slow_manifold} into \eqref{eq:gen_map_split} and expanding again in $\eps$ yields \eqref{eq:red_map_expanded} with $\Pi^{S_n}_{\mathcal N}$ given by the matrix \eqref{eq:proj}.
	\end{proof}
	
	\begin{remark}
		\label{rem:slow_manifold_expansion}
		The expression \eqref{eq:slow_manifold} obtained in the proof above is of interest its own right, as it provides a general formula for the slow manifold parameterization up to $O(\eps)$ in terms of $N$, $f$ and $G$ in the case that the critical manifold $S$ is given as a graph $y = \varphi_0(x)$.
	\end{remark}
	
	Proposition \ref{prop:reduced_map} allows for the following definition.
	
	\begin{definition}
		\label{def:red_map}
		\textup{(Reduced map)} The map
		\begin{equation}
			\label{eq:red_map}
			z \mapsto \bar z = z + \eps \Pi^{S_n}_{\mathcal N} G(z,0) , \qquad z \in S_n ,
		\end{equation}
		obtained by truncating \eqref{eq:red_map} at $O(\eps^2)$ is called the reduced map.
	\end{definition}
	
	
	Similarly to the theory for fast-slow ODEs, the reduced map describes the leading order dynamics on locally invariant slow manifolds perturbing from normally hyperbolic submanifolds of $S$. In contrast to the ODE theory, an $\eps-$independent expression for the right-hand-side is not available, since there is in general no way to `divide out' the factor of $\eps$ while preserving topological conjugacy. Fortunately, this is no hindrance in applications, since what is really needed is a calculable asymptotic approximation for the flow along locally invariant slow manifolds with $0 < \eps \ll 1$. This is provided by the formulation of the reduced map \eqref{eq:red_map}.
	
	\begin{remark}
		A similar formulation of the reduced problem also exists in the ODE setting, i.e.~there too, the reduced vector field is precisely the leading order vector field for the dynamics on Fenichel slow manifolds, occurring at $O(\eps)$ on the fast time-scale.
	\end{remark}
	
	
	\begin{remark}
		\label{rem:reduced_map_discretizations}
		In the special case of maps induced via discretization of a fast-slow ODE, a reduced map having an $\eps-$independent right-hand-side can typically be obtained by a rescaling of the discretization/step parameter $h>0$. For example, Euler-discretization of a general fast-slow ODE
		\[
		z' = N(z) f(z) + \eps G(z,\eps) , \qquad z \in \mathbb R^n, \qquad 0 < \eps \ll 1,
		\]
		leads to the map
		\[
		\frac{\bar z - z}{h} = N(z) f(z) + \eps G(z,\eps) .
		\]
		Rescaling $h = \eps^{-1} \tilde h$ leads to the reduced map
		\[
		\frac{\bar z - z}{\tilde h} = \Pi^{S_n}_{\mathcal N} G(z,0) , \qquad z \in S_n ,
		\]
		which has no $\eps-$dependence in the right-hand-side. Euler discretizations of this kind are considered further in Section \ref{sub:Euler_discretization}.
	\end{remark}
	
	Finally we note that similarly to the ODE case, the projection operator $\Pi^{S_n}_{\mathcal N}$ is not defined at fold-type singularities since the quantity $\Pi^{S_n}_{\mathcal N}G(z,0)$ blows up. It is however defined at flip and Neimark-Sacker-type singularities. We do not consider these issues (which relate to the loss of normal hyperbolicity) in detail in this work. See however Section \ref{sec:outlook} for further discussion.

	\subsubsection{The $m$'th iterate map}
	
	In a neighbourhood of the critical manifold $S$ it is also possible to derive the form of the $m$'th iterate map induced by repeated iteration of the map \eqref{eq:gen_maps} using local asymptotic self-similarity properties of the map. 
	As before, $S_n$ denotes a compact normally hyperbolic submanifold of $S$.
	
	
	\begin{proposition}
		\label{prop:red_jump_map}
		For all $\eps \in [0,\eps_0)$ with $\eps_0$ sufficiently small, the $m$'th iterate map induced by \eqref{eq:gen_maps} in a tubular neighbourhood of $S$ takes the form of a non-autonomous map
		\begin{equation}
			\label{eq:slow_formulation}
			z \mapsto \bar z^m = z + m \left(N(z)f(z) + \eps G(z,\eps) \right) + R(z,\eps,m) ,
		\end{equation}
		where for all fixed $m \in \mathbb N_+$ the higher order term $R(z,\eps,m)$ satisfies
		\begin{enumerate}
			\item[(i)] $R(z,\eps,m) = O(\eps^2,\eps f(z),f(z)^2)$;
			\item[(ii)] $R(z,\eps,1) \equiv O_n$;
			\item[(iii)] $R(z,0,m)|_{S_n} \equiv O_n$;
			\item[(iv)] $R(z,\eps,m) \to O_n$ as $(\eps,m) \to (0,\infty)$.
		\end{enumerate}
		In particular, the reduced $m$'th iterate map on $S_{n,\eps}$ is given by
		\begin{equation}
			\label{eq:mth_map}
			z \mapsto \bar z^m = z + \eps m \Pi^{S_n}_{\mathcal N} G(z,0) , \qquad z \in S_n ,
		\end{equation}
		where $\Pi^{S_n}_{\mathcal N}$ is the projection operator \eqref{eq:proj}.
	\end{proposition}
	
	\begin{proof}
		We work in a tubular neighbourhood of $S$ within which Taylor expansion about $f(z) = O_{n-k}$ is valid. Repeated iteration of the map \eqref{eq:gen_maps} leads to the following asymptotically self-similar sequence:
		\[
		\begin{split}
			\bar z &= z + N(z)f(z) + \eps G(z,\eps) , \\
			\bar z^2 &= z + 2 \left( N(z)f(z) + \eps G(z,\eps) \right) + R(z,\eps,2) , \\
			& \ \vdots \\
			\bar z^m &= z + m \left( N(z)f(z) + \eps G(z,\eps) \right) + R(z,\eps,m) ,
		\end{split}
		\]
		where $m\in\mathbb N_+$. For fixed $m$, the properties (i)-(iii) of the remainder term $R(z,\eps,m)$ are shown by induction on $m$, and omitted here for brevity. In order to prove property (iv), we adopt a sequential, componentwise notation
		\[
		R(z,\eps,m) := R_m(z,\eps) = (R_m^1(z,\eps), \ldots , R_m^n(z,\eps))^\textnormal{T} \in \mathbb R^n , \qquad m \in \mathbb N_+ ,
		\]
		and consider the real-valued sequences $(R^i_m)_{m \in \mathbb N_+}$ defined by the component functions $R_m^i : \mathbb R^n \times [0,\eps_0) \to \mathbb R$, $i = 1, \ldots , n$. Since by property (iii) (which holds for all fixed $m \in \mathbb N_+$) we have that
		\[
		\lim_{m \to \infty} R^i_m(z,0) = 0 , \qquad z \in S_n ,
		\]
		uniformly in $(z,\eps)$, and by continuity of $R$ we have that
		\[
		\lim_{\eps \to 0} R^i_m(z,\eps) = R^i_m(z,0) = 0 , \qquad z \in S_n ,
		\]
		pointwise in $(z,\eps)$, it follows by an application of the Moore-Osgood theorem \cite{Osgood1912} that the double limit is defined and commutes such that
		\[
		\lim_{\eps \to 0} \lim_{m \to \infty} R^i_m(z,\eps) = 
		\lim_{m \to \infty} \lim_{\eps \to 0} R^i_m(z,\eps) = 0 , \qquad z \in S_n .
		\]
		Since the above holds for each $i = 1, \ldots , n$, it follows that $R(z,\eps,m)|_{S_n} \to 0$ as $(\eps,m) \to (0,\infty)$, as required. 
		
		\
		
		It remains to show that the reduced $m$'th iterate map on $S_{n,\eps}$ is given by \eqref{eq:mth_map}. Since the expression \eqref{eq:mth_map} is pointwise, it suffices to show it in the local $z = (x,y)^\textnormal{T}$ coordinates of the map \eqref{eq:gen_map_split} used in proof of Propositions \ref{prop:EVs} and \ref{prop:reduced_map}. Specifically, the same preliminary assumptions used in the proof of Propositions \ref{prop:EVs} and \ref{prop:reduced_map} lead to the following local formulation of the $m$'th iterate map \eqref{eq:slow_formulation} near $S$:
		\[
		\begin{split}
			x &\mapsto \bar x^m = x + m \left( N^x(x,y) (y - \varphi_0(x)) + \eps G^x(x,y,\eps) \right) + R^x(x,y,\eps,m) ,  \\
			y &\mapsto \bar y^m = y + m \left( N^y(x,y) (y - \varphi_0(x)) + \eps G^y(x,y,\eps) \right) + R^y(x,y,\eps,m) ,
		\end{split}
		\]
		where we assume that $f(x,y) = y - \varphi_0(x)$ as in the proof of Proposition \ref{prop:reduced_map} (recall that this can be achieved by redefining $N$ if necessary), and where the remainder terms satisfy
		\begin{equation}
			\label{eq:rem_mth_map}
			\begin{split}
				R^x(x,y,\eps,m) &= O \left(\eps^2, \eps (y - \varphi_0(x)), (y - \varphi_0(x))^2 \right) , \\
				R^y(x,y,\eps,m) &= O \left(\eps^2, \eps (y - \varphi_0(x)), (y - \varphi_0(x))^2 \right) ,
			\end{split}
		\end{equation}
		for each fixed $m \in \mathbb N_+$ due to the remainder property (i). By invariance, we can simply restrict the map to $S_{n,\eps}$, which in the chosen coordinates has graph representation \eqref{eq:slow_manifold}. We obtain
		\begin{equation}
			\label{eq:mth_map_3}
			\begin{split}
				\begin{pmatrix}
					\bar x^m \\
					\bar y^m
				\end{pmatrix}
				&=
				\begin{pmatrix}
					x \\
					y
				\end{pmatrix}
				+
				\eps m
				\begin{pmatrix}
					- N^x (Df N)^{-1} Df G + G^x \\
					- N^y (Df N)^{-1} Df G + G^y
				\end{pmatrix}
				\bigg|_{y = \varphi_0(x), \eps = 0}
				+ O(\eps^2) \\
				&= 
				\begin{pmatrix}
					x \\
					y
				\end{pmatrix}
				+ \eps m \left( I_n - N (Df N)^{-1} Df \right) G \big|_{y = \varphi_0(x), \eps = 0} + O(\eps^2) \\
				&= 
				\begin{pmatrix}
					x \\
					\varphi_0(x)
				\end{pmatrix}
				+ \eps m \Pi_{\mathcal N}^{S_n} G(x, \varphi_0(x) ,0) + O(\eps^2) ,
			\end{split}
		\end{equation}
		after expanding in $\eps$, using the fact that $R^x(x,\varphi_\eps(x),\eps,m), R^y(x,\varphi_\eps(x),\eps,m) = O(\eps^2)$ for each fixed $m \in \mathbb N_+$ by equation \eqref{eq:rem_mth_map} above, and substituting \eqref{eq:proj} for $\Pi_{\mathcal N}^{S_n}$ in the final equality. The expression in \eqref{eq:mth_map} follows from \eqref{eq:mth_map_3} after truncating the $O(\eps^2)$ terms in line with Definition \ref{def:red_map}.
		%
	\end{proof}

	\begin{figure}[t!]
		\centering
		\subfigure[Reduced map on $S_{n,\eps}$.]{\includegraphics[width=.6\textwidth]{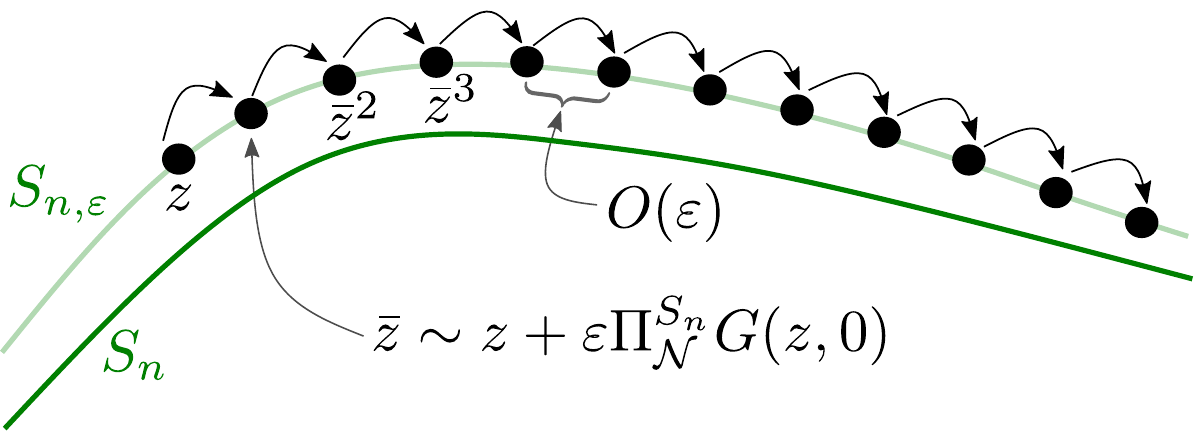}}
		\\
		\subfigure[$m$'th iterate map on $S_{n,\eps}$.]{\includegraphics[width=.64\textwidth]{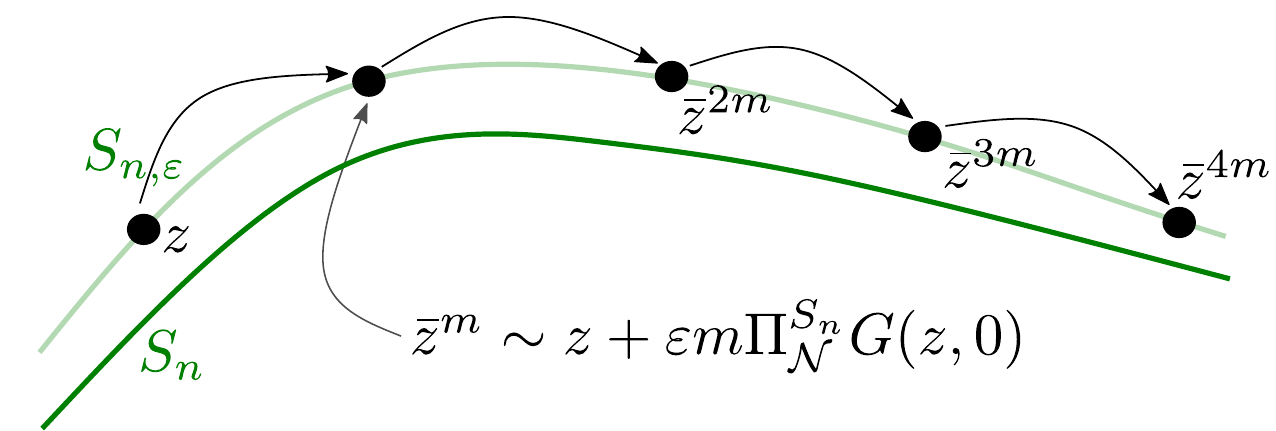}}
		\caption{Reduced and $m$'th iterate maps on $S_{n,\eps}$ defined in equations \eqref{eq:red_map} and \eqref{eq:mth_map} respectively. The normally hyperbolic critical manifold and the corresponding perturbed slow manifold are shown again in green and shaded green respectively. The dynamics of the map \eqref{eq:gen_maps} within $S_{n,\eps}$ are governed to leading order in $\eps$ by the reduced map \eqref{eq:red_map}. Iterates are generically separated by distances of $O(\eps)$, as shown in (a). Iterates of the $m$'th iterate reduced map \eqref{eq:mth_map} shown in (b) can be separated by $O(1)$ distances if the number of iterates $m$ is comparable to $\eps^{-1}$. For both reduced and $m$'th iterate maps, non-trivial dynamics is only possible for $0 < \eps \ll 1$. See however Remark \ref{rem:mth_map}, which describes a possible approach to identifying non-trivial dynamics on $S$ for $\eps = 0$ in the dual limit $(\eps,m) \to (0,\infty)$.}
		\label{fig:reduced_maps}
	\end{figure}
	
	See Figure \ref{fig:reduced_maps} for a comparison of reduced and $m$'th iterate maps. If the number of iterates $m \in \mathbb N_+$ is comparable to $\eps^{-1}$, then iterates of the $m$'th iterate map are generically separated by $O(m \eps) = O(1)-$distances. 
	In essence, the $\eps-$dependence cancels if one `speeds up' the map by considering only every $m$'th iterate with $m = \eps^{-1} \tilde m$ for some $\tilde m > 0$. 
	This has a similar effect to the time rescaling $t = \eps^{-1} \tau$ in fast-slow ODEs. In fact the $m$'th iterate map \eqref{eq:mth_map} relates directly to the ODE reduced problem $z' = \Pi_{\mathcal N}^{S_n} G(z,0)$ associated to general fast-slow systems \eqref{eq:nonstnd_form} (see \cite[Def.~3.8]{Wechselberger2019}), via the Euler discretization of the latter. Specifically, the map \eqref{eq:mth_map} can be rewritten as
	\begin{equation}
		\label{eq:reduced_Euler}
		\frac{{\bar z}^m - z}{\tilde m} = \Pi_{\mathcal N}^{S_n} G(z,0) ,
		\qquad z \in S_n ,
	\end{equation}
	which coincides with an Euler discretization of the ODE reduced problem if the step size $\tilde m > 0$ satisfying $\eps^{-1} \tilde m = m \in \mathbb N_+$; recall also Remark \ref{rem:reduced_map_discretizations}.
	
	\begin{remark}
		\label{rem:mth_map}
		The preceding discussion pertains to a relationship between the Euler discretization \eqref{eq:reduced_Euler} and the $m$'th iterate reduced map \eqref{eq:mth_map}, which is obtained by truncating the map \eqref{eq:mth_map_3} at $O(\eps^2)$. In order to rigorously prove a relationship between the discretized (continuous-time) reduced problem and an $m$'th iterate map on $S$ for $\eps = 0$, one must define a suitable embedding of the map \eqref{eq:mth_map_3} with $m \in \mathbb N_+$ into the parameter-dependent family
		\[
		\bar z_{\mathbb R} = z + \eps \zeta \Pi_{\mathcal N}^{S_n} G(z,0) + R(z,\eps,\zeta) , \qquad z \in S_n , \qquad \zeta > 0,
		\]
		and study the limit $\eps \to 0$ in the case that $\zeta = \eps^{-1} \tilde \zeta$ with fixed $\tilde \zeta > 0$ (this corresponds to a dual limit $(\eps,m) \to (0,\infty)$ in \eqref{eq:mth_map_3}). This is left for future work.
	\end{remark}
	
	%
	%
	

	\subsection{Fast-slow maps in standard form}
	\label{sub:fast-slow_maps_in_standard_form}
	
	In this section we consider a particularly important subclass of fast-slow maps \eqref{eq:gen_maps}, namely, fast-slow maps in the so-called \textit{standard form}
	\begin{equation}
		\label{eq:stnd_maps}
		\begin{split}
			x & \mapsto \bar x = x + \eps \tilde g(x,y,\eps) , \\
			y & \mapsto \bar y = y + \tilde f(x,y,\eps) ,
		\end{split}
	\end{equation}
	where the `fast-slow' structure is explicit in the $\eps-$factorisation of $\bar x - x = \eps \tilde g(x,y,\eps)$, which leads to a global separation into `slow variables' $x \in \mathbb R^k$, and `fast variables' $y \in \mathbb R^{n-k}$. Of course, time-scale terminology like `fast-slow' should be understood here only by analogy to the corresponding ODE systems: $x$ is a `slow variable' in the sense that successive iterates will generically be $O(\eps)-$close to one another, while $y$ is a `fast variable' in the sense that successive iterates will generically be separated by distances of $O(1)$. Many applications arise naturally in the standard form \eqref{eq:stnd_maps}, including those considered in detail in Sections \ref{sub:Chialvo_model} and \ref{sub:Poincare_maps}. In particular, allowing for slow evolution in a systems parameters leads to a fast-slow map in standard form \eqref{eq:stnd_maps}.
	
	\
	
	

	
	In the following we assume that $\tilde g:\mathbb R^n \times [0,\eps_0) \to \mathbb R^k$ and $\tilde f:\mathbb R^n \times [0,\eps_0) \to \mathbb R^{n-k}$ in \eqref{eq:stnd_maps} are $C^{r\geq 1}-$smooth in all arguments (although $\tilde g$ need only be $C^{r-1}$ in $\eps$), in which case the map \eqref{eq:stnd_maps} can be written in the more general form \eqref{eq:gen_maps} by expanding
	\[
	\tilde f(x,y,\eps) = \tilde f(x,y,0) + \eps \tilde f_{rem}(x,y,\eps)
	\]
	so that
	\[
	z =
	\begin{pmatrix}
		x \\ 
		y
	\end{pmatrix}
	\mapsto 
	\begin{pmatrix}
		\bar x \\ 
		\bar y
	\end{pmatrix}
	=
	\begin{pmatrix}
		x \\ 
		y
	\end{pmatrix}
	+
	\begin{pmatrix}
		O_{k,n-k} \\ 
		I_{n-k}
	\end{pmatrix}
	\tilde f(x,y,0) + \eps
	\begin{pmatrix}
		\tilde g(x,y,\eps) \\ 
		\tilde f_{rem}(x,y,\eps) 
	\end{pmatrix}
	.
	\]
	Thus, by writing $z = (x,y)^\textnormal{T} \in \mathbb R^n$, the class of maps \eqref{eq:stnd_maps} can be considered as an important subclass of fast-slow maps \eqref{eq:gen_maps} for which we have
	\begin{equation}
		\label{eq:NfG}
		N(z) =
		\begin{pmatrix}
			O_{k,n-k} \\ 
			I_{n-k}
		\end{pmatrix}
		, \qquad 
		f(z) = \tilde f(x,y,0), \qquad 
		G(z,\eps) =
		\begin{pmatrix}
			\tilde g(x,y,\eps) \\ 
			\tilde f_{rem}(x,y,\eps) 
		\end{pmatrix}
		.
	\end{equation}
	Assumptions \ref{ass:1}-\ref{ass:factorisation} are satisfied for fast-slow maps in standard form \eqref{eq:stnd_maps} if the level set $S = \{(x,y) : \tilde f(x,y,0) = O_{n-k}\}$ exists. In the following we provide corollaries of the more general notions derived so far for the special case of standard form maps \eqref{eq:stnd_maps}.

	\begin{remark}
		The class of fast-slow maps in standard form \eqref{eq:stnd_maps} should be distinguished from the class of maps in the form
		\begin{equation}
			\label{eq:2D_maps}
			\begin{split}
				x & \mapsto \bar x = \eps \tilde g(x,y,\eps) , \\
				y & \mapsto \bar y = \tilde f(x,y,\eps) ,
			\end{split}
		\end{equation}
		which also arise in applications; see e.g.~the H\'enon map in \cite[eqn.~(14.11)]{Kuehn2015} or the detailed analysis of invariant manifolds and corresponding foliations in 2-dimensional maps \eqref{eq:2D_maps} in \cite[Appendix A]{Szmolyan2004}. Generically, the maps \eqref{eq:2D_maps} are regularly (as opposed to singularly) perturbed in the sense of Definition \ref{def:singular_pert_map}, since the critical set $\{(O_k^{\textnormal T},y) : f(O_k^{\textnormal T},y,0) = O_{n-k} \}$ is generically empty or comprised entirely of isolated points $(O_k^{\textnormal T},y_\ast)$ for which $y_\ast$ satisfies $\tilde f(O_k^{\textnormal T},y_\ast,0) = O_{n-k}$. It follows that the theory developed herein does not in general apply for maps in the form \eqref{eq:2D_maps}, since they may not satisfy Assumption \ref{ass:1}.
	\end{remark}

	\subsubsection{Layer map}
	
\begin{figure}[t!]
		\centering
		\subfigure[Standard.]{\includegraphics[width=.45\textwidth]{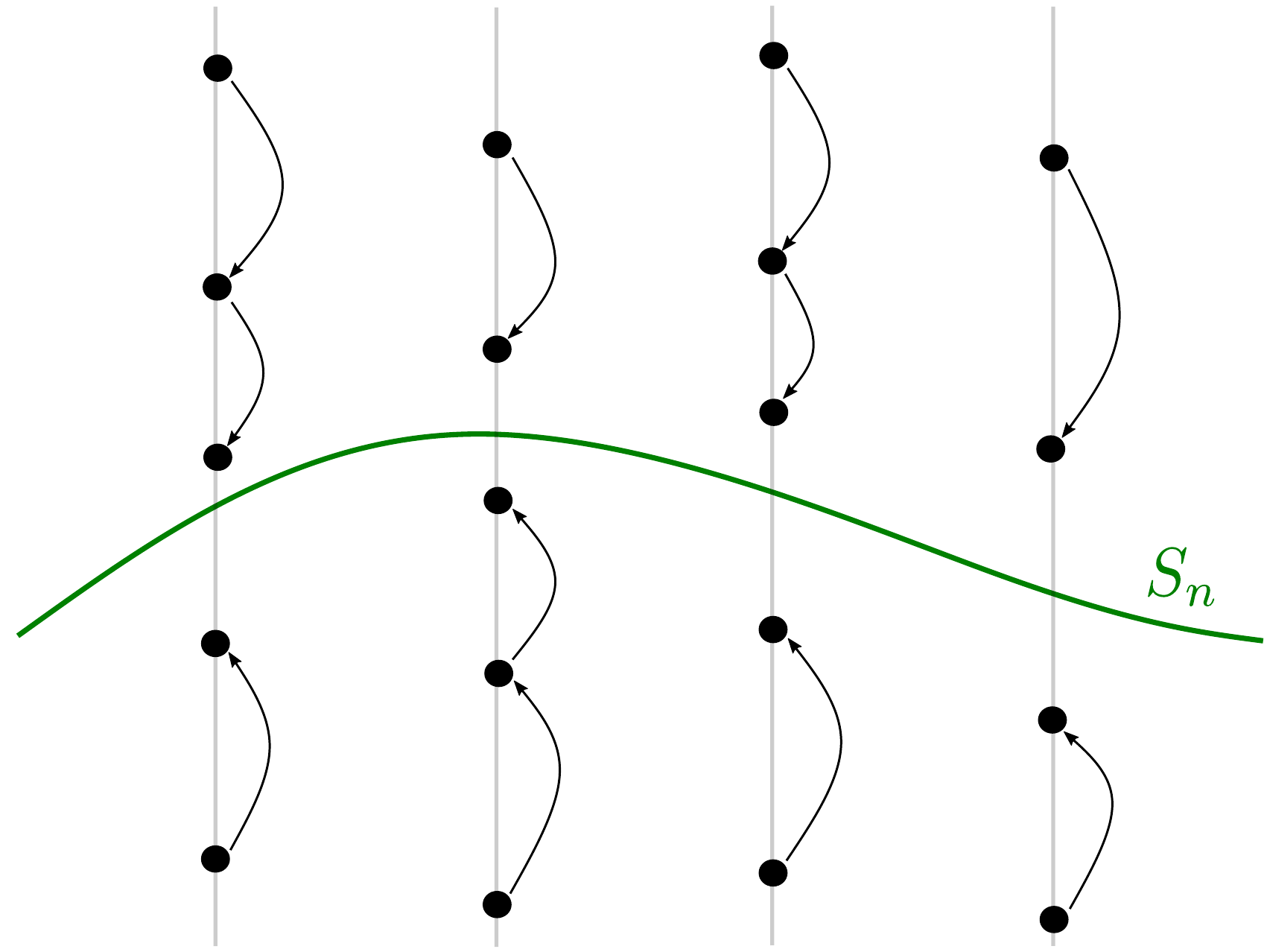}}
		\qquad
		\subfigure[Non-standard only.]{\includegraphics[width=.45\textwidth]{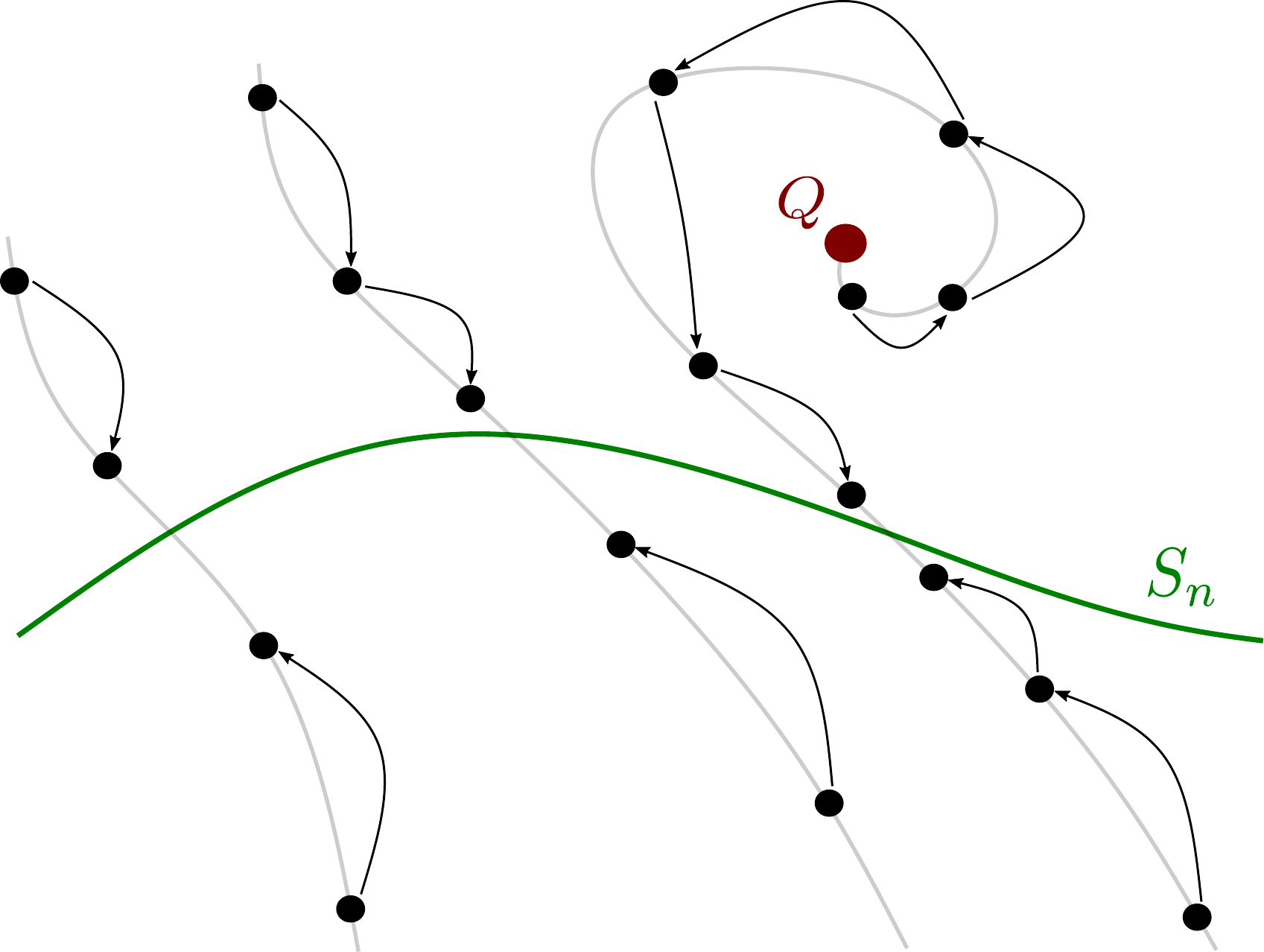}}
		\caption{Fast-slow maps in the standard form \eqref{eq:stnd_maps} can be characterised geometrically as the subclass of fast-slow maps in general form \eqref{eq:gen_maps} for which the fast foliation has been or can be globally rectified. Such a situation is sketched in (a). In (b) we show an example of a fast-slow map with an isolated fixed point $Q$. Such a map cannot be written in standard form \eqref{eq:stnd_maps} since the fast foliation cannot be globally rectified.}
		\label{fig:layer_problems}
	\end{figure}

	The layer map for \eqref{eq:stnd_maps} is
	\begin{equation}
		\label{eq:layer_map_stnd}
		\begin{split}
			& x \mapsto \bar x = x , \\
			& y \mapsto \bar y = y + \tilde f(x,y,0) .
		\end{split}
	\end{equation}
	Notice that the slow variables $x \in \mathbb R^k$ become parameters. This is a direct consequence of the separation of slow and fast variables in \eqref{eq:stnd_maps}, and not true generally for fast-slow maps \eqref{eq:gen_maps}, recall the layer map \eqref{eq:layer_map}. Geometrically, we see that the subclass of fast-slow maps in standard form \eqref{eq:stnd_maps} consists of the maps \eqref{eq:gen_maps} for which the fast foliation defined by the layer map \eqref{eq:layer_map} is globally rectified; see Figure \ref{fig:layer_problems}.
	
	The critical manifold is
	\[
	S = \left\{(x,y) : \tilde f(x,y,0) = f(x,y) = O_{n-k} \right\} ,
	\]
	and it follows from Proposition \ref{prop:EVs} that the $n-k$ non-trivial multipliers along $S$ coincide with the $n-k$ eigenvalues of the matrix
	\[
	I_{n-k} + N Df|_S = I_{n-k} + D_y f|_S .
	\]

	\subsubsection{Reduced map}
	
	It follows from equations \eqref{eq:proj} and \eqref{eq:NfG} that the projection operator $\Pi^S_{\mathcal N}$ associated with the map \eqref{eq:stnd_maps} has matrix representation
	\[
	\Pi^S_{\mathcal N} = I_n - N (Df N)^{-1} Df =
	\begin{pmatrix}
		I_k & O_{k,n-k} \\
		-(D_yf)^{-1}D_xf & O_{n-k,n-k}
	\end{pmatrix} .
	\]
	Hence the reduced map is given by
	\[
	\begin{split}
		& x \mapsto \bar x = x + \eps \tilde g(x,y,0) , \\
		& y \mapsto \bar y = y - \eps (D_yf)^{-1} D_xf \tilde g(x,y,0) ,
	\end{split}
	\qquad (x,y) \in S .
	\]
	and for each fixed $0 < \eps \ll 1$ there is a reduced $\eps^{-1}\tilde m = m$'th iterate map \eqref{eq:mth_map} defined by
	\[
	\begin{split}
		& x \mapsto \bar x^{m} = x + \tilde m \tilde g(x,y,0) , \\
		& y \mapsto \bar y^{m} = y - \tilde m (D_yf)^{-1} D_xf \tilde g(x,y,0) ,
	\end{split}
	\qquad (x,y) \in S .
	\]
	One can clearly see a correspondence with the well-known expression for the reduced problem for fast-slow ODEs in standard form, see e.g.~\cite{Jones1995,Kaper1999,Kuehn2015}. 

	\subsubsection{Local transformation to standard form}
	
	Similarly to the ODE case, there is a local equivalence between fast-slow maps in standard and non-standard form. We stress again that such an equivalence is only local, and valid only in sufficiently small neighbourhoods about points $p \in S$.
	
	\begin{proposition}
		\label{prop:stnd_general_relation}
		Consider the fast-slow map \eqref{eq:gen_maps} and let $p \in S$ be a normally hyperbolic point. For all $\eps \in [0,\eps_0)$ with $\eps_0$ sufficiently small, there exists a smooth invertible local coordinate transformation such that \eqref{eq:gen_maps} takes the standard form \eqref{eq:stnd_maps} in a sufficiently small neighbourhood about $p$.
	\end{proposition}
	
	\begin{proof}
		We may proceed via arguments similar to those applied in the proof of corresponding results for ODEs \cite[Section 3.7]{Wechselberger2019}. We assume the same local conditions as in the proof of Proposition \ref{prop:EVs}, in particular that the $(n-k) \times (n-k)$ matrix $D_yf|_S$ is regular so that $S$ is given locally as a graph $y = \varphi_0(x)$. Due to normal hyperbolicity, the nonlinear fast fibers $\mathcal F$ are given as constant level sets $L(x,y) = c$, for a $C^r-$smooth function $L:\mathbb R^k \times \mathbb R^{n-k} \to \mathbb R^{n-k}$ such that the $(n-k) \times k$ matrix $D_xL$ is locally regular and the following invariance property is satisfied:
		\begin{equation}
			\label{eq:fiber_inv}
			L(x + N^x(x,y) f(x,y), y + N^y(x,y) f(x,y) ) = L(x,y) .
		\end{equation}
		The function $L(x,y)$ is guaranteed to exist by the center manifold theorem, since for $\eps = 0$ each nonlinear fiber $\mathcal F$ is just the union of stable and unstable manifolds for $p \in S$ considered as a fixed point of the layer map \eqref{eq:layer_map}. The foliation for $\eps = 0$ can be rectified via the coordinate transformation
		\[
		s = L(x,y) ,
		\]
		which has a local inverse $x = K(s,y)$ by the inverse function theorem since $D_xL$ is locally regular. Using the invariance property \eqref{eq:fiber_inv} and expanding about $\eps = 0$ we obtain the map
		\begin{equation}
			\label{eq:stnd_form_after_trans}
			\begin{split}
				s &\mapsto s + \eps DLG(K(s,y),y,\eps) + O(\eps^2) , \\
				y &\mapsto y + N^y(K(s,y),y) f(K(s,y),y) + \eps G^y(K(s,y),y,\eps) ,
			\end{split}
		\end{equation}
		where $DLG = D_xL G^x + D_yL G^y$. The map \eqref{eq:stnd_form_after_trans} is in the standard form \eqref{eq:stnd_maps}.
	\end{proof}
	
	
	\begin{remark}
		\label{rem:Fenichel_normal_form}
		In the proof of Proposition \ref{prop:stnd_general_relation} we rectified the fast foliation in a neighbourhood of a normally hyperbolic point of $S$. By a subsequent application of the coordinate transformation used in the proof of Proposition \ref{prop:EVs}, one can also rectify the critical manifold $S$ in such neighbourhoods in order to obtain a local normal form analogous to the well-known (local) Fenichel normal form \cite{Fenichel1979,Jones1995,Kuehn2015} in the continuous-time setting.
	\end{remark}
	
	\begin{remark}
		Proposition \ref{prop:stnd_general_relation} asserts the existence of a local coordinate transformation putting general non-standard form maps \eqref{eq:gen_maps} locally (not globally) into standard form \eqref{eq:stnd_maps}. 
		It is worthy to note that the problem of obtaining an explicit form of the function $L(x,y)$ used in the transformation for a given application is typically very difficult or intractable. A direct application of the coordinate-independent theory developed in earlier sections is typically preferred on practical grounds in such cases.
	\end{remark}

	\section{Slow manifold theorems}
	\label{sec:slow_manifold_theorems}

	In this section we state the main results. Results pertaining to persistence of normally hyperbolic critical manifolds as locally invariant slow manifolds, as well as persistence of the stable/unstable manifolds \eqref{eq:layer_foliation} and their locally invariant foliations are given.  
	For expository reasons we have decided to state the results via a series of independent statements, similarly to the presentation of Fenichel's theorems for fast-slow ODEs in \cite{Jones1995}. We shall also -- again similarly again to \cite{Jones1995} -- present two versions of the result describing perturbations of the critical manifold, the latter being a specialisation to the case in which the obtained slow manifold/foliation has a graph representation. Though the latter (graph) formulation is less general and easily derived from the former (manifold) formulation, it is worthwhile to present both here for two reasons. Firstly, the critical manifold $S$ is frequently a graph in applications. Secondly, a graph formulation is always achievable locally. This fact is also leveraged in the proofs in Section \ref{sec:proof_of_theorem_fenichel}, where we shall typically restrict to the analysis of local graph formulations since the more global statements for compact manifolds are obtained via standard arguments based on a partition of unity.
	
	\
	
	We begin with the existence of locally invariant slow manifolds obtained as perturbations of the critical manifold. We shall assume throughout that the critical manifold $S$ is normally hyperbolic. This saves us from restricting to normally hyperbolic submanifolds $S_n \subseteq S$ as in Section \ref{sec:a_coordinate-independent_framework_for_fast-slow_maps}. All main results are stated for general maps \eqref{eq:gen_maps}, stated again here for convenience,
	\begin{equation}
		\label{eq:main_map}
		z \mapsto \bar z = H(z,\eps) = z + N(z) f(z) + \eps G(z,\eps) ,
	\end{equation}
	where $H(z,\eps)$ is $C^{r\geq 1}-$smooth on an open set $\mathcal U \times (-\eps_0,\eps_0)$ with $\mathcal U \subset \mathbb R^n$ and subject to Assumptions \ref{ass:1}-\ref{ass:factorisation}.

	\begin{thm}
		\label{thm:slow_manifolds}
		\textup{(Existence of slow manifolds)} Consider the map \eqref{eq:main_map} under Assumptions \ref{ass:1}-\ref{ass:factorisation}. Let $S$ be a compact, connected and normally hyperbolic critical manifold. Then there is an $\eps_0 > 0$ such that for all $\eps \in (0,\eps_0)$ there exists a compact connected manifold $S_{\eps}$ which is
		\begin{enumerate}
			\item[(i)] $O(\eps)-$close and diffeomorphic to $S$;
			\item[(ii)] $C^r-$smooth in both $z$ and $\eps$;
			\item[(iii)] locally invariant under the map \eqref{eq:main_map}, i.e.~the restricted map $H|_{S_\eps}$ is invertible and satisfies the following: if $z \in S_\eps$ and $H^j(z,\eps) \in \mathcal U$ for all $j = 1, \ldots , l$, then
			\[
			\bar z^j = H^j(z,\eps) \in S_\eps 
			\]
			for all $j = 1, \ldots , l$, and if $z \in S_\eps$ and $H^{-j}(z,\eps) \in \mathcal U$ for all $j = 1, \ldots , l$, then
			\[
			\bar z^{-j} = H^{-j}(z,\eps) \in S_\eps 
			\]
			for all $j = 1, \ldots , l$.
		\end{enumerate}
	\end{thm}
	
\begin{figure}[t!]
		\centering
		\includegraphics[scale=0.45]{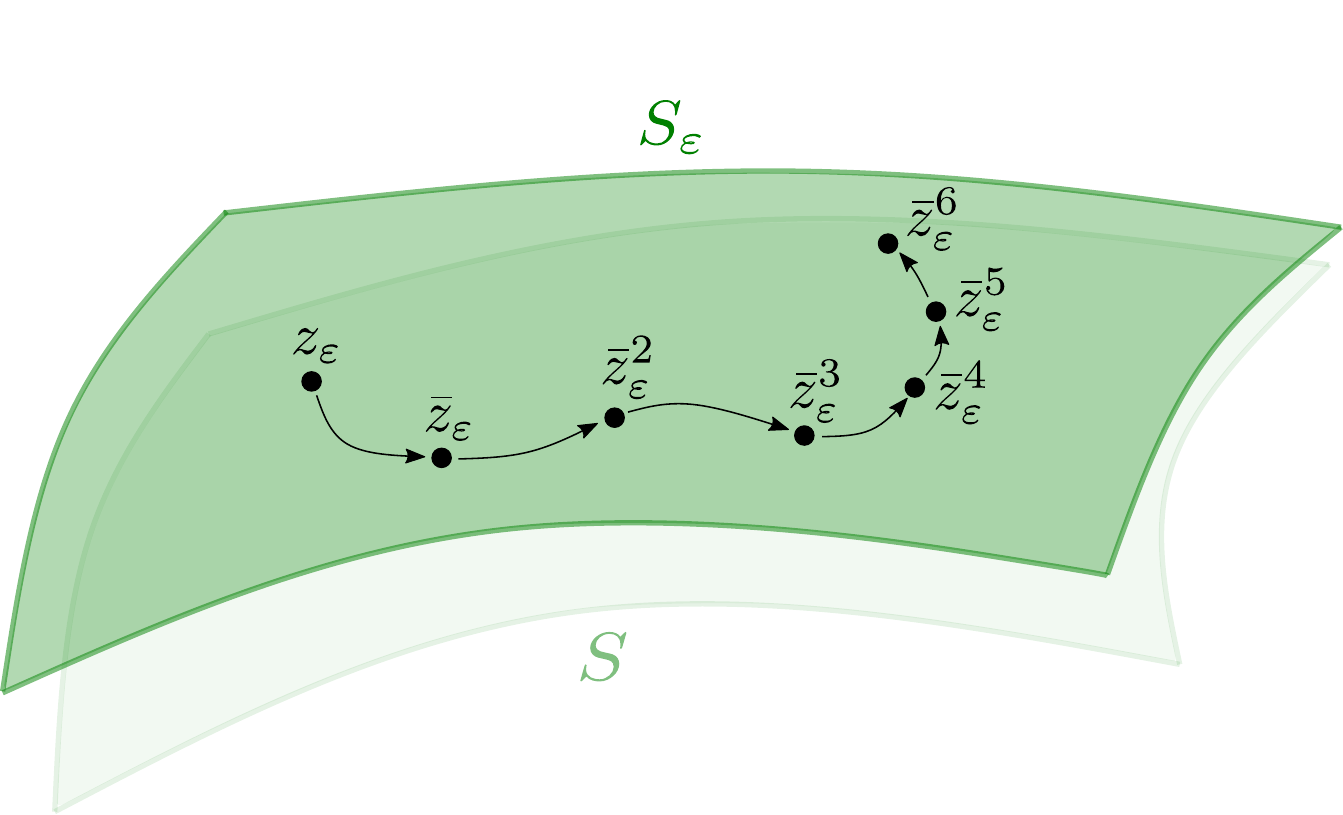}
		\caption{A normally hyperbolic critical manifold $S$ for the map \eqref{eq:main_map} perturbs to an $O(\eps)-$close locally invariant \textit{slow manifold} $S_\eps$ for $0 < \eps \ll 1$. This is described in Theorem \ref{thm:slow_manifolds}. Dynamics within $S_\eps$ are governed by the map \eqref{eq:red_map_expanded}, i.e.~approximated to leading order in $\eps$ by the reduced map \eqref{eq:red_map}. Several iterates $\bar z_\eps, \bar z_\eps^2, \ldots , \bar z_\eps^6 \in S_\eps$ starting from a point $z_\eps \in S_\eps$ are shown for illustrative purposes. Theorem \ref{thm:graph_slow_manifolds} provides further details for the case in which $S$ is has a graph representation $y = \varphi_0(x)$, in which case $S_\eps$ can be calculated up to $O(\eps^2)$ using the formula \eqref{eq:slow_manifold_approx}.}
		\label{fig:slow_manifold}
	\end{figure}

	The situation is sketched in Figure \ref{fig:slow_manifold}. Note that Theorem \ref{thm:slow_manifolds} also justifies the key assumption for the derivation of the reduced map \eqref{eq:red_map} via Proposition \ref{prop:reduced_map}, namely, the existence of a locally invariant slow manifold $S_{\eps}$ perturbing from $S$. Taken together, Theorem \ref{thm:slow_manifolds} implies the existence of a locally invariant slow manifold $S_\eps$ that is $O(\eps)-$close to $S$, and Proposition \ref{prop:reduced_map} provides an explicit form for the map governing the dynamics on $S_\eps$.
	
	It is also worthy to note that in general, slow manifolds described by Theorem \ref{thm:slow_manifolds} are non-unique. This fact, which is also true for Fenichel slow manifolds in the ODE setting and for center-type manifolds more generally \cite{Fenichel1979,Jones1995,Kuehn2015}, follows from the use of $C^\infty$ cutoff functions in order to control the dynamics entering or leaving a neighbourhood of $S$ along the center directions. 
	As a consequence, any two slow (resp.~Fenichel, center) manifolds must have the same Taylor series in $\eps$, however they may still be separated by $O(e^{-c/\eps})$ distances beyond all orders, where $c>0$. Because of the closeness of slow manifolds, we shall frequently refer to any fixed choice of slow manifold $S_{\eps}$ described by Theorem \ref{thm:slow_manifolds} as `the' slow manifold, as has become common nomenclature in fast-slow ODE theory. 
	
	
		
	
	\
	
	In order to provide a graph formulation of Theorem \ref{thm:slow_manifolds}, we impose an additional assumption.
	
	\begin{assumption}
		\label{ass:graph_form}
		The critical manifold $S$ is given as a graph of a $C^r-$smooth function $\varphi_0: K \to \mathbb R^k$ over a compact, simply connected domain $\mathcal K \subset \mathbb R^{n-k}$, the boundary of which is a $C^\infty$ $(k-1)-$dimensional submanifold $\partial \mathcal K$. We write
		\[
		S = \left\{(x,\varphi_0(x)) : x \in \mathcal K \right\}.
		\]
	\end{assumption}
	
	Assumption \ref{ass:graph_form} is actually quite natural, in the sense that it can always be satisfied locally after a change of notation if necessary since the matrix $Df|_S$ has full row rank by Assumption \ref{ass:1}. In order to consider maps \eqref{eq:main_map} satisfying Assumption \ref{ass:graph_form}, we reintroduce the componentwise notation used already in the proofs of Propositions \ref{prop:EVs}, \ref{prop:reduced_map}, \ref{prop:red_jump_map} and \ref{prop:stnd_general_relation}, i.e.~we write \eqref{eq:main_map} as
	\begin{equation}
		\label{eq:main_map_graph}
		\begin{pmatrix}
			x \\
			y
		\end{pmatrix}
		\mapsto
		\begin{pmatrix}
			\bar x \\
			\bar y
		\end{pmatrix}
		=
		\begin{pmatrix}
			H^x(x,y,\eps) \\
			H^y(x,y,\eps)
		\end{pmatrix}
		=
		\begin{pmatrix}
			x \\
			y
		\end{pmatrix}
		+
		\begin{pmatrix}
			N^x(x,y) \\
			N^y(x,y)
		\end{pmatrix}
		f(x,y) + \eps
		\begin{pmatrix}
			G^x(x,y,\eps) \\
			G^y(x,y,\eps)
		\end{pmatrix} ,
	\end{equation}
	where as before $x \in \mathbb R^k$, $y\in \mathbb R^{n-k}$, $N^x(x,y)$ and $N^y(x,y)$ are matrices of dimensions $k \times (n-k)$ respectively $(n-k) \times (n-k)$, $G^x(x,y,\eps)$, $G^y(x,y,\eps)$ are column vectors of length $k$ respectively $n-k$, and $f(x,y) = y - \varphi_0(x)$ after redefining $N(x,y)$ if necessary.

	\begin{thm}
		\label{thm:graph_slow_manifolds}
		\textup{(Existence of slow manifolds as graphs)} Consider the map \eqref{eq:main_map_graph} under Assumptions \ref{ass:1}, \ref{ass:factorisation} and \ref{ass:graph_form}.	Then there exists a $C^r-$smooth function $\varphi_\eps : \mathcal K \times [0,\eps_0) \to \mathbb R^k$ such that the slow manifold is also a graph
		\[
		S_\eps = \left\{(x,\varphi_\eps(x)) : x \in \mathcal K \right\} ,
		\]
		where in particular
		\begin{equation}
			\label{eq:slow_manifold_approx}
			\varphi_\eps(x) = \varphi_0(x) - \eps (D_yf)^{-1} (Df N)^{-1} Df G + O(\eps^2) .
		\end{equation}
		Moreover, $S_\eps$ is locally invariant in the following sense: for each $(x,y) \in S_\eps$ the invariance equation
		\begin{equation}
			\label{eq:graph_sm_forward_inv}
			\bar y = \varphi_\eps\left(H^x(x,\varphi_\eps(x),\eps) \right) = H^y\left(x, \varphi_\eps(x),\eps \right) 
		\end{equation}
		is satisfied as long as $H(x,y,\eps) \in \mathcal U$, and the invariance equation
		\begin{equation}
			\label{eq:graph_sm_backward_inv}
			\bar y^{-1} = \varphi_\eps\left((H^x)^{-1}(x,\varphi_\eps(x),\eps) \right) = (H^y)^{-1}\left(x, \varphi_\eps(x),\eps \right) 
		\end{equation}
		is satisfied as long as $H^{-1}(x,y,\eps) \in \mathcal U$.
	\end{thm}
	
	Theorem \ref{thm:slow_manifolds} can be shown to follow from the graph formulation in Theorem \ref{thm:graph_slow_manifolds} using compactness and a partition of unity. Hence in Section \ref{sec:proof_of_theorem_fenichel} it is sufficient to prove Theorem \ref{thm:graph_slow_manifolds}. Slow manifold properties such as smoothness follow from the existence of the $C^r-$smooth function $\varphi_\eps$. Equation \eqref{eq:slow_manifold_approx} in Theorem \ref{thm:graph_slow_manifolds} provides the asymptotics for the slow manifold $S_\eps$. In particular, this expression constitutes a local formula for $S_\eps$ up to $O(\eps^2)$ in terms of the initial data $N$, $f$ and $G$.
	
	\
	
	We turn now to persistence of stable and unstable manifolds $W_{loc}^s(S)$ and $W_{loc}^u(S)$ respectively, as defined for $\eps = 0$ in \eqref{eq:layer_foliation}. Similarly to the case of fast-slow ODEs, both stable and unstable manifolds can be shown to perturb to nearby manifolds with certain local invariance properties. 
	
	\begin{thm}
		\label{thm:stable_manifolds}
		\textup{(Persistence of stable/unstable manifolds)} Consider the map \eqref{eq:main_map} under Assumptions \ref{ass:1}-\ref{ass:factorisation}. Then for all $\eps \in (0,\eps_0)$ with $\eps_0 > 0$ sufficiently small there exists manifolds $W_{loc}^s(S_{\eps})$ and $W_{loc}^u(S_{\eps})$ that are
		\begin{enumerate}
			\item[(i)] $O(\eps)-$close and diffeomorphic to $W_{loc}^s(S)$ and $W_{loc}^u(S)$ respectively;
			\item[(ii)] $C^r-$smooth in both $z$ and $\eps$;
			\item[(iii)] locally positively and negatively invariant under the map \eqref{eq:gen_maps}, respectively. More precisely, the restricted map $H|_{W^u_{loc}(S_\eps)}$ is invertible and satisfies the following: 
			if $z \in W_{loc}^s(S_{\eps})$ and $H^j(z,\eps) \in \mathcal U$ for all $j = 1, \ldots , l$, then
			\[
			\bar z^j = H^j(z,\eps) \in W_{loc}^s(S_{\eps}) 
			\]
			for all $j = 1, \ldots , l$. Similarly, if 
			$z \in W_{loc}^u(S_{\eps})$ and $H^{-j}(z,\eps) \in \mathcal U$ for all $j = 1, \ldots , l$, then
			\[
			\bar z^{-j} = H^{-j}(z,\eps) \in W_{loc}^u(S_{\eps}) ,
			\]
			for all $j = 1, \ldots , l$.
		\end{enumerate}
		Moreover, the slow manifold $S_\eps$ lies within the intersection of the perturbed stable and unstable manifolds, i.e.
		\[
		S_\eps = W_{loc}^{s}(S_\eps)\cap W_{loc}^{u}(S_\eps).
		\]
	\end{thm}
	
\begin{figure}[t!]
		\centering
		\includegraphics[scale=0.45]{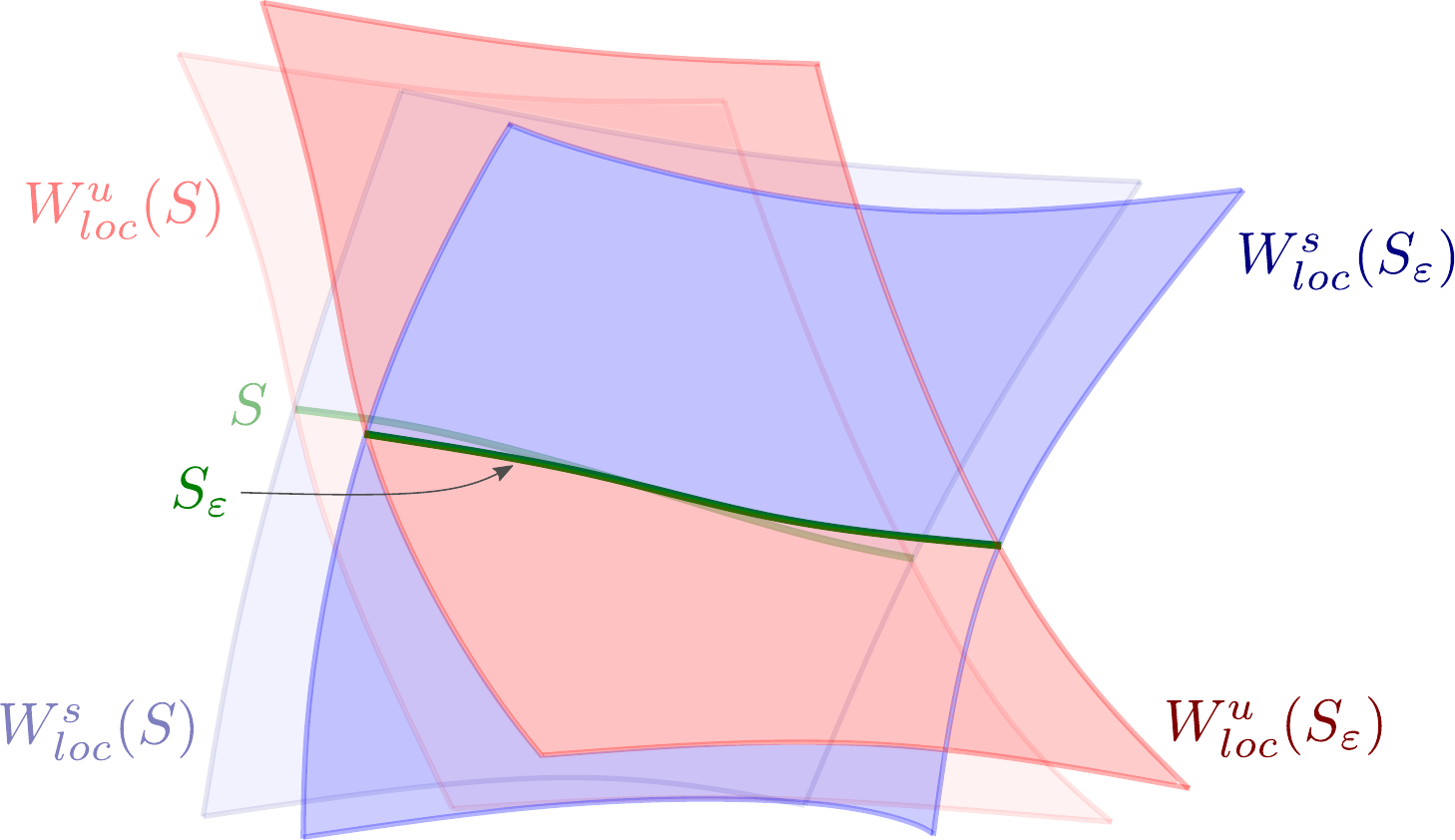}
		\caption{The stable/unstable manifolds $W^{u/s}_{loc}(S)$ (here in shaded blue/red) of a normally hyperbolic critical manifold $S$ (shaded green) perturb to $O(\eps)-$close positively/negatively invariant manifolds $W^{s/u}_{loc}(S_\eps)$ (blue/red) which intersect along the slow manifold $S_\eps$ (green), see Theorem \ref{thm:stable_manifolds}.}
		\label{fig:Wsu_perturbed}
	\end{figure}

	The situation is sketched in Figure \ref{fig:Wsu_perturbed}. The persistence of stable and unstable manifolds $W_{loc}^{s/u}(S)$ as nearby locally invariant manifolds $W_{loc}^{s/u}(S_\eps)$ described by Theorem \ref{thm:stable_manifolds} is essentially analogous to the persistence of normally hyperbolic critical manifolds $S$ as nearby locally invariant slow manifolds $S_\eps$ in Theorem \ref{thm:slow_manifolds}. In particular (and for similar reasons), the manifolds $W_{loc}^{s/u}(S_\eps)$ are also non-unique but $O(e^{-c/\eps})-$close for some constant $c > 0$. The final assertion that $S_\eps$ lies in the intersection $W_{loc}^{s}(S_\eps)\cap W_{loc}^{u}(S_\eps)$ explains the similarities between Theorems \ref{thm:slow_manifolds} and \ref{thm:stable_manifolds}. In fact, as is typical in proofs for the existence of normally hyperbolic manifolds more generally, Theorem \ref{thm:slow_manifolds} will be derived (via Theorem \ref{thm:graph_slow_manifolds}) in Section \ref{sec:proof_of_theorem_fenichel} as a direct consequence of Theorem \ref{thm:stable_manifolds}.
	
	\
	
	Finally and similarly to the case for ODEs, it turns out the manifolds $W_{loc}^{s/u}(S_{\eps})$ admit an invariant foliation by smooth fibers with base points on $S_{\eps}$, i.e.~the invariant foliation of the manifolds $W_{loc}^{u/s}(S)$ by unstable/stable manifolds of points on $S$, see again \eqref{eq:foliations_perturbed}, also persist. The family of stable/unstable fibers also exhibit useful invariance properties, and a quantitative estimate for the contraction/expansion rate along fibers can be given in terms of spectral information. 
	
	We first require a little more notation. Recall from Definition \ref{def:nh} that the non-trivial multipliers of a $k-$dimensional normally hyperbolic critical manifold $S$ are denoted by $\mu_j(z)$, where $z \in S$ and $j = 1, \ldots, n-k$. By normal hyperbolicity, $|\mu_j(z)| \neq 1$ for all $j$. Denote the $n_a \leq n-k$ stable multipliers with $|\mu_{j}(z)| < 1$ by $\mu_{a,j}(z)$, and the $n_r \leq n-k$ unstable multipliers with $|\mu_{j}(z)| > 1$ by $\mu_{r,j}(z)$, where $n_a + n_r = n-k$. Assuming that $S$ is compact, we may define spectral bounds
	\begin{equation}
		\label{eq:spectral_bounds_original}
		\nu_A := \sup_{z \in S , j = 1, \ldots , n_a} | 
		\mu_{a,j}(z) | < 1 , \qquad 
		\nu_R := \inf_{z \in S , j = 1, \ldots , n_r} | 
		\mu_{r,j}(z) | > 1 .
	\end{equation}
	Contraction and repulsion along stable and unstable fibers respectively can be quantified in terms of $\nu_A$ and $\nu_R$. We now state the main result on the persistence of stable and unstable foliations.

	\begin{figure}[t!]
		\centering
		\subfigure[Orientation preserving case.]{\includegraphics[width=.7\textwidth]{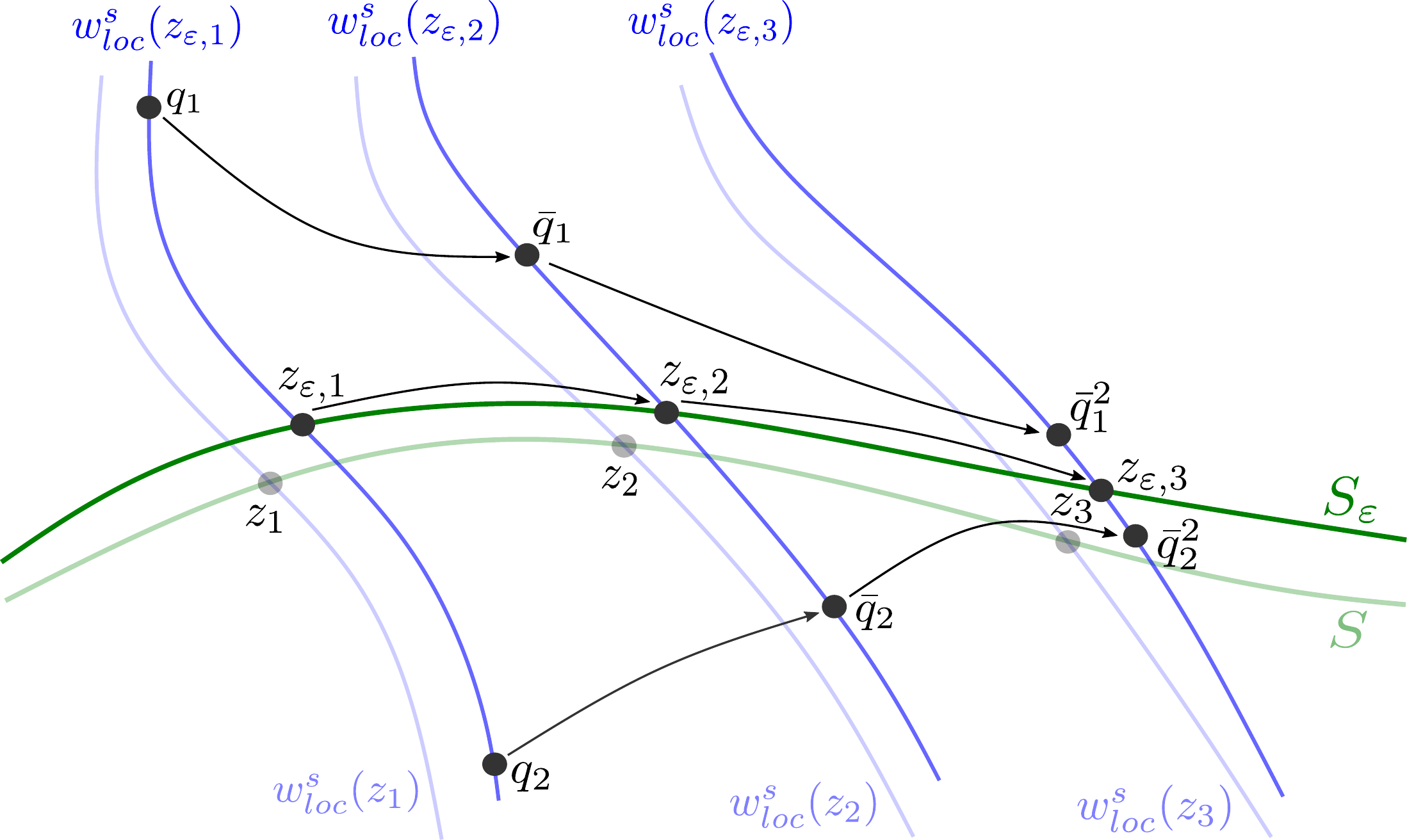}}
		
		\
		
		\subfigure[Orientation reversing case.]{\includegraphics[width=.7\textwidth]{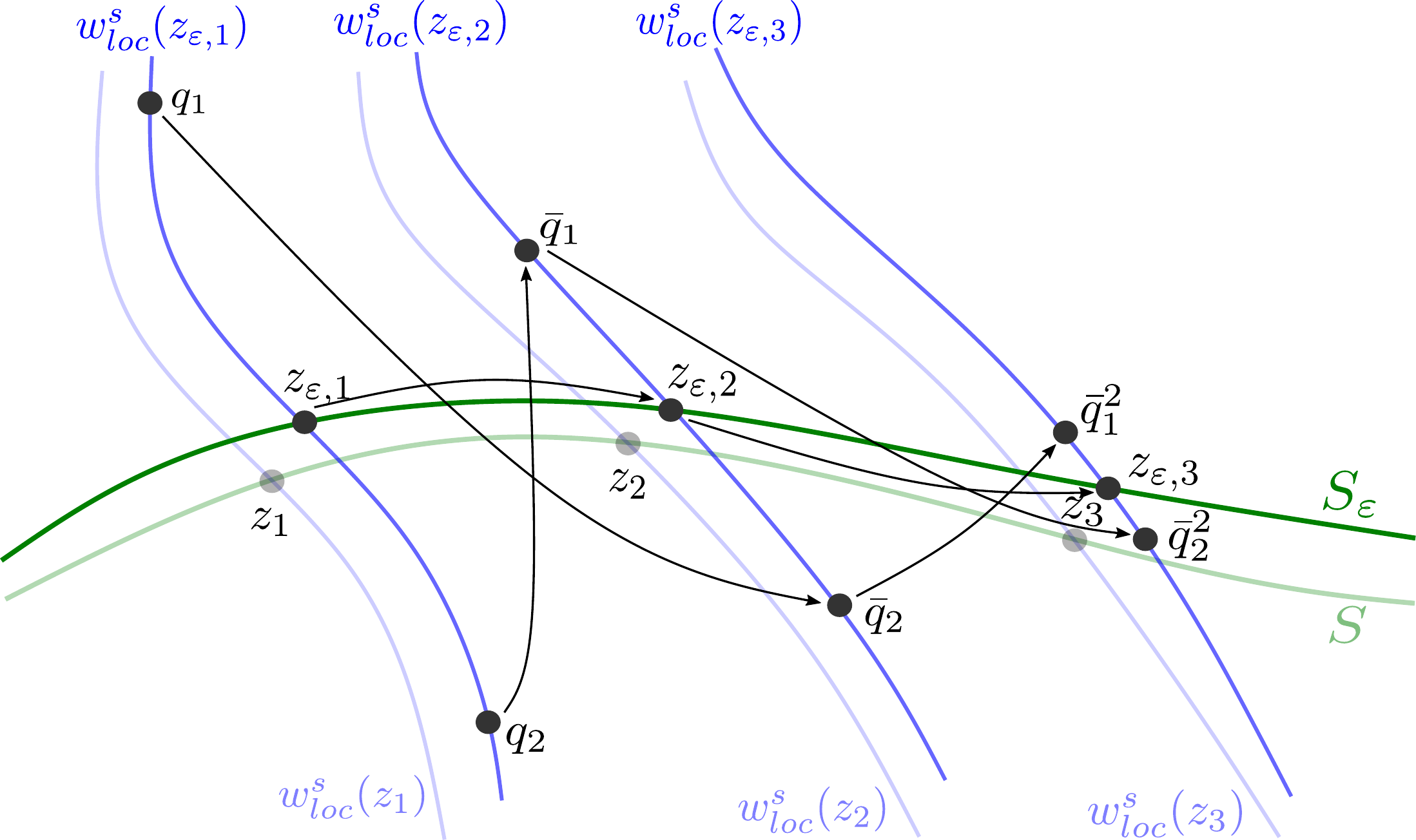}}
		\caption{Positive invariance of the perturbed stable fibers $\{w_{loc}^s(z_\eps)\}_{z_\eps \in S_\eps}$ of a $1-$dimensional critical manifold in dimension $n=2$. By the invariance property in Theorem \ref{thm:foliations} (i), initial conditions $q_1, q_2$ in the stable fiber $w^s_{loc}(z_{\eps,1})$ with base point $z_{\eps,1} \in S_\eps$ are mapped by $H$ into the stable fiber $w^s_{loc}(z_{\eps,2})$, then $w^s_{loc}(z_{\eps,3})$, with base points $z_{\eps,2}$ and $z_{\eps,3}$ corresponding to iterates of $z_{\eps,1}$. Moreover, iterates are exponentially contracted towards their base points on $S_\eps$ by Theorem \ref{thm:foliations} (iii). In (a): the orientation preserving case in which $S$ has a non-trivial multiplier $\mu \in (0,1)$. In (b): the orientation reversing case with $\mu \in (-1,0)$.}
		\label{fig:foliation_invariance}
	\end{figure}
	
	\begin{thm}
		\label{thm:foliations}
		\textup{(Persistence of stable/unstable foliations)} 
		Assume the same hypotheses as Theorem \ref{thm:stable_manifolds}. Then for $\eps \in (0,\eps_0)$ with $\eps_0 > 0$ sufficiently small, 
		the perturbed stable/unstable manifolds $W^{s/u}_{loc}(S_{\eps})$ of Theorem \ref{thm:stable_manifolds} admit foliations by stable/unstable fibers $w^{s/u}_{loc}(z_\eps)$ with base points $z_\eps \in S_{\eps}$, i.e.
		\begin{equation}
			\label{eq:foliations_perturbed}
			W^s_{loc}(S_{\eps}) = \bigcup_{z_\eps \in S_\eps} w^s_{loc}(z_\eps) , \qquad
			W^u_{loc}(S_{\eps}) = \bigcup_{z_\eps \in S_\eps} w^u_{loc}(z_\eps) .
		\end{equation}
		These foliations converge fiberwise to the foliations \eqref{eq:layer_foliation} in the layer map as $\eps \to 0$, and 
		satisfy the following:
		\begin{enumerate}
			\item[(i)] The stable fibers $w^s_{loc}(z_\eps)$ form a locally positively invariant family, i.e.~if $z_\eps \in S_{\eps}$ and $H(z_\eps, \eps) \in \mathcal U$, then
			\[
			H\left( w_{loc}^s(z_\eps), \eps \right) \subset w_{loc}^s \left( H(z_\eps, \eps) \right) .
			\]
			\item[(ii)] The unstable fibers $w^u_{loc}(z_\eps)$ form a locally negatively invariant family, i.e.~if $z_\eps \in S_{\eps}$ and $H^{-1}(z_\eps, \eps) \in \mathcal U$, then
			\[
			H^{-1} \left(w_{loc}^u (z_\eps), \eps \right) \subset w_{loc}^u \left( H^{-1}(z_\eps, \eps) \right) .
			\]
			\item[(iii)] For each $z_\eps \in S_{\eps}$ such that $H(z_\eps, \eps) \in \mathcal U$ for all $j = 1, \ldots , l$, contraction along stable fibers $w_{loc}^s(z_\eps)$ is exponential with rate faster than $\chi_A^j$, for a fixed constant $\chi_A \in (\nu_A, 1)$. 
			More precisely, for all $\xi_\eps \in w_{loc}^s(z_\eps)$ such that $H^j(\xi_\eps, \eps)$ and $H^j(z_\eps, \eps)$ stay in $\mathcal U$ for all $j = 1, \ldots , l$, we have that
			\[
			\| H^j(\xi_\eps, \eps) - H^j(z_\eps, \eps) \| \leq \chi_A^j \| \xi_\eps - z_\eps \| .
			\]
			\item[(iv)] For each $z_\eps \in S_{\eps}$ such that $H^{-1}(z_\eps, \eps) \in \mathcal U$ for all $j = 1, \ldots , l$, contraction along unstable fibers $w_{loc}^s(z_\eps)$ is exponential with rate faster than $\chi_R^{-j}$, for a fixed constant $\chi_R \in (1, \nu_R)$.  
			More precisely, for all $\xi_\eps \in w_{loc}^u(z_\eps)$ such that $H^{-j}(\xi_\eps, \eps)$ and $H^{-j}(z_\eps, \eps)$ stay in $\mathcal U$ for all $j = 1, \ldots , l$, we have that
			\[
			\| H^{-j}(\xi_\eps, \eps) - H^{-j}(z_\eps, \eps) \| \leq \chi_R^{-j} \| \xi_\eps - z_\eps \| .
			\]
			\item[(v)] The stable fibers $w^s_{loc}(z_\eps)$ are Lipschitz manifolds. If additionally degenerate superstability in the sense that
			\[
			\mu_{a,j}(z) = 0 
			\]
			for all $j = 1, \ldots , n_a$ only occurs at isolated points $z \in S$, then the stable fibers $w^s_{loc}(z_\eps)$ are $C^{r-1}-$smooth.
			\item[(vi)] The unstable fibers $w^u_{loc}(z_\eps)$ are $C^{r-1}-$smooth.
		\end{enumerate}
	\end{thm}

	Theorem \ref{thm:foliations} asserts that the foliation of stable/unstable manifolds $W_{loc}^{s/u}(S)$ by lower dimensional stable/unstable manifolds $w^{s/u}_{loc}(z)$ associated to points $z \in S$, also perturb in a regular fashion. The local positive/negative invariance of the families $\{w_{loc}^{s/u}(z_\eps)\}_{z_\eps \in S_\eps}$ is described in assertions (i)-(ii), and sketched for the stable foliation in Figure \ref{fig:foliation_invariance}. Statements (iii)-(iv) assert the exponential contraction and repulsion along stable and unstable fibers respectively. This explains the common use of terms like ``stable", ``unstable" and ``normally hyperbolic" when describing perturbed fibers or manifolds like $w_{loc}^{s/u}(z_\eps)$, $W^{s/u}_{loc}(S_\eps)$ and $S_\eps$. Notice by assertions (v)-(vi) that stable/unstable fibers $w_{loc}^{s/u}(z_\eps)$ are only $C^{r-1}-$smooth with respect to variation of the base point $z_\eps \in S_\eps$, even though by Theorem \ref{thm:stable_manifolds} the stable/unstable manifolds $W^{s/u}_{loc}(S_\eps)$ themselves are $C^r-$smooth. This is also true of foliations with Fenichel slow manifolds as base in the ODE setting \cite{Fenichel1979,Jones1995,Wechselberger2019,Wiggins2013}. In contrast to the ODE setting, however, an additional nondegeneracy condition is required in assertion (v) in order to infer $C^{r-1}-$smoothness of the stable fibers in particular.
	
	\begin{remark}
		\label{rem:foliation_smoothness}
		The proof of Theorem \ref{thm:foliations} assertion (v) in Section \ref{ssec:foliation_proof} relies on a smoothness result from \cite{Nipp2013} which cannot be applied in the (very) degenerate superstable case in which the critical manifold $S$ has a connected component $U_{ss}$ with $\mu_a(z) \equiv 0$ for $z \in S \cap U_{ss}$. The question of whether the nondegeneracy in Theorem \ref{thm:foliations} is a necessary condition for smoothness of the stable fibers is not considered further in this work.
	\end{remark}

	\section{Proofs}
	\label{sec:proof_of_theorem_fenichel}
	
	In the following we prove the main results from Section \ref{sec:slow_manifold_theorems}. Our results can be viewed as a specialisations of pre-existing results on the existence of normally hyperbolic manifolds and their foliations more generally. Specifically in this work, the majority of our results will be derived as specialisations of more general theorems in \cite{Nipp2013} to the case of fast-slow maps in the general nonstandard form \eqref{eq:main_map}. In order to apply the results in \cite{Nipp2013}, the equations must first be `prepared' to a certain extent, and herein lies much of the analysis. The aim is to formulate necessary and sufficient conditions in \cite{Nipp2013} in terms of coordinate-invariant properties of the map, e.g.~the multipliers along $S$, so that results obtained for normal forms can be directly related to the original map \eqref{eq:main_map}. Further work must be done to prove the local invariance of slow manifolds and foliations described by our results, since the results in \cite{Nipp2013} and related works rely on some simplifying assumptions which yield only to globally invariant objects. Finally, additional work is needed in order to derive those features which are characteristic of or specific to fast-slow maps, for example, calculations leading to direct estimates for the slow manifold parameterization, or for the asymptotic decay rates along invariant fibers.
	
	This section is structured as follows: In Section \ref{sub:preparatory_results} we introduce notation, identify suitable neighbourhoods and derive local coordinates such that the relevant inflowing/outflowing requirements for the application of results in \cite{Nipp2013} are satisfied. In Section \ref{sub:NS_conds} we verify (a suitable reformulation of) necessary and sufficient conditions for the application of these results. Finally in Section \ref{sub:NS_thms}, we apply the results from \cite{Nipp2013} and use them to prove our main results from Section \ref{sec:slow_manifold_theorems}.

	\subsection{Preparatory results}
	\label{sub:preparatory_results}
	
	We begin by identifying neighbourhoods and local coordinates suited to our needs. 
	The main task is to separate directions along which iterates of \eqref{eq:main_map} are inflowing and outflowing respectively.

	For local analyses we can assume the map \eqref{eq:main_map} is in the form \eqref{eq:main_map_graph}. Recall that this form is always achievable locally, i.e.~does not rely on Assumption \ref{ass:graph_form}. The first simplification is to rectify $S$ along the $x-$axes via the transformation $v = f(x,y)$, which has a locally defined inverse $y = K(x,v)$. As shown in the proof of Proposition \ref{prop:EVs}, this leads to
	\begin{equation}
		\label{eq:straight_S_map2}
		\begin{split}
			x & \mapsto \bar x = x + \tilde N^x(x,v) v + \eps \tilde G^x(x,v,\eps) , \\
			v & \mapsto \bar v = v + Df N(x,K(x,v)) v + \eps Df G(x,K(x,v),\eps) +  O(v^2,v\eps,\eps^2) ,
		\end{split}
	\end{equation} 
	where we have defined $\tilde N^x(x,v) := N^x(x,K(x,v))$, $\tilde G^x(x,v,\eps) := G^x(x,K(x,v),\eps)$, and $Df = (D_xf,D_yf)$. In these coordinates the critical manifold is simply
	\[
	S = \{(x,0) : x \in \mathcal K \} ,
	\]
	where we may (and will) assume that $\mathcal K \subset \mathbb R^k$ is compact, simply connected, and that $S$ is normally hyperbolic. Fast directions are encoded in (but not identified with) the variables $v \in \mathbb R^{n-k}$, since the $n-k$ multipliers of the linear part $I_{n-k} + DfN|_S$ coincide by Proposition \ref{prop:EVs} with the $n-k$ non-trivial multipliers $\mu_j$ associated with the fast directions. 
	
	\begin{lemma}
		\label{lem:partial_fenichel_normal_form}
		Fix $\eps_0 > 0$ sufficiently small. For all $\eps \in [0,\eps_0)$, there exists a local coordinate transformation transforming \eqref{eq:straight_S_map2} into
		\begin{equation}
			\label{eq:partial_fenichel_normal_form}
			\begin{split}
				x & \mapsto \bar x = x + \tilde N^x(x,N^y(a,b)^\textnormal{T}) N^y (a, b)^\textnormal{T} + \eps \tilde G^x(x,N^y(a,b)^\textnormal{T},\eps) , \\
				a & \mapsto \bar a = a + \Lambda_a(x,a,b) a + \eps (N_a^y)^{-1} Df G(x,N^y(a,b)^\textnormal{T},\eps) , \\
				b & \mapsto \bar b = b + \Lambda_r(x,a,b) b + \eps (N_r^y)^{-1} Df G(x,N^y(a,b)^\textnormal{T},\eps) ,
			\end{split}
		\end{equation}
		where $\Lambda_a(x,0,0)$ is a diagonal $n_a \times n_a$ matrix with eigenvalues $\lambda_{i}(x)$ such that
		\[
		| \mu_i(x) | = | 1 + \lambda_i(x) | < 1
		\]
		uniformly in $x$ for all $i = 1, \ldots , n_a$, and $\Lambda_r(x,0,0)$ is a diagonal $n_r \times n_r$ matrix with eigenvalues $\lambda_{i}(x)$ such that
		\[
		| \mu_i(x) | = | 1 + \lambda_i(x) | > 1
		\]
		uniformly in $x$ for all $i = 1, \ldots , n_r$.
	\end{lemma}
	
	\begin{proof}
		We use the following facts:
		\begin{itemize}
			\item For any $z \in S$,
			\[
			\imag Dh(z) = \mathcal N_z = E^s(z) \cup E^u(z) = \textrm{span\,} \{N^i(z) : i = 1, \ldots , n-k\} ,
			\]
			where the $N^i(z)$ are columns of $N(z)$.
			\item For any $z \in S$,
			\[
			\ker Dh(z) = T_zS = \textrm{span\,} \{Df_i(z) : i = 1, \ldots , n-k\}^\perp ,
			\]
			where the $Df_i(z)$ are the rows of $Df(z)$, and $\perp$ denotes the orthogonal complement.
		\end{itemize}
		It follows from these facts that for each $z \in S$, the $n \times n$ matrix $A(z) := [Df(z)^\perp \ N(z)]$ forms the matrix with eigenvector columns needed to make a Jordan decomposition of the layer map \eqref{eq:straight_S_map2}$|_{\eps = 0}$. 
		%
		%
		Here is suffices to Jordan decompose only the $v$ variables which, in the case of saddle-type critical manifolds, will allow for a splitting of attracting and repelling components in the leading order. Since column vectors of $N(z)$ can be chosen to be eigenvectors spanning $E^u(z) \cup E^s(z)$ at each $z \in S$, the (linearly independent) column vectors of the $(n-k) \times (n-k)$ matrix $N^y(z)$ are eigenvectors of the matrix $DfN(z)$. This motivates the coordinate transformation
		\[
		u = (N^y)^{-1} v \in \mathbb R^{n-k} ,
		\]
		which leads to the map
		\[
		\begin{split}
			x & \mapsto \bar x = x + \tilde N^x(x,N^yu) N^y u + \eps \tilde G^x(x,N^yu,\eps) , \\
			u & \mapsto \bar u = u + \Lambda(x,u) u + \eps (N^y)^{-1} Df G(x,N^yu,\eps) + O(u^2,\eps u, \eps^2) , 
		\end{split}
		\]
		where $I_{n-k} + \Lambda(x,0) = \textrm{diag} \{\mu_1(x), \ldots , \mu_{n-k}(x) \}$. Since it is achievable by a simple permutation of notation if necessary, we may assume without loss of generality that $|\mu_j(x)| < 1$ for $j = 1, \ldots , n_a$ and $|\mu_j(x)| > 1$ for the remaining $n_r$ multipliers with $j = n_a + 1, \ldots , n-k$ (recall that $n_a + n_r = n-k$). In this case the matrix $\Lambda$ has a block-diagonal structure, and we may write
		\[
		\begin{split}
			x & \mapsto \bar x = x + \tilde N^x(x,N^yu) N^y u + \eps \tilde G^x(x,N^yu,\eps) , \\
			a & \mapsto \bar a = a + \Lambda_a(x,u) a + \eps (N_a^y)^{-1} Df G(x,N^yu,\eps) , \\
			b & \mapsto \bar b = b + \Lambda_r(x,u) b + \eps (N_r^y)^{-1} Df G(x,N^yu,\eps) , \\
		\end{split}
		\]
		for suitably defined $N^y_a$, $N^y_r$, where $u = (a, b)^\textnormal{T} \in \mathbb R^{n_a} \times \mathbb R^{n_r}$, and the diagonal matrices $I_{n_a} + \Lambda_a(x,0)$ and $I_{n_r} + \Lambda_r(x,0)$ encode all and only the attracting and repelling multipliers respectively. 
	\end{proof}

	In order to meet the relevant inflowing/outflowing and invariance criteria, 
	the set $\mathcal K$ needs to be extended. Specifically, results in \cite{Nipp2013} apply for open neighbourhoods $X \times Y$ for which $X$ is inflowing and $Y$ is outflowing (or visa-versa) with respect to the map. As in all proofs of the center manifold theorem, escape along the slow (i.e.~center) directions poses a problem for the satisfaction of certain global invariance properties. The usual approach to controlling the slow directions is to enlarge the neighbourhood $\mathcal K$, and augment the slow dynamics by the addition of a suitable cutoff/bump function which prevents iterates from escaping this larger neighbourhood. Results in \cite{Nipp2013} which only apply for neighbourhoods with global inflowing/outflowing and invariance properties can be applied on this enlarged neighbourhood, and subsequently restricted to $\mathcal K$.
	
	The existence of a suitable enlargement $\hat{\mathcal K}$ such that $\mathcal K \subset \text{Int\ } \hat{\mathcal K}$ is guaranteed by smoothness. In particular, following the analogous setup for ODEs in \cite{Fenichel1979,Jones1995}, $\hat {\mathcal K}$ can be chosen so that $S$ is defined/extended over $\hat{\mathcal K}$ as a $C^r-$smooth graph
	\begin{equation}
		\label{eq:S_enlarged}
		\hat S = \left\{(x,\phi_0(x)) : x \in \hat{\mathcal K} \right\} ,
	\end{equation}
	such that $\phi_0|_{\mathcal K}(x) \equiv 0$, with boundary $\partial \hat{\mathcal K}$ given by $\hat \omega(x) = 0$ for a $C^\infty-$smooth function $\hat \omega$ such that $\nabla \hat \omega|_{\partial \hat{\mathcal K}} \neq 0$. In fact, let $\hat \omega(x)$ be normalised so that $n_x := \nabla \hat \omega(x)$ defines an outward pointing normal to $\partial \hat{\mathcal K}$. Now consider the map
	\begin{equation}
		\label{eq:partial_fenichel_normal_form_enlarged}
		\begin{split}
			x & \mapsto \bar x = x + \tilde N^x(x,N^yu) N^y (a, b)^\textnormal{T} + \eps \tilde 	G^x(x,N^y(a,b)^\textnormal{T},\eps) + \hat \delta \rho(x) n_x , \\
			a & \mapsto \bar a = a + \Lambda_a(x,a,b) a + \eps (N_a^y)^{-1} Df G(x,N^y(a,b)^\textnormal{T},\eps) , \\
			b & \mapsto \bar b = b + \Lambda_r(x,a,b) b + \eps (N_r^y)^{-1} Df G(x,N^y(a,b)^\textnormal{T},\eps) ,
		\end{split}
	\end{equation}
	where $\hat \delta > 0$ is a constant to be specified later on, and $\rho(x)$ is a $C^\infty$ function contributing only on $\mathbb R^k \setminus \hat{\mathcal K}$, i.e.
	\[
	\rho(x) = 
	\begin{cases}
		1 , & x \in \mathbb R^k \setminus \hat{\mathcal K} , \\
		0 , & x \in \mathcal K ,
	\end{cases}
	\]
	and $\rho(x) \in (0,1)$ for $x \in \hat{\mathcal K} \setminus \mathcal K$. Due to compactness, such a function can always be constructed using local $C^\infty$ bump functions and a partition of unity. It is important to note that the maps \eqref{eq:partial_fenichel_normal_form} and \eqref{eq:partial_fenichel_normal_form_enlarged} agree on $\mathcal K$.
	
	\
	
	Finally, it will be helpful to have a notation for spectral bounds. We define
	\begin{equation}
		\label{eq:spectral_bounds}
		\nu_A := \sup_{x \in \hat{\mathcal K} , j = 1, \ldots , n_a} | 
		\mu_{a,j}(x) | < 1 , \qquad 
		\nu_R := \inf_{x \in \hat{\mathcal K} , j = 1, \ldots , n_r} | 
		\mu_{r,j}(x) | > 1 ,
	\end{equation}
	where $\mu_{a,j}(x)$ and $\mu_{r,j}(x)$ denote multipliers of the matrices $I_{n_a} + \Lambda_a|_{\hat S}$ and $I_{n_r} + \Lambda_r|_{\hat S}$ respectively. Note that $\nu_A$ and $\nu_R$ agree with their previous definition in \eqref{eq:spectral_bounds_original} since the multipliers $\mu_{a,j}$ and $\mu_{r,j}$ are invariant under $C^{r}$ coordinate transformations.

	\subsection{Checking necessary and sufficient conditions}
	\label{sub:NS_conds}

	We now show that the relevant invariance properties for the application of results in \cite{Nipp2013} are satisfied for the map \eqref{eq:partial_fenichel_normal_form_enlarged}. This will be easier with a little extra notation. We first rewrite the map \eqref{eq:partial_fenichel_normal_form_enlarged} as
	\begin{equation}
		\label{eq:map_bxa}
		\begin{split}
			b & \mapsto \bar b = \widehat H^b(b,x,a,\eps) , \\
			x & \mapsto \bar x = \widehat H^x(b,x,a,\eps) , \\
			a & \mapsto \bar a = \widehat H^a(b,x,a,\eps) , 
		\end{split}
	\end{equation}
	where in particular, we have permuted the order of the equations and the arguments in order for simplicity in the notation to follow. Our analysis here is purely local, and may be restricted to
	\[
	b \in B := \{b \in \mathbb R^{n_r} : |b| \leq \delta_0 \} , \quad a \in A := \{a \in \mathbb R^{n_a} : |a| \leq \delta_0 \} , \quad \eps \in [0, \eps_0) ,
	\]
	for fixed $\delta_0, \eps_0 > 0$ which we shall frequently choose to be sufficiently small for the validity of estimates. We are also free to choose $\mathcal K$ and $\hat{\mathcal K}$ in such a way that
	\[
	x \in \mathcal K \ \implies \ |x| < \Delta , \qquad
	x \in \hat{\mathcal K} \ \implies \ |x| < \hat \Delta ,
	\]
	for fixed $\hat \Delta > \Delta > 0$.
	
	Following \cite{Nipp2013}, we now introduce two more equivalent formulations that are better suited for the analyses of attracting and repelling invariant objects respectively. Consider the attracting formulation first. In this case, we write
	\[
	x_A := (b,x)^\textnormal{T} \in X_A, \qquad y_A = a \in Y_A ,
	\]
	define
	\[
	\widehat H_A^{x_A}(x_A,y_A,\eps) := 
	\begin{pmatrix}
		\widehat H^b(x_A,y_A,\eps) \\
		\widehat H^x(x_A,y_A,\eps)
	\end{pmatrix} ,
	\qquad 
	\widehat H^{y_A}_A(x_A,y_A,\eps) := \widehat H^a(x_A,y_A,\eps) ,
	\]
	and consider the map \eqref{eq:map_bxa} expressed as
	\begin{equation}
		\label{eq:attracting_formulation}
		\widehat H_A :
		\begin{pmatrix}
			x_A \\ 
			y_A 
		\end{pmatrix}
		\mapsto 
		\begin{pmatrix}
			\bar x_A \\
			\bar y_A
		\end{pmatrix}
		=
		\begin{pmatrix}
			\widehat H_A^{x_A}(x_A,y_A,\eps) \\
			\widehat H_A^{y_A}(x_A,y_A,\eps)
		\end{pmatrix} .
	\end{equation}
	For the repelling formulation, we write
	\[
	x_R := b \in X_R , \qquad y_R := (x,a)^\textnormal{T} \in Y_R ,
	\]
	define
	\[
	\widehat H_R^{x_R}(x_R,y_R,\eps) :=  \widehat H^b(x_R,y_R,\eps) , \qquad
	\widehat H_R^{y_R}(x_R,y_R,\eps) :=
	\begin{pmatrix}
		\widehat H^x(x_R,y_R,\eps) \\
		\widehat H^a(x_R,y_R,\eps)
	\end{pmatrix} ,
	\]
	and consider the map \eqref{eq:map_bxa} expressed as
	\begin{equation}
		\label{eq:repelling_formulation}
		\widehat H_R :
		\begin{pmatrix}
			x_R \\ 
			y_R 
		\end{pmatrix}
		\mapsto 
		\begin{pmatrix}
			\bar x_R \\
			\bar y_R
		\end{pmatrix}
		=
		\begin{pmatrix}
			\widehat H_R^{x_R}(x_R,y_R,\eps) \\
			\widehat H_R^{y_R}(x_R,y_R,\eps)
		\end{pmatrix} .
	\end{equation}

	\subsubsection{Invariance conditions}
	
	The relevant inflowing, outflowing and invariance conditions 
	are summarized in the following result.

	\begin{lemma}
		\label{lem:inflowing_outflowing_conds}
		\label{lem:HM_ab}
		\textup{(c.f.~\cite[Hypothesis HM]{Nipp2013})} For each
		\begin{equation}
			\label{eq:hat_delta_bounds}
			0 < \hat \delta < \min \left\{ \frac{\hat \Delta - (\Delta + c_0)}{C} ,  \frac{\Delta^2}{\hat \Delta (1 + M \hat \Delta)} \right\} ,
		\end{equation}
		where $M := \sup_{x \in \hat{\mathcal K}} |x \hat \rho'(x)|$, there exists $\delta_0, \eps_0 > 0$ such that the map \eqref{eq:map_bxa} has the following invariance properties for all $\eps \in [0,\eps_0)$:
		\begin{enumerate}
			\item $\widehat H_A$ is inflowing with respect to $Y_A$, i.e.
			\[
			(x_A,y_A) \in X_A \times Y_A
			\quad \implies \quad 
			\widehat H^{y_A}_A(x_A,y_A,\eps) \in Y_A ,
			\]
			and outflowing with respect to $X_A$, i.e.~for all $(\bar x_A,y_A) \in X_A \times Y_A$ there exists an $x_A \in X_A$ such that
			\[
			\bar x_A = \widehat H_A^{x_A}(x_A,y_A,\eps) .
			\]
			\item $\widehat H_R$ is inflowing with respect to $Y_R$, i.e.
			\[
			(x_R,y_R) \in X_R \times Y_R
			\quad \implies \quad 
			\widehat H^{y_R}_A(x_R,y_R,\eps) \in Y_R ,
			\]
			and outflowing with respect to $X_R$, i.e.~for all $(\bar x_R,y_R) \in X_R \times Y_R$ there exists an $x_R \in X_R$ such that
			\[
			\bar x_R = \widehat H_R^{x_R}(x_R,y_R,\eps) .
			\]
		\end{enumerate}
	\end{lemma}
	
	\begin{proof}
		We start with the statement (a). 
		Since $\widehat H_A$ is smooth in all arguments and we have the linear contractivity property
		\[
		\sup_{x \in \mathcal K} \norm{ I_{n_a} + \Lambda_a|_{\hat S} } = \nu_A < 1 ,
		\]
		we obtain the estimate
		\begin{equation}
			\label{eq:Ha_bound}
				| \widehat H_A^{x_A}(x_A,y_A,\eps) | = | \widehat H^a(x_A,y_R,\eps) |
				\leq \nu_A \delta_0 + \delta_0^2 c_1 + \eps_0 c_2 < \delta_0 \\
		\end{equation}
		for all $(x_A,y_A) \in X_A \times Y_A$, for some constants $c_1, c_2 \geq 0$, and where the rightmost inequality is obtained by choosing $\delta_0, \eps_0 > 0$ sufficiently small. Hence $\widehat H_A$ is inflowing with respect to $Y_A$.

		
		In order to show that $\widehat H_A$ is outflowing with respect to $X_A$, let $(\bar x_A, y_A) = ((\bar b, \bar x), a) \in X_A \times Y_A$ so that in particular, $|\bar x| \leq \hat \Delta$ and $|\bar b| < \delta_0$. We need to show that there exists $(x_A,y_A) = ((b,x)^\textnormal{T},a) \in X_A \times Y_A$ such that
		\[
		\bar x = \widehat H^x(b,x,a,\eps)  \qquad \text{and} \qquad \bar b = \widehat H^b(b,x,a,\eps) .
		\]
		The latter follows immediately, since the linear expansion property
		\[
		\inf_{x \in \hat{\mathcal K}} \norm{ I_{n_r} + \Lambda_r|_{\hat S}} \geq \nu_R > 1
		\] 
		implies the existence of an inverse $(I_{n_r} + \Lambda_r(x,a,b))^{-1}$ for $\delta_0 > 0$ sufficiently small. In particular, $\bar b = \widehat H^b(b,x,a,\eps)$ can be solved via the implicit function theorem and the norm of $b$ can be estimated directly via the implicit equation
		\begin{equation}
			\label{eq:Hb_soln}
			b = (I_{n_r} + \Lambda_r(x,a,b))^{-1} \bar b + O(\eps) ,
		\end{equation}
		which satisfies $|b| < \delta_0$ as required for sufficiently small $\delta_0, \eps_0 > 0$.
		
		Now consider $\bar x = \widehat H^x(b,x,a,\eps)$. We need to show that this equation
		can be solved for an $x$ such that $x \in \hat{\mathcal K}$. Towards this end, we define
		\[
		\hat \rho(x) := \rho(x) n_x \frac{|x|}{x} 
		\]
		and consider
		\[
		x = A(x)^{-1} \left(\bar x - \tilde N^x(x,N^yu) N^y (a, b)^\textnormal{T} - \eps \tilde 	G^x(x,N^y(a,b)^\textnormal{T},\eps) \right) =: R(x) ,
		\]
		where
		\[
		A(x)^{-1} := \left(1 +  \hat \delta \frac{\hat \rho(x)}{|x|} \right)^{-1} .
		\]
		Fixed points $R(x) = x$ correspond to solutions of $\bar x = \widehat H^x(b,x,a,\eps)$, and may be identified by a contraction mapping argument if we can show that $R$ is a map $R: \hat{\mathcal K} \to \hat{\mathcal K}$ satisfying $|R'(x)| < 1$ for all $x \in \hat{\mathcal K}$. First, we show that $R: \hat{\mathcal K} \to \hat{\mathcal K}$. Observe that
		\[
		|R(x)| \leq \frac{1}{1 + \hat \delta \hat \rho(x) / |x|} \left(|\bar x| + c_1 \delta_0 + c_2 \eps_0 \right) ,
		\]
		for some (new) constants $c_1, c_2 \geq 0$. We need to show that $|R(x)| < \hat \Delta$ for all $\bar x$ and $x$ such that $|x|,|\bar x| < \hat \Delta$. There are two cases to consider:
		\begin{itemize}
			\item $|x| < \Delta + c_0 < \hat \Delta$ for some $c_0 > 0$, and
			\item $|x| \in [\Delta + c_0, \hat \Delta)$.
		\end{itemize}
		In the first case, solutions $\bar x$ (if they do exist) have to satisfy
		\[
		|\bar x| = |\widehat H^x(b,x,a,\eps) | < \Delta + c_0 + c_1 \delta_0 + c_2 \eps_0 + \hat \delta C ,
		\]
		where $C > 0$ is is a constant in the interval $(\sup_{|x|<\Delta+c_0} \rho(x),1)$. Hence
		\[
		|R(x)| \leq \Delta + c_0 + 2 c_1 \delta_0 + 2 c_2 \eps_0 + \hat \delta C < \hat \Delta ,
		\]
		is satisfied after choosing
		\begin{equation}
			\label{eq:delta_upper_bound}
			\hat \delta < \frac{\hat \Delta - (\Delta + c_0)}{C} - 2 \frac{c_1 \delta_0 + c_2 \eps_0}{C} .
		\end{equation}
		Now assume that $|x| \in [\Delta + c_0, \hat \Delta)$. We use the fact that $\rho(x) \geq C > 0$ uniformly with respect to all such $x$. It follows that
		\[
			|R(x)| \leq \frac{1}{1 + \hat \delta \hat \rho(x) / |x|} \left(\hat \Delta + c_1 \delta_0 + c_2 \eps_0 \right) 
			\leq \frac{1}{1 + \hat \delta C / \hat \Delta} \left(\hat \Delta + c_1 \delta_0 + c_2 \eps_0 \right) 
			< \hat \Delta
		\]
		after choosing
		\begin{equation}
			\label{eq:delta_lower_bound}
			\hat \delta > \frac{c_1 \delta_0 + c_2 \eps_0}{C} .
		\end{equation}
		Combining \eqref{eq:delta_upper_bound} and \eqref{eq:delta_lower_bound} shows that $R:\hat{\mathcal K} \to \hat{\mathcal K}$ for all $\hat \delta$ such that
		\[
		\frac{c_1 \delta_0 + c_2 \eps_0}{C} < \hat \delta_0 < \frac{\hat \Delta - (\Delta + c_0)}{C} - 2 \frac{c_1 \delta_0 + c_2 \eps_0}{C} ,
		\]
		which, by choosing sufficiently small $\delta_0, \eps_0 > 0$, can satisfied for any fixed
		\begin{equation}
			\label{eq:hat_delta_bounds1}
			0 < \hat \delta < \frac{\hat \Delta - (\Delta + c_0)}{C} .
		\end{equation}
		It remains to show that $R$ is a contraction. A direct calculation gives
		\[
		\begin{split}
			|R'(x)| &\leq A(x)^{-1} (\tilde c_1 \delta_0 + \tilde c_2 \eps_0) +
			(A(x)^{-1})' \left( |\bar x| + c_1 \delta_0 + c_2 \eps_0 \right) \\
			&\leq \tilde c_1 \delta_0 + \tilde c_2 \eps_0 + \left( |\bar x| + c_2 \delta_0 + c_2 \eps_0 \right) \frac{\hat \delta |x \hat \rho'(x) - \rho(x) | }{|x|^2} ,
		\end{split}
		\]
		where we used $A(x) \geq 1$. Notice that the right-hand-side is well-defined for $x = 0$, since $0 \in \mathcal K$ and $\hat \rho'|_{\mathcal K} = \hat \rho|_{\mathcal K} \equiv 0$. In particular,
		\[
		\frac{|x \hat \rho'(x) - \rho(x) | }{|x|^2} \leq \frac{1}{\Delta^2} \sup_{|x| \in [\Delta, \hat \Delta)} |x \hat \rho'(x) - \rho(x) | \leq \frac{\hat \Delta M + 1}{\Delta^2} ,
		\]
		where $M := \sup_{x \in \hat{\mathcal K}} |x \hat \rho'(x)|$. Hence
		\[
		|R'(x)| \leq \tilde c_1 \delta_0 + \tilde c_2 \eps_0 + \hat \delta \left( |\bar x| + c_2 \delta_0 + c_2 \eps_0 \right) \left(\frac{\hat \Delta M + 1}{\Delta^2}\right) < 1 ,
		\]
		where the rightmost inequality is satisfied for sufficiently small $\delta_0 , \eps_0 > 0$ and
		\begin{equation}
			\label{eq:hat_delta_bounds2}
			\hat \delta < \frac{\Delta^2}{\hat \Delta (1 + M \hat \Delta)} .
		\end{equation}
		It follows that $R$ is a contraction as required. Finally, combining \eqref{eq:hat_delta_bounds1} and \eqref{eq:hat_delta_bounds2} yields the bounds for $\hat \delta$ in \eqref{eq:hat_delta_bounds}, completing the proof of statement (a).
		
		\
		
		The proof of statement (b) follows from arguments similar to those given in the proof of (a) above. Specifically, $\widehat H_R$ is inflowing with respect to $Y_R$ since
		\begin{itemize}
			\item $(x_R, y_R) \in X_R \times Y_R \ \implies \ |\widehat H^x(x_R,y_R,\eps)| < \hat \Delta$;
			\item $(x_R, y_R) \in X_R \times Y_R \ \implies \ |\widehat H^a(x_R,y_R,\eps)| < \delta_0$.
		\end{itemize}
		The former inequality was shown in the proof that $R:\hat{\mathcal K} \to \hat{\mathcal K}$ above, and the latter follows from \eqref{eq:Ha_bound}. The fact that $\widehat H_R$ is outflowing with respect to $X_R$ follows immediately from (in fact, is the same as) the fact that by \eqref{eq:Hb_soln}, $\hat b = \widehat H^b(b,x,a,\eps)$ has a unique solution with $|\bar b| < \delta_0$.
		%
		%
	\end{proof}
	
	\begin{remark}
		Lemma \ref{lem:HM_ab} confirms the inflowing/outflowing hypotheses HMa)-b) in \cite{Nipp2013} for the map \eqref{eq:map_bxa}, i.e.~for the map \eqref{eq:partial_fenichel_normal_form_enlarged}. Note the crucial role played by the cutoff function $\hat \delta \rho(x) n_x$. The fact that $\hat \delta$ must be fixed and positive in accordance with \eqref{eq:hat_delta_bounds}, or more precisely, in accordance with \eqref{eq:delta_lower_bound}, shows that such cutoff procedures are necessary to satisfy the hypotheses.
	\end{remark}

	\subsubsection{Conditions on relative contraction rates}

	We turn now to the verification of conditions on the relative contraction rates. Such conditions are typically given in a quite general setting in terms of Lipschitz-type bounds for the components of the map \eqref{eq:map_bxa}. Our aim here is to restate these conditions in terms of spectral properties like normal hyperbolicity which, as a spectral condition, is invariant under $C^1-$smooth coordinate transformations and thus detectable in the original coordinates of \eqref{eq:main_map}. The cost of this reformulation is that we require at least $C^1-$smoothness of the map, while only Lipschitz continuity is required in \cite{Nipp2013}. A similar price is paid in Fenichel theory \cite{Fenichel1979,Wiggins2013}.
	
	Let $i \in \{A,R\}$ in order to streamline sub/superscript notations wherever possible. Since the maps $\widehat H_i$ are $C^1-$smooth in $(x_i,y_i,\eps)$, there exist constants $\Gamma^i_{11}, L^i_{12}, L^i_{13}, L^i_{21}, L^i_{22}, L^i_{23} \geq 0$ such that the following Lipschitz-type estimates are satisfied for any fixed choice of $x_{i,1}, x_{i,2} \in X_i$, $y_{i,1}, y_{i,2} \in Y_i$ and $\eps_1, \eps_2 \in [0,\eps_0)$:
	\begin{equation}
		\label{eq:Lipschitz_consts}
		\begin{split}
			|\widehat H^{x_i}_i(x_{i,1},y_i,\eps) - \widehat H^{x_i}_i(x_{i,2},y_i,\eps)| & \geq \Gamma^i_{11} |x_{i,1} - x_{i,2}| , \\
			|\widehat H^{x_i}_i(x_i,y_{i,1},\eps_1) - \widehat H^{x_i}_i(x_i,y_{i,2},\eps_2)| & \leq L^i_{12} |y_{i,1} - y_{i,2}| + L^i_{13} |\eps_1 - \eps_2| , \\
			|\widehat H^{y_i}_i(x_{i,1},y_{i,1},\eps_1) - \widehat H^{y_i}_i(x_{i,2},y_{i,2},\eps_2)| & \leq L^i_{21} |x_{i,1} - x_{i,2}| + L^i_{22} |y_{i,1} - y_{i,2}| + L^i_{23} |\eps_1 - \eps_2| .
		\end{split}
	\end{equation}
	Smoothness allows for a characterisation of these Lipschitz-type bounds in terms of partial derivatives. 
	
	\begin{lemma}
		\label{lem:constants}
		The following choices satisfy the conditions in \eqref{eq:Lipschitz_consts} for the map $\widehat H_A$:
		\begin{equation}
			\label{eq:Lipschitz_bounds_A}
			\begin{aligned}
				\Gamma_{11}^A &= 1 + r^A_{11}(\delta_0,\eps_0) , \\
				L_{12}^A &= \sup_{x \in \hat{\mathcal K}} \norm{ (\tilde N^x N^y_a)|_{\hat S} } + r_{12}^A(\delta_0,\eps_0) , \\
				L_{13}^A &= \sup_{x \in \hat{\mathcal K}} \sqrt{ |\tilde G^x|_{\hat S}|^2 + | ((N_r^y)^{-1} DfG)|_{\hat S}|^2 } + r_{13}^A(\delta_0, \eps_0) , \\
				L_{21}^A &= r_{21}^A(\delta_0,\eps_0) , \\
				L_{22}^A &= \nu_A + r_{22}^A(\delta_0, \eps_0) , \\
				L_{23}^A &= \sup_{x \in \hat{\mathcal K}} |((N^y_a)^{-1} DfG)|_{\hat S}| + r_{23}^A(\eps_0) ,
			\end{aligned}
		\end{equation}
		where the functions $r_{kj}^A(\delta_0,\eps_0)$ and $r_{23}^A(\eps_0)$ are continuous such that $r_{kj}^A(0,0) = r_{23}^A(0) = 0$.
		
		The following choices satisfy the conditions in \eqref{eq:Lipschitz_consts} for the map $\widehat H_R$:
		\begin{equation}
			\label{eq:Lipschitz_bounds_R}
			\begin{aligned}
				\Gamma_{11}^R &= \nu_R + r_{11}^R(\delta_0, \eps_0) , \\
				L_{12}^R &= r_{12}^R(\delta_0,\eps_0) , \\
				L_{13}^R &= \sup_{x \in \hat{\mathcal K}} |((N^y_r)^{-1} DfG)|_{\hat S}| + r_{13}^R(\eps_0) , \\
				L_{21}^R &= \sup_{x \in \hat{\mathcal K}} \norm{(\tilde N^x N^y_r)|_{\hat S}} + r_{21}^R(\delta_0,\eps_0) , \\
				L_{22}^R &= 1 + \hat \delta + r_{22}^R(\delta_0,\eps_0) , \\
				L_{23}^R &= \sup_{x \in \hat{\mathcal K}} \sqrt{ | \tilde G^x|_{\hat S}|^2 + | ((N_a^y)^{-1} DfG)|_{\hat S}|^2 } + r_{23}^R(\delta_0, \eps_0) , 
			\end{aligned}
		\end{equation}
		where the functions $r^R_{kj}(\delta_0,\eps_0)$ and $r^R_{13}(\eps_0)$ are continuous such that $r_{kj}^R(0,0) = r_{13}^R(0) = 0$.
	\end{lemma}
	
	\begin{proof}
		Since both $\widehat H_i$ are at least $C^1-$smooth in all arguments, partial derivatives are locally bounded and the following choices for $\Gamma^i_{11}, L^i_{12}, L^i_{13}, L^i_{21}, L^i_{22}$, $L^i_{23} \geq 0$ satisfy \eqref{eq:Lipschitz_consts}:
		\begin{equation}
			\label{eq:Lipschitz_consts_2}
			\begin{split}
				\Gamma_{11}^i = \inf_{({x_i},{y_i},\eps) \in U_i} \norm{ D_{x_i} \widehat H^{x_i}_i({x_i},{y_i},\eps) } , \qquad 
				& L_{21}^i = \sup_{({x_i},{y_i},\eps) \in U_i} \norm{ D_{x_i} \widehat H^{y_i}_i({x_i},{y_i},\eps) } , \\
				L_{12}^i = \sup_{({x_i},{y_i},\eps) \in U_i} \norm{ D_{y_i} \widehat H^{x_i}_i({x_i},{y_i},\eps) } , \qquad 
				& L_{22}^i = \sup_{({x_i},{y_i},\eps) \in U_i} \norm{ D_{y_i} \widehat H^{y_i}_i({x_i},{y_i},\eps)} , \\
				L_{13}^i = \sup_{({x_i},{y_i},\eps) \in U_i} \norm{ D_\eps \widehat H^{x_i}_i({x_i},{y_i},\eps)} , \qquad 
				& L_{23}^i = \sup_{({x_i},{y_i},\eps) \in U_i} \norm{ D_\eps \widehat H^{y_i}_i({x_i},{y_i},\eps)} ,
			\end{split}
		\end{equation}
		where $U_i := X_i \times Y_i \times [0,\eps_0)$. We omit the subsequent derivations for all expressions except $\Gamma_{11}^A$ and $L_{22}^R$, since these can be derived by a direct evaluation of the corresponding expression in \eqref{eq:Lipschitz_consts_2}, followed (if necessary) by the application of some standard triangle-type inequalities.
		
		To obtain the expression for $\Gamma_{11}^A$, notice that Taylor expanding about $\hat{S} : a = b = \eps = 0$ gives
		\[
		\begin{split}
			\Gamma_{11}^A &= 
			\inf_{(x_A,y_A,\eps) \in U_A} \norm{
				\begin{pmatrix}
					D_b \widehat H^b(b,x,a,\eps) & D_x \widehat H^b(b,x,a,\eps) \\
					D_b \widehat H^x(b,x,a,\eps) & D_x \widehat H^x(b,x,a,\eps)
				\end{pmatrix}
			} \\
			&= \inf_{x \in \hat{\mathcal K}} \norm{
				\begin{pmatrix}
					D_b \widehat H^b & D_x \widehat H^b \\
					D_b \widehat H^x & D_x \widehat H^x
				\end{pmatrix}
				\bigg|_{\hat S}
			} + r_{11}^A(\delta_0, \eps_0) \\
			&= \inf_{x \in \hat{\mathcal K}} \norm{
				\begin{pmatrix}
					I_{n_r} + \Lambda_r & O_{n_r \times k} \\
					\tilde N^x N^y_r & I_k + \hat \delta D_x (\rho(x) n_x)
				\end{pmatrix}
				\bigg|_{\hat S}
			} + r_{11}^A(\delta_0, \eps_0) ,
		\end{split}
		\]
		where $r_{11}^A(\delta_0, \eps_0)$ is continuous and satisfies $r_{11}^A(0,0) = 0$. Since $\norm{I_k} = 1$, the matrix norm in the last line must be greater than or equal to $1$. Hence $\Gamma_{11}^A \geq 1 + r_{11}^A(\delta_0,\eps_0)$, thereby justifying the choice for $\Gamma_{11}^A$ in \eqref{eq:Lipschitz_bounds_A}.
		
		Now consider $L_{22}^R$. Similar calculations lead to
		\[
		L_{22}^R = \sup_{x \in \hat{\mathcal K}} \norm{
			\begin{pmatrix}
				I_k + \hat \delta D_x (\rho(x) n_x) & (\tilde N^x N^y_a) (x,0,0) \\
				O_{n_a \times k} & I_{n_a} + \Lambda_a(x,0,0)
			\end{pmatrix}
		} + r_{22}^R(\delta_0, \eps_0) =: \norm{B} + \tilde r_{22}^R(\delta_0, \eps_0) ,
		\]
		for a continuous function $\tilde r_{22}^R(\delta_0, \eps_0)$ such that $\tilde r_{22}^R(0,0) = 0$. We need a sufficiently sharp estimate for $\norm{B}$ which, by construction, has an associated operator norm property
		\begin{equation}
			\label{eq:op_norm_cond}
			\norm{B
				\begin{pmatrix}
					x \\
					b
				\end{pmatrix}
			}
			\leq L^R_{22}
			\left\lvert
			\begin{pmatrix}
				x \\
				b
			\end{pmatrix}
			\right\rvert
		\end{equation}
		for sufficiently small but fixed $\delta_0, \eps_0 > 0$, uniformly with respect to $(b,x) \in Y_R$. Direct estimates yield
		\[
		\begin{split}
			\norm{B
				\begin{pmatrix}
					x \\
					b
				\end{pmatrix}
			}^2 &= \left\lvert
			\begin{pmatrix}
				(I_k + \hat \delta D_x (\rho(x) n_x))x + (\tilde N^x N^y_a) (x,0,0) b \\
				(I_{n_a} + \Lambda_a(x,0,0)) b
			\end{pmatrix}
			\right\rvert^2  \\
			& \leq
			(1 + \hat \delta)^2 \hat \Delta^2 + c_1 \hat \Delta \delta_0 + c_2 \delta_0^2 + \nu_A^2 \delta_0^2 \\
			&\leq \tilde C^2 \left\lvert
			\begin{pmatrix}
				x \\
				b
			\end{pmatrix}
			\right\rvert^2 \\
			&\leq \tilde C^2 (\hat \Delta^2 + \delta_0^2) ,
		\end{split}
		\]
		where $c_1, c_2 \geq 0$ are constants and we require that
		\[
		\tilde C \geq \sqrt{ \frac{(1 + \hat \delta)^2 \hat \Delta^2 + c_1 \hat \Delta \delta_0 + c_2 \delta_0^2 + \nu_A^2 \delta_0^2}{\hat \Delta ^2 + \delta_0^2} } = 
		1 + \hat \delta + \hat r_{22}^R(\delta_0) ,
		\]
		where the function $\hat r_{22}^R(\delta_0)$ is continuous with $\hat r_{22}^R(0) = 0$. Choosing the minimal such $\tilde C$ and combining with \eqref{eq:op_norm_cond} leads to the expression for $L_{22}^R$ in \eqref{eq:Lipschitz_bounds_R}.
	\end{proof}
	
	We now use the expressions in Lemma \ref{lem:constants} in order to check the relevant contractivity/repulsivity conditions. 
	In order to simplify notation in the following we define
	\[
	\omega_i := \frac{2 L^i_{12} L^i_{21}}{\Gamma^i_{11} - L^i_{22} + \sqrt{(\Gamma^i_{11} - L^i_{22})^2 - 4 L^i_{12} L^i_{21} }}
	\]
	for both $i \in \{A,R\}$.
	
	\begin{lemma}
		\label{lem:CM_ar}
		\textup{(c.f.~\cite{Nipp2013} conditions CM, CMA and CMR)}
		For sufficiently small $\delta_0, \eps_0 > 0$ and $\hat \delta$ such that
		\begin{equation}
			\label{eq:hat_delta_bounds3}
			\nu_R > 1 + \hat \delta ,
		\end{equation}
		the following conditions are satisfied:
		\begin{itemize}
			\item[(i)] In both cases $i \in \{A,R\}$, we have
			\[
			2 \sqrt{L^i_{12} L^i_{21}} < \Gamma^i_{11} - L_{22}^i .
			\]
			\item[(ii)] In case $i = A$ we have
			\[
			L_{22}^A + \omega_A < 1.
			\]
			\item[(iii)] In case $i = R$ we have
			\[
			\Gamma^R_{11} - \omega_R > 1 .
			\]
		\end{itemize}
	\end{lemma}
	
	\begin{proof}
		This is immediate given the expressions in Lemma \ref{lem:constants}. Since $L^i_{12} L^i_{21} \to 0$ as $(\delta_0,\eps_0) \to (0,0)$, the condition (i) can always be satisfied as long as $0 < \Gamma^i_{11} - L_{22}^i$. The latter conditions can be checked directly. One the one hand we have
		\[
		\Gamma^A_{11} - L_{22}^A = 1 - \nu_A - (r_{11}^A(\delta_0,\eps_0) + r_{22}^A(\delta_0,\eps_0)) > 0 ,
		\]
		for sufficiently small $\delta_0, \eps_0 > 0$ since $\nu_A < 1$. On the other hand we have
		\[
		\Gamma^R_{11} - L_{22}^R = \nu_R - (1 + \hat \delta) - (r_{11}^R(\delta_0,\eps_0) + r_{22}^R(\delta_0,\eps_0)) > 0 ,
		\]
		for sufficiently small $\delta_0, \eps_0 > 0$ and $\hat \delta$ in the interval \eqref{eq:hat_delta_bounds3}. In order to verify (ii)-(iii), observe that $\omega_i \to 0$ as $(\delta_0,\eps_0) \to (0,0)$ in both cases $i \in \{A,R\}$. The results in (ii) and (iii) follow immediately since $L_{22}^A < 1$ and $\Gamma^R_{11} > 1$ for sufficiently small $\delta_0, \eps_0 > 0$.
	\end{proof}

	\subsubsection{Smoothness conditions}
	
	It remains to check a number of conditions relating to smoothness of the invariant slow manifolds and foliations. We start with conditions for the smoothness of slow manifolds.

	\begin{lemma}
		\label{lem:conds_smoothness_slow_manifolds}
		\textup{(c.f.~\cite{Nipp2013} conditions CMR(k) and CMA(k) for smoothness of slow manifolds)} 
		Consider the $C^{r\geq1}$-smooth map \eqref{eq:map_bxa}. Then for sufficiently small $\delta_0, \eps_0 > 0$, and $\hat \delta$ such that
		\begin{equation}
			\label{eq:hat_delta_bounds4}
			\nu_R > (1 + \hat \delta)^r > 1 ,
		\end{equation}
		the following conditions are satisfied:
		\begin{itemize}
			\item[(i)] In case $i = A$ we have
			\[
			L_{22}^A + \omega_A < (\Gamma_{11}^A - \omega_A)^r .
			\]
			\item[(ii)] In case $i = R$ we have
			\[
			(L_{22}^R + \omega_R)^r < \Gamma^R_{11} - \omega_R .
			\]
		\end{itemize}
	\end{lemma}
	
	\begin{proof}
		Straightforward calculations using 
		the expressions in Lemma \ref{lem:constants}. The requirement \eqref{eq:hat_delta_bounds4} is necessary in (ii).
	\end{proof}

	Finally we need to verify a number of conditions relating to smoothness of the invariant foliations corresponding to persisting stable and unstable manifolds. These can be checked directly using the expressions in Lemma \ref{lem:constants} as well as the quantities
	\[
	L_{11}^A := \sup_{({x_A},{y_A},\eps) \in U_A} \norm{ D_{x_A} \widehat H^{x_A}_A({x_A},{y_A},\eps) } , \qquad 
	\Gamma_{22}^i := \inf_{({x_i},{y_i},\eps) \in U_i} \norm{ D_{y_i} \widehat H^{y_i}_i({x_i},{y_i},\eps) } .
	\]
	For $L_{11}^A$ and $\Gamma_{22}^R$, arguments based on direct estimates and operator norm properties similar to those applied for $\Gamma_{11}^A$ in the proof of Lemma \ref{lem:constants} lead to concrete (though potentially less than optimal) expressions
	\begin{equation}
		\label{eq:L11}
		L_{11}^A = 1 + \hat \delta + \tilde r_{11}^A(\delta_0, \eps_0) , \qquad 
		\Gamma_{22}^R = 1 + \tilde r_{22}^R(\delta_0, \eps_0) ,
	\end{equation}
	where $\tilde r^A_{11}(\delta_0, \eps_0)$ and $\tilde r^R_{22}(\delta_0, \eps_0)$ are continuous such that $\tilde r^A_{11}(0,0) = \tilde r^R_{22}(0,0) = 0$. For $\Gamma_{22}^R$, direct calculations yield
	\begin{equation}
		\label{eq:Gamma22}
		\Gamma_{22}^A = 
		\inf_{x \in \hat{\mathcal K}} \norm{ I_{n_a} + \Lambda_a|_{\hat S} } + \tilde r^A_{22}(\delta_0, \eps_0) ,
	\end{equation}
	where $\tilde r^A_{22}(\delta_0, \eps_0)$ is continuous such that $\tilde r^A_{22}(0,0) = 0$. It is important to note that $\Gamma_{22}^A = 0$ is possible under superstable conditions in which all $n_a$ multipliers of $I_{n_a} + \Lambda_a|_{\hat S}$ have real part zero.

	\begin{lemma}
		\label{lem:conds_smoothness_foliation}
		\textup{(c.f.~\cite{Nipp2013} conditions CMB, CMAB($k-1$) and CMRB($k-1$) for smoothness of the foliations)} 
		Consider the $C^{r\geq1}$-smooth map \eqref{eq:map_bxa}, and assume that $I_{n_a} + \Lambda_a|_{\hat S}$ has at least one multiplier $\mu_{a,j}(x)$ such that
		\begin{equation}
			\label{eq:nondegeneracy_fiber_smoothness}
			|\mu_{a,j}(x)| > \tilde c > 0 
		\end{equation}
		for some $x \in \hat{\mathcal K}$, for some constant $\tilde c > 0$. Then for sufficiently small $\delta_0, \eps_0 > 0$, and $\hat \delta$ such that
		\begin{equation}
			\label{eq:hat_delta_bounds5}
			\nu_A < \frac{1}{(1 + \hat \delta)^{r-1}} < 1 + \hat \delta < \nu_R ,
		\end{equation}
		the following conditions are satisfied:
		\begin{itemize}
			\item[(i)] In both cases $i \in \{A,R\}$ we have
			\[
			\Gamma_{22}^i - \omega_i > 0.
			\]
			\item[(ii)] In case $i = A$ we have
			\[
			\frac{\Gamma^A_{11} - \omega_A}{L_{22}^A + \omega_A} > (L_{11}^A + \omega_A)^{r-1} .
			\]
			\item[(iii)] In case $i = R$ we have
			\[
			\frac{L^R_{22} + \omega_R}{\Gamma_{11}^R - \omega_R} < (\Gamma_{22}^R - \omega_R)^{r-1} .
			\]
		\end{itemize}
	\end{lemma}
	
	\begin{proof}
		The condition (i) is immediate in case $i = R$. In case $i = A$, it follows from the requirement $|\mu_{a,j}(x)| > \tilde c > 0$ for some $j \in \{1, \ldots, n_a\}$ and $x \in \hat{\mathcal K}$ since in this case,
		\[
		\norm{ I_{n_a} + \Lambda_a|_{\hat S} } \geq |\mu_{a,j}(x)| > \tilde c > 0 .
		\]
		
		Conditions (ii)-(iii) are verified directly using equations \eqref{eq:L11}-\eqref{eq:Gamma22} and the expressions in Lemma \ref{lem:constants}. Specifically, the inequality in (ii) follows from the fact that
		\[
		\frac{\Gamma_{11}^A}{L_{22}^A} \sim \frac{1}{\nu_A} > (1 + \hat \delta)^{r-1} 
		\]
		under \eqref{eq:hat_delta_bounds5} as $(\delta_0, \eps_0) \to (0,0)$, and the inequality in (iii) follows from the fact that
		\[
		\frac{L_{22}^R}{\Gamma_{11}^R} \sim \frac{1 + \hat \delta}{\nu_R} <
		1 
		\]
		as $(\delta_0, \eps_0) \to (0,0)$.
	\end{proof}
	
	\begin{remark}
		\label{rem:foliaiton_smoothness_local}
		The local condition on the stable multipliers in \eqref{eq:nondegeneracy_fiber_smoothness} rules out the possibility of degenerate local superstability in the sense that we disallow the case where $\mu_{a,j}(x) \equiv 0$ for all $x \in \hat{\mathcal K}$ and $j = 1, \ldots, n_a$.
	\end{remark}


	\subsection{Proof of the main results}
	\label{sub:NS_thms}
	
	We are now in a position to prove the main results in Section \ref{sec:slow_manifold_theorems}. Since it allows for a simpler proof of the slow manifold Theorems \ref{thm:slow_manifolds} and \ref{thm:graph_slow_manifolds}, we first prove persistence of the stable and unstable manifolds $W^{s/u}_{loc}(S)$ defined in \eqref{eq:layer_foliation} for $0 < \eps \ll 1$.
	
	We shall consider the map \eqref{eq:map_bxa} throughout, for any $\hat \delta$ fixed in the interval
	\begin{equation}
		\label{eq:delta_bounds6}
		0 < \hat \delta < \min \left\{ \nu_R^{1/r} - 1 ,  \frac{\hat \Delta - (\Delta + c_0)}{C} ,  \frac{\Delta^2}{\hat \Delta (1 + M \hat \Delta)} \right\} .
	\end{equation}
	This ensures that the results of Section \ref{sub:NS_conds} apply for $\eps_0, \delta_0 > 0$ sufficiently small.

	\subsubsection{Proof of Theorem \ref{thm:stable_manifolds}}
	\label{ssub:proof_of_stable_manifold_persistence}

	We first prove persistence of the local stable manifold $W^s_{loc}(S)$. It follows from Lemmas \ref{lem:HM_ab}, \ref{lem:constants} and \ref{lem:CM_ar} that Theorem \cite[Theorem 1.3]{Nipp2013} applies to the map \eqref{eq:map_bxa} in the repelling formulation $\widehat H_R$, see again \eqref{eq:repelling_formulation}. This yields the existence of a positively invariant manifold
	\[
	\mathcal M^s := \left\{(\varphi_{R}(y_R,\eps),y_R) : y_R \in Y_R \right\} \subseteq X_R \times Y_R , 
	\]
	where the function $\varphi_{R}(y_R,\eps)$ is (at least) uniformly Lipschitz continuous and satisfies the invariance equation
	\begin{equation}
		\label{eq:global_invariance_Ms}
		\bar x_R = \widehat H_R^{x_R}\left( \varphi_{R}(y_R,\eps), y_R, \eps \right) = \varphi_{R} \left( \widehat H_R^{y_R}\left(\varphi_{R}(y_R,\eps), y_R, \eps\right), \eps \right) ,
	\end{equation}
	for all $y_R = (x,a)^\textnormal{T}$ such that $\widehat H_R^{y_R}(\varphi_{R}(y_R,\eps), y_R, \eps ) \in Y_R$. Combining this with Lemma \ref{lem:CM_ar}, it follows after an application of \cite[Theorem 3.1]{Nipp2013} that the function $\varphi_{R}$, and hence the manifold $\mathcal M^s$, is $C^r-$smooth in both $y_R$ and $\eps$.
	
	An analogous application of the preceding arguments in the case that $\eps = 0$ leads to
	\[
	\mathcal M^s|_{\eps = 0} = W^s_{loc}(\hat S) ,
	\]
	where $W^s_{loc}(\hat S)$ denotes the smooth extension of the stable manifold $W^s_{loc}(S)$ defined in \eqref{eq:layer_foliation} onto the enlarged domain $X_A \times Y_A$ with $x \in \hat{\mathcal K}$. From this observation, one may conclude that $\mathcal M^s$ and $W^s_{loc}(\hat S)$ are $O(\eps)-$close for $\eps \in (0,\eps_0)$ by a direct application of \cite[Theorem 2.4]{Nipp2013}. Finally, a locally (as opposed to globally) invariant manifold is obtained by restricting $\mathcal M^s$ to the original domain with $x \in \mathcal K \subset \hat{\mathcal K}$. This yields a new manifold given in $(b,x,a)-$coordinates as a graph
	\begin{equation}
		\label{eq:perturbed_stable_manifold}
		\mathcal M^s|_{B \times \mathcal K \times A} =: W^s_{loc}(S_\eps) = \left\{ (\varphi_{R}(x,a,\eps), x, a) : (x,a) \in \mathcal K \times A \right\} ,
	\end{equation}
	where in particular we have
	$\varphi_{R}(x,a,\eps) = O(\eps)$, 
	since $W^s_{loc}(S)$ and $W^s_{loc}(S_\eps)$ are $O(\eps)-$close and $W^s_{loc}(S) = \{(0,x,a):(x,a) \in \mathcal K \times A \}$ for the map \eqref{eq:map_bxa}$|_{\eps=0}$. The local invariance property of $W^s_{loc}(S_\eps)$ in Theorem \ref{thm:stable_manifolds} follows from the global (positive) invariance of $\mathcal M^s$, i.e.~the property that $\widehat H_R(\mathcal M^s,\eps) \subseteq \mathcal M^s$, after restricting to $B \times \mathcal K \times A$. 
	The preceding arguments show that the perturbed stable manifold $W^s_{loc}(S_\eps)$ has a graph representation which locally satisfies the relevant properties (i), (ii) and (iii) of Theorem \ref{thm:stable_manifolds}. 
	
	\
	
	Persistence of the local unstable manifold $W^u_{loc}(S)$ follows using similar arguments to those presented for the persistence of $W^s_{loc}(S)$, so we shall restrict ourselves to an overview of the proof. Applying \cite[Theorem 1.5]{Nipp2013} to \eqref{eq:map_bxa} in the attracting formulation $\widehat H_A$ in \eqref{eq:attracting_formulation} yields the existence of a negatively invariant manifold
	\[
	\mathcal M^u := \left\{ \left(x_A, \varphi_{A}(x_A,\eps) \right) : x_A \in X_A \right\} \subseteq X_A \times Y_A ,
	\]
	where the function $\varphi_{A}(x_A,\eps)$ is $C^r-$smooth by \cite[Theorem 3.6]{Nipp2013} and satisfies the invariance equation
	\[
	\bar y_A = \widehat H_A^{y_A} \left( x_A, \varphi_A(x_A,\eps), \eps \right) = \varphi_A \left( \widehat H_A^{x_A} \left( x_A, \varphi_A(x_A,\eps), \eps \right) , \eps \right) .
	\]
	Consideration of the $\eps = 0$ case shows that $\mathcal M^u|_{\eps = 0} = W^u_{loc}(\hat S)$, which is $O(\eps)-$close to $\mathcal M^u$ when $\eps \in (0,\eps_0)$ due to \cite[Theorem 2.1]{Nipp2013}. Restricting to $\mathcal K \subset \hat{\mathcal K}$ yields the (locally invariant) perturbed unstable manifold
	\[
	W^u_{loc}(S_\eps) := \left\{ (b, x, \varphi_{A}(b,x,\eps)) : (b,x) \in B \times \mathcal K \right\} ,
	\]
	described in Theorem \ref{thm:stable_manifolds}, where $\varphi_{A}(b,x,\eps) = O(\eps)$. Note that we also gain invertibility of the restricted map $\widehat H_A|_{\mathcal M^u}$ via \cite[Theorem 1.5]{Nipp2013}, implying invertibility of $H|_{W^u_{loc}(S_\eps)}$ as required in Theorem \ref{thm:stable_manifolds}. This allows for a straightforward derivation of the negative invariance condition in assertion (iii) using the fact that $\mathcal M^u \subseteq \widehat H_A(\mathcal M^u,\eps)$.
	
	\
	
	Aside from the final assertion about the intersection $S_\eps = W^s_{loc}(S_\eps) \cap W^u_{loc}(S_\eps)$, which will be considered in the next section, the preceding arguments combine to prove a local graph formulation of Theorem \ref{thm:stable_manifolds}. Standard arguments using compactness and a partition of unity complete the proof.
	\qed

	\subsubsection{Proof of Theorems \ref{thm:slow_manifolds} and \ref{thm:graph_slow_manifolds}}
	
	We now prove Theorem \ref{thm:graph_slow_manifolds}, which is sufficient to prove Theorem \ref{thm:slow_manifolds}.
	
	Much of the content of Theorem \ref{thm:graph_slow_manifolds} can be derived using Theorem \ref{thm:stable_manifolds} and its proof in Section \ref{ssub:proof_of_stable_manifold_persistence} by defining the `slow manifold' $S_\eps$ as the intersection of perturbed stable and unstable manifolds $W^s_{loc}(S_\eps)$ and $W^u_{loc}(S_\eps)$. For the map \eqref{eq:map_bxa}, the proof of Theorem \ref{thm:stable_manifolds} in Section \ref{ssub:proof_of_stable_manifold_persistence} together with \cite[Theorem 1.7]{Nipp2013} implies that this intersection is well-defined and given by
	\begin{equation}
		\label{eq:S_eps_intersection}
		S_\eps := W^s_{loc}(S_\eps) \cap W^u_{loc}(S_\eps) = \left\{ (r_b(x),x,r_a(x)) : x \in \mathcal K \right\} ,
	\end{equation}
	where $r_b(x)$ and $r_a(x)$ are $C^r$-smooth functions which satisfy
	\begin{equation}
		\label{eq:ra_rb}
		r_b(x) = O(\eps) , \qquad r_a(x) = O(\eps) .
	\end{equation}
	Invertibility of the restricted map $H|_{S_\eps}$ follows from the fact that $H|_{\mathcal M^u}$ is invertible, since $S_\eps \subset \hat S_\eps \subseteq \mathcal M^u$. The fact that $S_\eps$ is $O(\eps)-$close and diffeomorphic to the critical manifold $S = \{(x,0) : x \in \mathcal K\}$ follows by \eqref{eq:S_eps_intersection} and \eqref{eq:ra_rb}. The local invariance of $S_\eps$ described in Theorem \ref{thm:graph_slow_manifolds} follows from the global invariance $H(\hat S_\eps,\eps) = \hat S_\eps$ implied by \cite[Theorem 1.7]{Nipp2013} after restriction to $\mathcal K \subset \widehat{\mathcal K}$, and 
	the particular form of the invariance equations \eqref{eq:graph_sm_forward_inv} and \eqref{eq:graph_sm_backward_inv} in Theorem \ref{thm:graph_slow_manifolds} follow after expressing the local invariance requirement
	\[
	z_\eps \in S_\eps , \ H(z_\eps,\eps) \in \mathcal U \quad 
	\implies \quad
	\bar z_\eps = H(z_\eps,\eps) \in S_\eps ,
	\]
	in the local $(x,y)$ coordinates of \eqref{eq:main_map_graph}. Finally, the form of the expansion \eqref{eq:slow_manifold_approx} follows from its derivation in the proof of Proposition \ref{prop:reduced_map}, see again equation \eqref{eq:slow_manifold} and Remark \ref{rem:slow_manifold_expansion}.
	
	Thus we have proved Theorem \ref{thm:graph_slow_manifolds} and, correspondingly, the local expression of Theorem \ref{thm:slow_manifolds}. Extending the local result via compactness and a partition of unity completes the proof of Theorem \ref{thm:slow_manifolds}.
	\qed

	\subsubsection{Proof of Theorem \ref{thm:foliations}}
	\label{ssec:foliation_proof}
	
	Applying \cite[Theorem 4.1]{Nipp2013} to the map \eqref{eq:map_bxa} in the repelling formulation $\widehat H_R$, 
	we find that the perturbed stable manifold $W^s_{loc}(\hat S_\eps)$ of the map \eqref{eq:map_bxa} admits a foliation by stable fibers
	$w^s_{loc}(z_\eps)$ 
	with base points $z_\eps \in \hat S_\eps$, which 
	are $C^r-$smooth in $\eps$ and satisfy the (global) positive invariance property
	\[
	\widehat H_R(w^s_{loc}(z_\eps)) \subseteq w^s_{loc} \left(\widehat H_R(z_\eps) \right) ,
	\]
	for all $z_\eps \in \hat S_\eps$. 
	This motivates the definition
	\begin{equation}
		\label{eq:stable_foliation_enlarged}
		W^s_{loc}(\hat S_\eps) := \bigcup_{z_\eps \in \hat S_\eps} w_{loc}^s(z_\eps) .
	\end{equation}
	Theorem 4.1 in \cite{Nipp2013} also ensures that the fibers $w_{loc}^s(z_\eps)$ are continuous and identical to the set of points $z_R \in X_R \times Y_R$ for which iterates $\widehat H_R^j(z_R,\eps)$ are exponentially attracted to iterates $\widehat H_R^j(z_\eps,\eps)$ along $\hat S_\eps$. 
	More precisely, $w_{loc}^s(z_\eps)$ is the set of points $z_R \in X_R \times Y_R$ such that for each $j \in \mathbb N$ 
	we have
	\[
	\big| \widehat H_R^j(z_R,\eps) - \widehat H_R^j(z_\eps,\eps) \big| \leq c \chi_A^j ,
	\]
	for some $c > 0$ and any fixed $\chi_A \in (\nu_A,1)$. Since the limit $\eps \to 0$ is defined (fibers are $C^r-$smooth in $\eps$) and the stable fibers for the layer map 
	also admit such a characterisation, it follows that $W_{loc}^s(\hat S_\eps)$ converges fiberwise to $W_{loc}^s(\hat S)$ as $\eps \to 0$. The fact that the fibers are $C^{r-1}-$smooth overall, i.e.~with respect to coordinate and base point variation, follows from \cite[Theorem 5.1]{Nipp2013}. The nondegeneracy condition in assertion (v) is a consequence of the requirement \eqref{eq:nondegeneracy_fiber_smoothness} imposed in order to prove Lemma \ref{lem:conds_smoothness_foliation}, see also Remark \ref{rem:foliaiton_smoothness_local}. Smoothness and convergence for $\eps \to 0$ together imply that $W_{loc}^s(\hat S_\eps)$ and $W_{loc}^s(\hat S)$ are $O(\eps)-$close. Restricting to the original domain $\mathcal K \subset \hat{\mathcal K}$ yields the desired results. In particular, restricting to $\mathcal K$ amounts to a restriction of \eqref{eq:stable_foliation_enlarged} to $S_\eps \subset \hat S_\eps$, yielding the stable foliation $W^s_{loc}(S_\eps) \subset W^s_{loc}(\hat S_\eps)$ from equation \eqref{eq:foliations_perturbed}. The local invariance and contraction properties in assertions (i) and (iii) follow from the corresponding global invariance and contraction properties identified above on the enlarged domain.

	\
	
	Similar arguments based on \cite[Theorems 4.2 and 5.2]{Nipp2013} prove the corresponding results for the unstable foliation; here we omit the details for brevity, noting only that the additional requirement that $H$ is invertible along unstable fibers $w^u_{loc}(z_\eps)$ follows from the invertibility of $H$ on $W^u_{loc}(S_\eps)$ and equation \eqref{eq:foliations_perturbed}, see again Theorem \ref{thm:stable_manifolds} (iii). The preceding arguments prove Theorem \ref{thm:foliations} locally, and can be extended over compact domains in order to complete the proof using a partition of unity.
	\qed

	\section{Applications}
	\label{sec:examples}
	
	Having proved the main results, we consider a number of applications of the DGSPT formalism developed in Sections \ref{sec:a_coordinate-independent_framework_for_fast-slow_maps}-\ref{sec:slow_manifold_theorems}. We begin in Section \ref{sub:Chialvo_model} with a geometric analysis of a two-dimensional map-based model for neuronal dynamics originally presented in \cite{Chialvo1995}. This example will help to demonstrate the application of the theory, and is intended to provide a kind of `benchmark application' similar to the Van der Pol oscillator \cite{vdP1920,vdP1926} for continuous-time fast-slow ODE systems. Sections \ref{sub:Euler_discretization}-\ref{sub:Poincare_maps} are more theoretical, and demonstrate the utility of the theory for analysing discretizations and Poincar\'e maps induced by fast-slow ODE systems respectively.

	\subsection{A map-based neural model}
	\label{sub:Chialvo_model}
	
	
	%
	%
	%
	%

	We consider a map-based model for neuronal bursting originally due to Chialvo \cite{Chialvo1995} and considered further in e.g.~\cite{Jing2006,Muni2022,Trujillo2021,Wang2018}. The model is given by
	\begin{equation}
		\label{eq:Chialvo_model}
		\begin{split}
			w \mapsto & \ \bar w = a w - b v + c, \\
			v \mapsto & \ \bar v = v^2 \exp(w-v) + k, 
		\end{split}
	\end{equation}
	where $v$ denotes membrane potential voltage, $w$ is a recovery variable, $k \geq 0$ models a time-dependent perturbation of the voltage, and $a \in (0,1)$, $b \in (0,1)$, $c>0$ are parameters relating to the recovery process; see \cite{Chialvo1995} for a physical interpretation.
	
	
	Recently in \cite{Trujillo2021}, the authors considered a reduced $1-$dimensional model obtained by considering the recovery variable $w$ as a parameter. Although not mathematically justified, this assumption is motivated by the observation in \cite{Chialvo1995} that \eqref{eq:Chialvo_model} exhibits fast-slow dynamics in particular parameter regimes. This assumption can be justified in suitable parameter regimes using the theory developed in Sections \ref{sec:a_coordinate-independent_framework_for_fast-slow_maps}-\ref{sec:slow_manifold_theorems}, i.e.~using DGSPT. Specifically, the map in which $w$ is treated as a parameter can be viewed as a layer map obtained in the singular limit over regions of $(a,b,c)-$parameter space such that
	\begin{equation}
		\label{eq:fast_slow_cond}
		a w - b v + c = w + \eps g(w,v,\eps) ,
	\end{equation}
	where $0 < \eps \ll 1$ and $g : \mathbb R \times \mathbb R \times [0,\eps_0) \to \mathbb R$ is $C^r-$smooth. In order for the condition \eqref{eq:fast_slow_cond} to be satisfied for $\eps \to 0$ uniformly in $(w,v)-$phase space, we impose the following model assumption.
	
	\begin{assumption}
		\label{ass:Chialvo_fast_slow}
		\textup{(Restriction to fast-slow parameter regime)} The parameters $a \in (0,1)$, $b \in (0,1)$ and $c > 0$ depend smoothly on $\eps$ and satisfy
		\[
		\lim_{\eps \to 0} a(\eps) = 1^- , \qquad 
		\lim_{\eps \to 0} b(\eps) = 0^+ , \qquad
		\lim_{\eps \to 0} c(\eps) = 0^+ ,
		\]
		such that in particular we have
		\[
		a(\eps) = 1 - \tilde a \eps + O(\eps^2) , \qquad
		b(\eps) = \tilde b \eps + O(\eps^2) , \qquad 
		c(\eps) = \tilde c \eps + O(\eps^2) ,
		\]
		for some $\tilde a > 0$, $\tilde b > 0$ and $\tilde c > 0$.
	\end{assumption} 
	
	Thus we consider
	\begin{equation}
		\label{eq:Chialvo}
		\begin{split}
			w \mapsto & \hspace{2pt} \bar{w} = w + \eps g(w,v,\eps) , \\
			v \mapsto & \hspace{2.5pt} \bar{v} \hspace{1pt} = v + f(w,v) , \\
		\end{split}
	\end{equation}
	where
	\[
	f(w,v) = - v + v^2 \exp(w-v) + k, \qquad 
	g(w,v,\eps) = \tilde c - \tilde b v - \tilde a w + O(\eps) .
	\]
	The map \eqref{eq:Chialvo} is in the standard form for fast-slow maps, and can be written in the general form \eqref{eq:gen_maps} after setting $N(v,v) = (0,1)^\textnormal{T}$ and $G(w,v,\eps) = (g(w,v,\eps),0)^\textnormal{T}$ so that \eqref{eq:Chialvo} becomes
	\begin{equation}
		\label{eq:Chialvo_nonstandard}
		\begin{pmatrix}
			w \\
			v
		\end{pmatrix}
		\mapsto 
		\begin{pmatrix}
			\bar w \\
			\bar v
		\end{pmatrix}
		=
		\begin{pmatrix}
			w \\
			v
		\end{pmatrix}
		+
		\begin{pmatrix}
			0 \\
			1
		\end{pmatrix}
		f(w,v) + \eps
		\begin{pmatrix}
			g(w,v,\eps) \\
			0
		\end{pmatrix}
		.
	\end{equation}
	In what follows, we refer to \eqref{eq:Chialvo} or equivalently \eqref{eq:Chialvo_nonstandard} as the \textit{fast-slow Chialvo map}.
	
	\begin{remark}
		The authors in \cite{Trujillo2021} point out that the map \eqref{eq:Chialvo_model} is fast-slow for $ 0 < 1 - a = b = c := \eps \ll 1$. This corresponds to the particular case with $\tilde a = \tilde b = \tilde c = 1$ in \eqref{eq:Chialvo}. 
	\end{remark}
	
	\begin{remark}
		Further examples of fast-slow maps for which the theory developed herein applies include the so-called \textit{Izhikevich model} \cite{Izhikevich2004} obtained by Euler discretization of a continuous-time model in \cite{Izhikevich2003}, different variants of the \textit{Rulkov model} \cite{Rulkov2001,Rulkov2002,Shilnikov2004}, and the \textit{Courbage-Nekorkin-Vdovin model} appearing in \cite{Nekorkin2007} and later modified in \cite{Courbage2007}. We refer to the review article \cite{Ibarz2011} and the many references therein.
	\end{remark}

	\subsubsection{Layer map and slow manifolds}
	
	Due to the explicit separation of `time-scales' in \eqref{eq:Chialvo}, the layer map
	\begin{equation}
		\label{eq:Chialvo_layer}
		\begin{pmatrix}
			w \\
			v
		\end{pmatrix}
		\mapsto 
		\begin{pmatrix}
			\bar w \\
			\bar v
		\end{pmatrix}
		=
		\begin{pmatrix}
			w \\
			v
		\end{pmatrix}
		+
		\begin{pmatrix}
			0 \\
			1
		\end{pmatrix}
		f(w,v)
	\end{equation}
	reduces to the $1-$dimensional map $v \mapsto v + f(w,v)$ with parameters $w \in \mathbb R$ and $k \geq 0$ considered in detail in \cite{Trujillo2021}. The critical manifold $S=\{(w,v) \in \mathbb R^2 : f(w,v) = 0\}$ can be written as a graph
	\[
	S = \left\{ (\varphi_0(v),v) : v > k \right\} , \qquad 
	\varphi_0(v) := v + \ln \left( \frac{v-k}{v^2} \right) ,
	\]
	and by Proposition \ref{prop:EVs}, the (unique) non-trivial multiplier along $S$ is given by
	\begin{equation}
		\label{eq:Chialvo_mu}
		\mu(v) = 1 + Df N |_S = 1 + D_vf|_S = \frac{(v-k)(2-v)}{v} , \qquad v > k.
	\end{equation}
	Non-normally hyperbolic points on $S$ are points $(\varphi_0(v),v)$ such that $|\mu(v)| = 1$. Since $\mu(v) \in \mathbb R$, this can only occur for fold-type singularities with $\mu(v) = 1$, or flip-type singularities with $\mu(v) = -1$.
	
	\
	
	We obtain the following result, which can also (for the most part) be found in \cite{Trujillo2021}; we simply restate it here in a manner more suited to our formalism.
	
\begin{figure}[t!]
		\centering
		\includegraphics[scale=0.5]{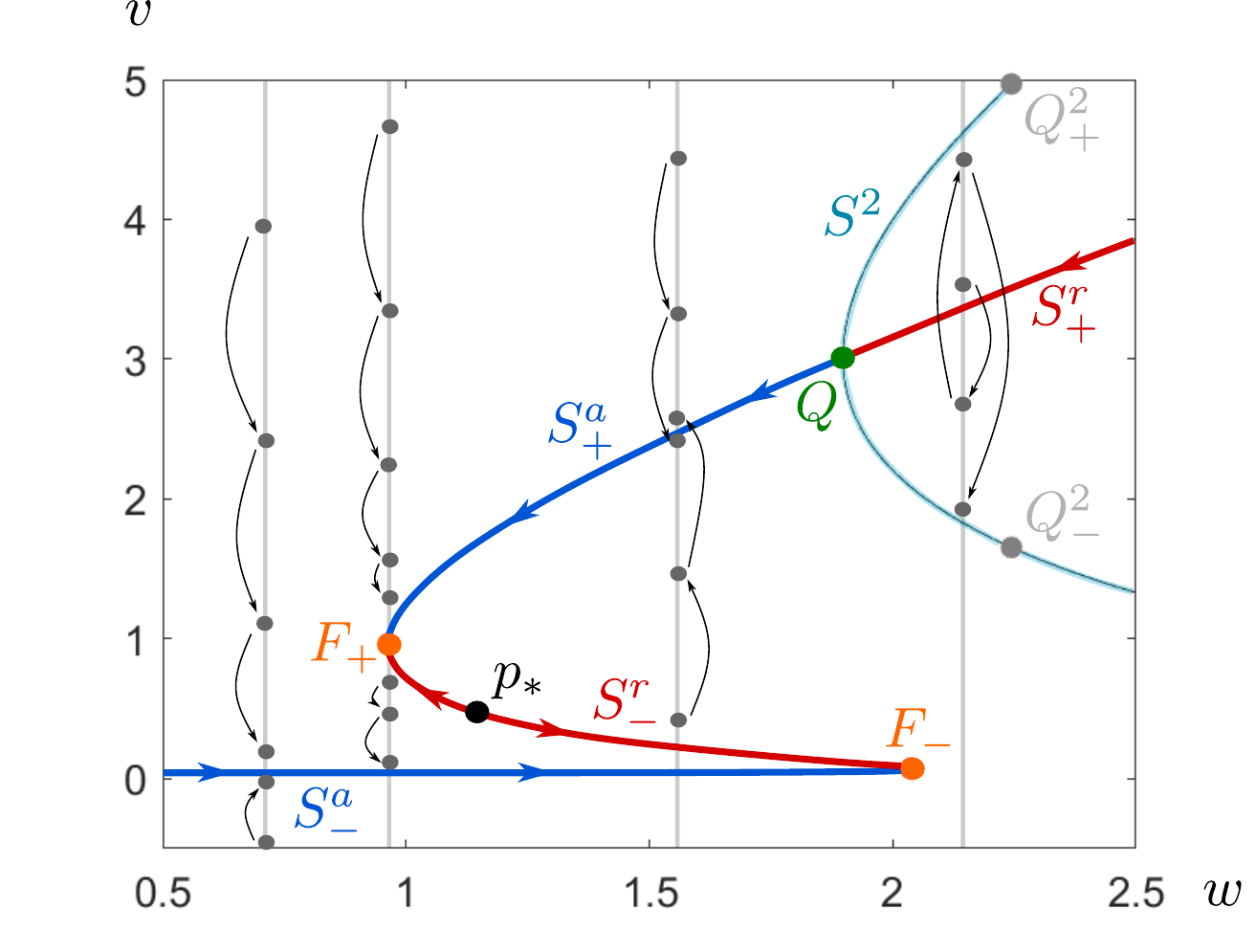}
		\caption{Geometry and dynamics of the fast-slow Chialvo map \eqref{eq:Chialvo} in the singular limit $\eps = 0$, with parameter values $\tilde a = 1$, $\tilde b = 5$, $\tilde c = 3.5$ and $k = 0.035$. The critical manifold $S$ has a cubic-like profile, with normally hyperbolic and attracting (repelling) critical manifolds $S^a_{\pm}$ ($S^r_{\pm}$) shown in blue (red) separated by non-hyperbolic fold points $F_\pm$ (orange) and a supercritical flip point $Q$ (green); see Lemma \ref{lem:Chialvo_layer}. Iterates of the layer map \eqref{eq:Chialvo_layer} along vertical fast fibers (shaded grey) are sketched along with the direction of iteration, in order to illustrate the fast dynamics. The reduced map, i.e.~\eqref{eq:Chialvo_slow} truncated at $O(\eps^2)$, has a unique unstable fixed point $p_\ast \in S_-^r$, see Proposition \ref{prop:Chialvo_fixed_points}. The direction of `slow iteration' under the reduced map is indicated with single arrows. Also shown in cyan is the (numerically computed) curve of period-2 points $S^2$ bifurcating from the flip point $Q$, as well as numerically identified flip bifurcations $Q_{\pm}^2$ in the second iterate map $v \mapsto \bar v^2$. A period-doubling route to chaos occurs with increasing $w$; see \cite{Trujillo2021} and Remark \ref{rem:Chialvo_cascade}.}
		\label{fig:Chialvo_singular}
	\end{figure}

	\begin{lemma}
		\label{lem:Chialvo_layer}
		Fix $k \in (0,3-2\sqrt{2})$. Then the critical manifold decomposes like
		\[
		S = S_-^a \cup F_- \cup S_-^r \cup F_+ \cup S_+^a \cup Q \cup S_+^r ,
		\]
		where $S_-^a$, $S_+^a$ are normally hyperbolic and attracting, $S_-^r$, $S_+^r$ are normally hyperbolic and repelling, $F_-$, $F_+$ and $Q$ are respectively fold, fold and supercritical flip points for the layer map \eqref{eq:Chialvo_layer}.
		
		The fold points are located at $F_\pm : (\varphi_0(v_{\pm}), v_\pm)$, where $v_\pm = ( 1 + k \mp \sqrt{k^2 - 6k + 1} ) / 2$, 
		%
		the flip point is located at $Q : (\varphi_0(v_{flip}) , v_{flip})$, where $v_{flip} = ( 3 + k + \sqrt{k^2 - 2k + 9} ) / 2$, and
		\[
		\begin{aligned}
			S_-^a &= S \cap \{ v \in (k,v_-) \} , \qquad 
			&&S_-^r = S \cap \{ v \in (v_-,v_+) \} , \\
			S_+^a &= S \cap \{ v \in (v_+,v_{flip}) \} , \qquad 
			&&S_+^r = S \cap \{ v > v_{flip} \} ,
		\end{aligned}
		\]
		see Figure \ref{fig:Chialvo_singular}.
		
		If $k=0$, there is only one fold point $F_+ : (1,1)$, and the critical manifold decomposes like $S = S_-^r \cup F_+ \cup S_+^a \cup Q \cup S_+^r$.
		%
		%
	\end{lemma}
	
	\begin{proof}
		The existence of fold and flip bifurcations can be verified directly by solving $\mu(v) = 1$ and $\mu(v) = -1$ respectively for $v$, where $\mu(v)$ is given by \eqref{eq:Chialvo_mu}, and by checking the relevant non-degeneracy conditions in e.g.~\cite[Theorem 4.1]{Kuznetsov2013} and \cite[Theorem 4.3]{Kuznetsov2013} respectively. Details and explicit calculations are given in \cite[Theorems 2.2-2.3]{Trujillo2021}.
	\end{proof}
	
	In the following we restrict to $k \in (0,3-2\sqrt{2})$, such that the critical manifold is `$S-$shaped' as in Figure \ref{fig:Chialvo_singular}. For further details on the dynamics in the layer problem for $k = 0$ and $k \geq 3-2\sqrt{2}$, we refer again to \cite{Trujillo2021}.
	
	\
	
	Persistence of (compact submanifolds of) normally hyperbolic branches of $S$ as slow manifolds for $0 < \eps \ll 1$ follows by the results in Section \ref{sec:slow_manifold_theorems}. 
	
	\begin{proposition}
		\label{prop:Chialvo_slow_manifolds}
		Fix $k \in (0,3 - 2\sqrt{2})$. For $0 < \eps \ll 1$ sufficiently small, compact submanifolds of $S_-^a$, $S_+^a$, $S_-^r$ and $S_+^r$ persist as nearby slow manifolds $S_{-,\eps}^a$, $S_{+,\eps}^a$, $S_{-,\eps}^r$ and $S_{+,\eps}^r$ respectively. Each slow manifold is given by the graph of
		\begin{equation}
			\label{eq:Chialvo_slow_manifolds}
			\varphi_\eps(v) = v + \ln \left( \frac{v-k}{v^2} \right) - \eps \left( \frac{ \tilde c - (\tilde b + \tilde a) v - \tilde a \ln \left(\frac{v - k}{v^2} \right) }{\mu(v)-1} \right) + O(\eps^2) ,
		\end{equation}
		over a suitable compact $v-$interval, e.g.~$S_{+,\eps}^a = \{(\varphi_\eps(v),v) : v \in [v_+ + \delta, v_{flip} - \delta ] \}$ for an arbitrarily small but fixed $\delta > 0$.
	\end{proposition}
	
	\begin{proof}
		This follows directly from Theorem \ref{thm:graph_slow_manifolds} with notation $(x,y) = (v,w)$, and $N$, $f$ and $G$ defined as in \eqref{eq:Chialvo_nonstandard}. Substituting the relevant expressions into \eqref{eq:slow_manifold_approx} gives
		\[
		\begin{split}
			\varphi_\eps(v) &= \varphi_0(v) - \eps \left[ (D_wf)^{-1} (DfN)^{-1} Df G \right] (\varphi_0(v), v) + O(\eps^2) \\
			&= \varphi_0(v) - \eps \left( \frac{g(\varphi_0(v),v,0)}{\mu(v)-1} \right) + O(\eps^2) .
		\end{split}
		\]
		The expression \eqref{eq:Chialvo_slow_manifolds} follows after substituting the expressions for $g$ and $\varphi_0(v)$.
	\end{proof}

	The slow manifolds $S_{-,\eps}^a$, $S_{+,\eps}^a$, $S_{-,\eps}^r$ and $S_{+,\eps}^r$ identified in Proposition \ref{prop:Chialvo_slow_manifolds} have all the properties described in Theorems \ref{thm:slow_manifolds}, \ref{thm:graph_slow_manifolds}, \ref{thm:stable_manifolds} and \ref{thm:foliations}. 
	
	\begin{remark}
		\label{rem:Chialvo_cascade}
		Classical bifurcation theory implies the existence of a locally quadratic, normally hyperbolic and attracting critical manifold $S^2$ for the second iterate (layer) map $v \mapsto \bar v^2$. The manifold $S^2$ 
		intersects $S$ transversally at the supercritical flip point $Q$, see Figure \ref{fig:Chialvo_singular}. Dynamically, $S^2$ forms a manifold of stable period-2 cycles for the layer map \eqref{eq:Chialvo_layer}. In fact there is a period-doubling cascade for increasing $w$ which leads to chaotic dynamics in the layer map. We do not consider chaotic dynamics in detail in this work, but refer to \cite{Trujillo2021} for details on chaotic dynamics in the layer map \eqref{eq:Chialvo_layer} and \cite{Baesens1991} for details on the slow passage through a period-doubling cascade.
	\end{remark}

	\subsubsection{Reduced map and dynamics on slow manifolds}
	
	
	To leading order in $\eps$, the dynamics on slow manifolds $S_{-,\eps}^a$, $S_{+,\eps}^a$, $S_{-,\eps}^r$ and $S_{+,\eps}^r$ is determined by the reduced map \eqref{eq:red_map}. This is obtained directly using Proposition \ref{prop:reduced_map}.
	
	\begin{lemma}
		\label{lem:Chialvo_slow}
		Slow iteration along $S_{-,\eps}^a$, $S_{+,\eps}^a$, $S_{-,\eps}^r$ and $S_{+,\eps}^r$ is governed by the map
		\begin{equation}
			\label{eq:Chialvo_slow}
			\begin{pmatrix}
				w \\
				v
			\end{pmatrix}
			\mapsto 
			\begin{pmatrix}
				\bar w \\
				\bar v
			\end{pmatrix}
			=
			\begin{pmatrix}
				\varphi_0(v) \\
				v
			\end{pmatrix}
			+ \eps
			\begin{pmatrix}
				1  \\
				- \tfrac{v-k}{\mu(v) - 1}
			\end{pmatrix}
			\left( \tilde c - (\tilde b + \tilde a) v - \tilde a \ln \left(\frac{v - k}{v^2} \right) \right) + O(\eps^2) ,
		\end{equation}
		for all $(w,v) = (\varphi_0(v),v) \in S \setminus (F_- \cup F_+)$ restricted to the corresponding (compact) $v-$interval. By invariance, this can be represented by the $1-$dimensional map
		\begin{equation}
			\label{eq:Chialvo_slow_1D}
			v \mapsto \bar v = v - \eps \left( \frac{v-k}{\mu(v) - 1} \right) \left( \tilde c - (\tilde b + \tilde a) v - \tilde a \ln \left(\frac{v - k}{v^2} \right) \right) + O(\eps^2) .
		\end{equation}
	\end{lemma}
	
	\begin{proof}
		Let $S_n := S \setminus (F_- \cup F_+)$. The form of the map \eqref{eq:Chialvo_slow} is obtained via equation \eqref{eq:red_map_expanded} in Proposition \ref{prop:reduced_map}, noting that the projection operator $\Pi_\mathcal{N}^{S_n}$ is given by
		\[
		\Pi_{\mathcal N}^{S_n} = I_2 - N (Df N)^{-1} Df |_{S_n} =
		\begin{pmatrix}
			0 & 1 \\
			0 & - (D_v f)^{-1} D_wf
		\end{pmatrix} \bigg|_{S_n}
		=
		\begin{pmatrix}
			0 & 1 \\
			0 & - \tfrac{v-k}{\mu(v)-1}
		\end{pmatrix} \bigg|_{S_n} ,
		\]
		where we used the fact that $\mu(v) = 1 + D_vf|_{S_n}$.
	\end{proof}
	
	Fixed points on $S_{-,\eps}^a$, $S_{+,\eps}^a$, $S_{-,\eps}^r$ and $S_{+,\eps}^r$ occur for values $v_{\ast,\eps} = v_\ast + O(\eps)$, where $v_\ast$ satisfies
	\begin{equation}
		\label{eq:Chialvo_fixed_points}
		g(\varphi_0(v_\ast),v_\ast,0) = \tilde c - (\tilde b + \tilde a) v_\ast - \tilde a \ln \left(\frac{v_\ast - k}{v_\ast^2} \right) = 0 .
	\end{equation}
	The intermediate value theorem guarantees that equation \eqref{eq:Chialvo_fixed_points} has at least one solution $v_\ast \in (k,\infty)$, since $g(\varphi_0(v),v,0) \sim - \tilde a \ln ((v-k)/v^2) > 0$ as $v \to k^+$ and $g(\varphi_0(v),v,0) \sim -(\tilde b + \tilde a) v < 0$ as $v \to \infty$. For simplicity we restrict to the case in which $v_\ast$ is the only solution and, correspondingly, the reduced map \eqref{eq:Chialvo_slow} has precisely one fixed point.
	
	\begin{assumption}
		\label{ass:Chialvo_fixed_point}
		\textup{(Unique equilibrium on $S$)} 
		The parameters $\tilde a > 0$, $\tilde b > 0$, $\tilde c > 0$ and $k \in (0,3-2\sqrt{2})$ are chosen such that the function $g(\varphi_0(v),v,0)$ is strictly decreasing, i.e.~$D_v g(\varphi_0(v),v,0) < 0$ for all $v > k$. This is true for parameter values satisfying
		\begin{equation}
			\label{eq:fixed_point_inequality}
			- 2 k \tilde a + (\tilde a + \tilde a k + \tilde b k) v - (\tilde a + \tilde b) v^2 < 0
		\end{equation}
		for all $v > k$.
	\end{assumption}
	
	Assumption \ref{ass:Chialvo_fixed_point} is satisfied for the parameter values in Figure \ref{fig:Chialvo_singular}, for example, on the condition that $-0.07 + 1.21 v - 6 v^2 < 0$ for all $v > k = 0.035$. This is indeed the case, since this expression is negative at e.g.~$v=1$, and the discriminant $\Delta \approx -0.22 < 0$.
	
	\begin{remark}
		For other parameter choices the graph of $g(\varphi_0(v),v,0)$ has a `cubic' profile, with a unique local minimum and a unique local maximum. In this case equation \eqref{eq:Chialvo_fixed_points} can have either $1$, $2$ or $3$ solutions depending on whether the local minimum is positive, zero or negative respectively. A complete description of the partitioning of the phase portrait in $(\tilde a, \tilde b, \tilde c, k)-$space is beyond the scope of this article.
	\end{remark}
	
	The following result classifies the stability of the perturbed fixed point $p_{\ast,\eps}$ depending on whether $p_\ast : (\varphi_0(v_\ast),v_\ast)$ lies on $S_-^a$, $S_+^a$, $S_-^r$ or $S_+^r$.
	
	\begin{proposition}
		\label{prop:Chialvo_fixed_points}
		Consider the fast-slow Chialvo map \eqref{eq:Chialvo} with parameter values in the set defined by Assumption \ref{ass:Chialvo_fixed_point}. The following assertions hold for all $0 < \eps \ll 1$ sufficiently small:
		\begin{enumerate}
			\item[(i)] For parameter values $\tilde a, \tilde b, \tilde c$ and $k$ such that $p_\ast \in S_-^a$, $S_+^a$ or $S_+^r$, there exists a unique fixed point $p_{\ast,\eps} \in S_{-,\eps}^a$, $S_{+,\eps}^a$ or $S_{+,\eps}^r$ respectively, which is $O(\eps)-$close to $p_\ast$ and asymptotically stable within $S_{-,\eps}^a$, $S_{+,\eps}^a$ or $S_{+,\eps}^r$ respectively.
			\item[(ii)] For parameter values $\tilde a, \tilde b, \tilde c$ and $k$ such that $p_\ast \in S_-^r$, there exists a unique fixed point $p_{\ast,\eps} \in S_{-,\eps}^r$ which is $O(\eps)-$close to $p_\ast$ and asymptotically unstable within $S_{-,\eps}^r$.
		\end{enumerate}
	\end{proposition}
	
	\begin{proof}
		First we show that if $p_\ast$ is a normally hyperbolic point on $S$ then it persists as an $O(\eps)-$close fixed point of the fast-slow Chialvo map \eqref{eq:Chialvo} for $0 < \eps \ll 1$. By Theorem \ref{thm:slow_manifolds}, $p_\ast$ perturbs to a point $p_{\ast,\eps}$ on one of the slow manifolds $S_{-,\eps}^a$, $S_{+,\eps}^a$, $S_{-,\eps}^r$ or $S_{+,\eps}^r$. It therefore suffices to prove the result for the map \eqref{eq:Chialvo_slow_1D}. Since $g(\varphi_0(v_\ast),v_\ast,0) = 0$ and by Assumption \ref{ass:Chialvo_fixed_point} $D_v g(\varphi_0(v_\ast), v_\ast, 0) < 0$, the implicit function theorem implies persistence of $p_\ast$ as an $O(\eps)-$close fixed point $p_{\ast,\eps}$ as required.
		
		It remains to prove the stability claims in statements (i)-(ii). The Jacobian associated to the right-hand-side of \eqref{eq:Chialvo_slow_1D} has a unique multiplier
		\[
		\mu(v_{\ast,\eps}) = 1 - \eps \left( \frac{v_\ast-k}{\mu(v_\ast) - 1} \right) D_v g(\varphi_0(v_\ast), v_\ast, 0) + O(\eps^2) ,
		\]
		where $v_{\ast,\eps} = v_\ast + O(\eps)$ is the $v-$coordinate of $p_{\ast,\eps}$. Since $v_\ast - k > 0$ and by Assumption \ref{ass:Chialvo_fixed_point} $D_v g(\varphi_0(v_\ast), v_\ast, 0) < 0$, the stability claims in statements (i)-(ii) follow from the fact that $\mu(v) - 1 < 0$ on $S_-^a \cup S_+^a \cup S_+^r$ and $\mu(v) - 1 > 0$ on $S_-^r$.
	\end{proof}
	
	We note that each possibility in Proposition \ref{prop:Chialvo_fixed_points} is realised on open parameter sets. For the particular choice of parameters in Figure \ref{fig:Chialvo_singular}, a unique asymptotically unstable fixed point $p_\ast \in S_-^r$ perturbs to a unique asymptotically unstable fixed point on $S_{-,\eps}^r$ for $0 < \eps \ll 1$.

	\subsubsection{Global dynamics for $0 < \eps \ll 1$}

	Using the preceding analysis, we can provide a partial geometric description of the global dynamics of the fast-slow Chialvo map in suitable regions of phase space, specifically, away from non-normally hyperbolic points on $S$. We consider four possibilities, distinguished on the `singular level' according to (i) which branch the fixed point $p_\ast$ is on, and (ii) the relative positioning of the lower fold point $F_-$ and the flip point $Q$. We restrict to the case of a single fixed point on $S$ in accordance with Assumption \ref{ass:Chialvo_fixed_point} throughout. Numerical computations are carried out using the software package MatContM \cite{MatContM}. 
	
	
	
	\begin{remark}
		A full classification of bursting mechanisms in the fast-slow Chialvo map \eqref{eq:Chialvo} is beyond the scope of this work. We refer to \cite{Izhikevich2004} for more on the classification of bursting mechanisms in two-dimensional fast-slow maps.
	\end{remark}

	\begin{figure}[t!]
		\centering
		\subfigure[Case I: Excitability.]{\includegraphics[width=.49\textwidth]{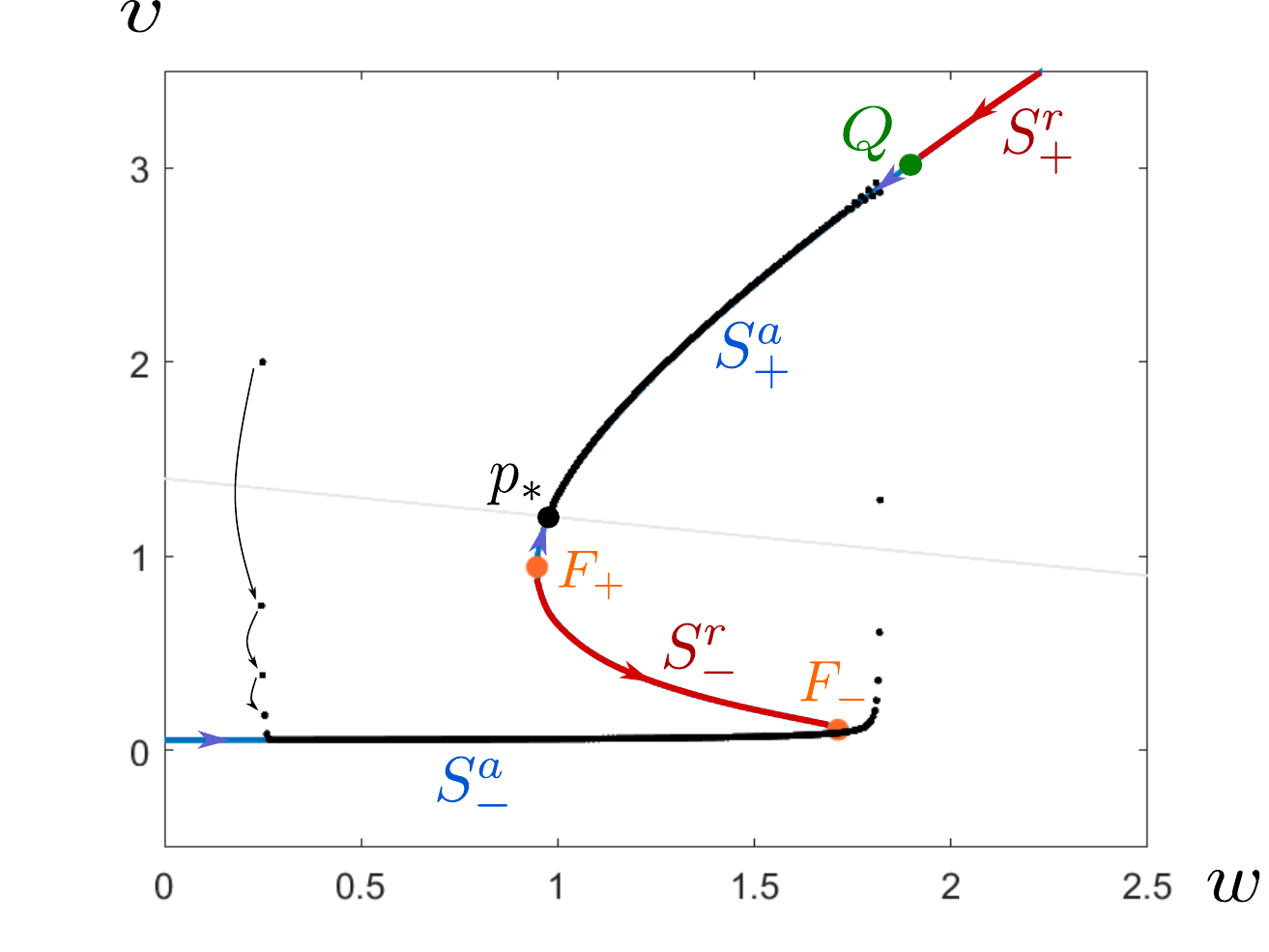}}
		\ 
		\subfigure[Case II: Relaxation.]{\includegraphics[width=.49\textwidth]{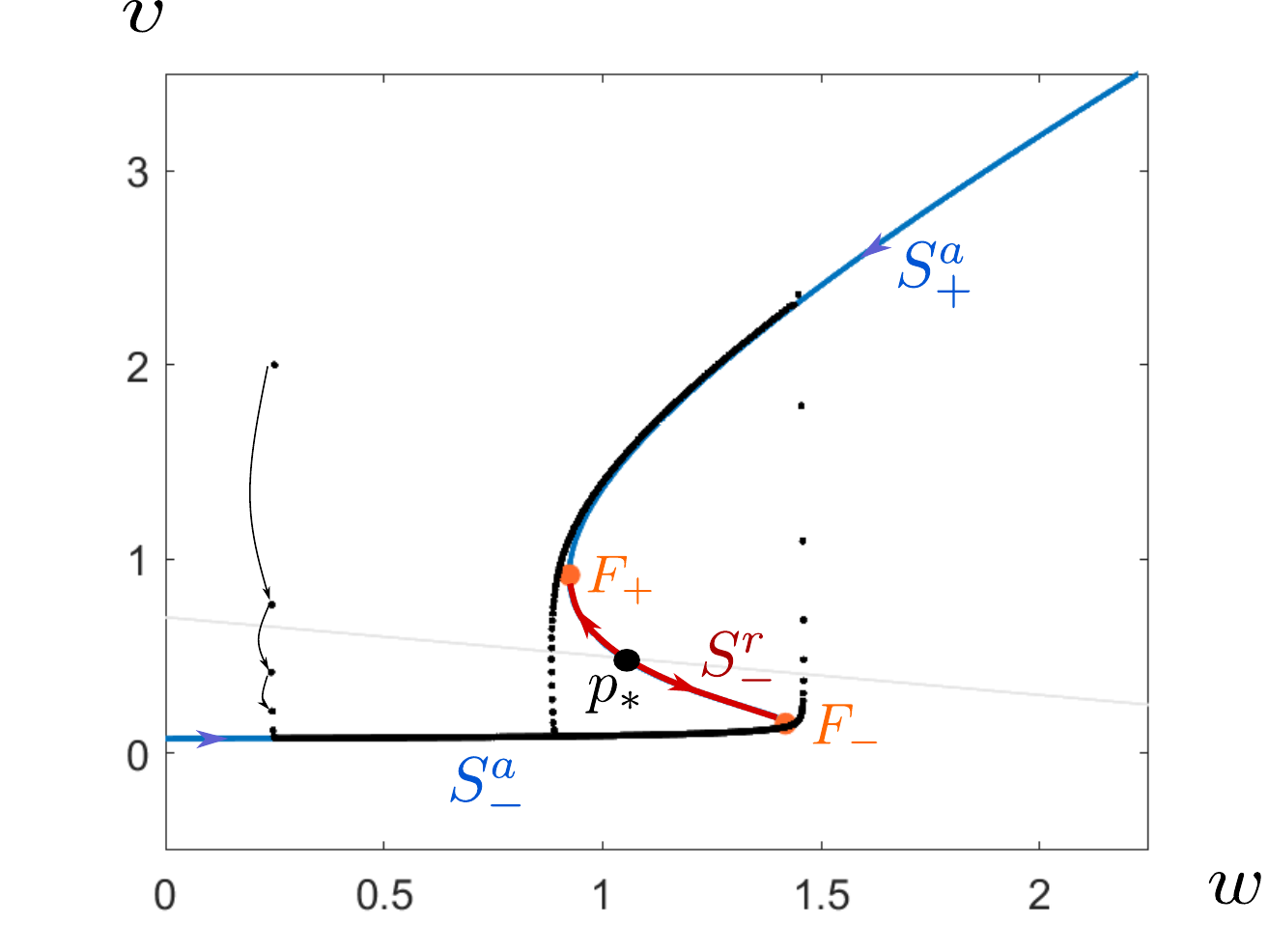}}
		\\
		\subfigure[Case III: Non-chaotic bursting.]{\includegraphics[width=.49\textwidth]{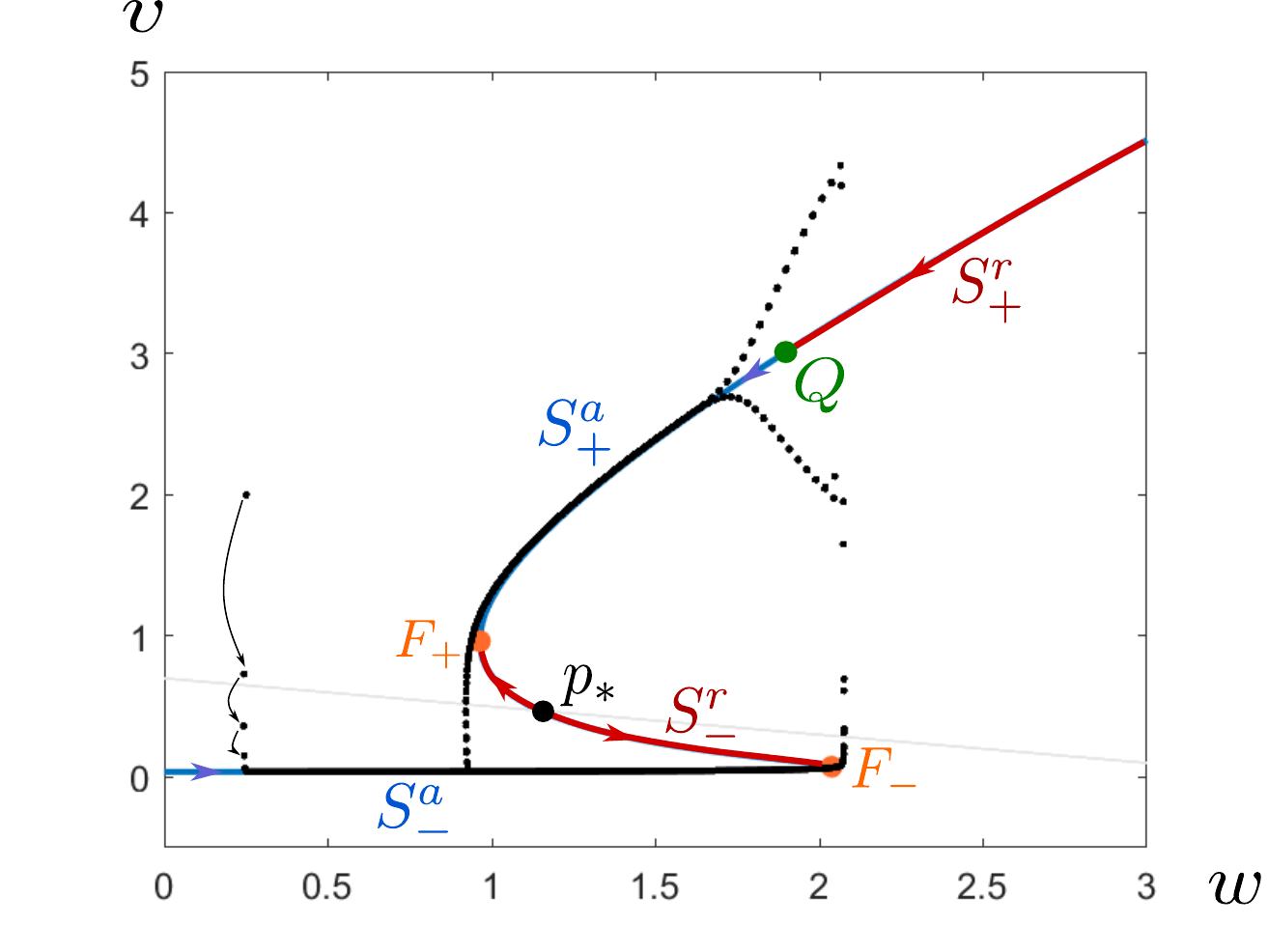}}
		\ 
		\subfigure[Case IV: Chaotic bursting.]{\includegraphics[width=.49\textwidth]{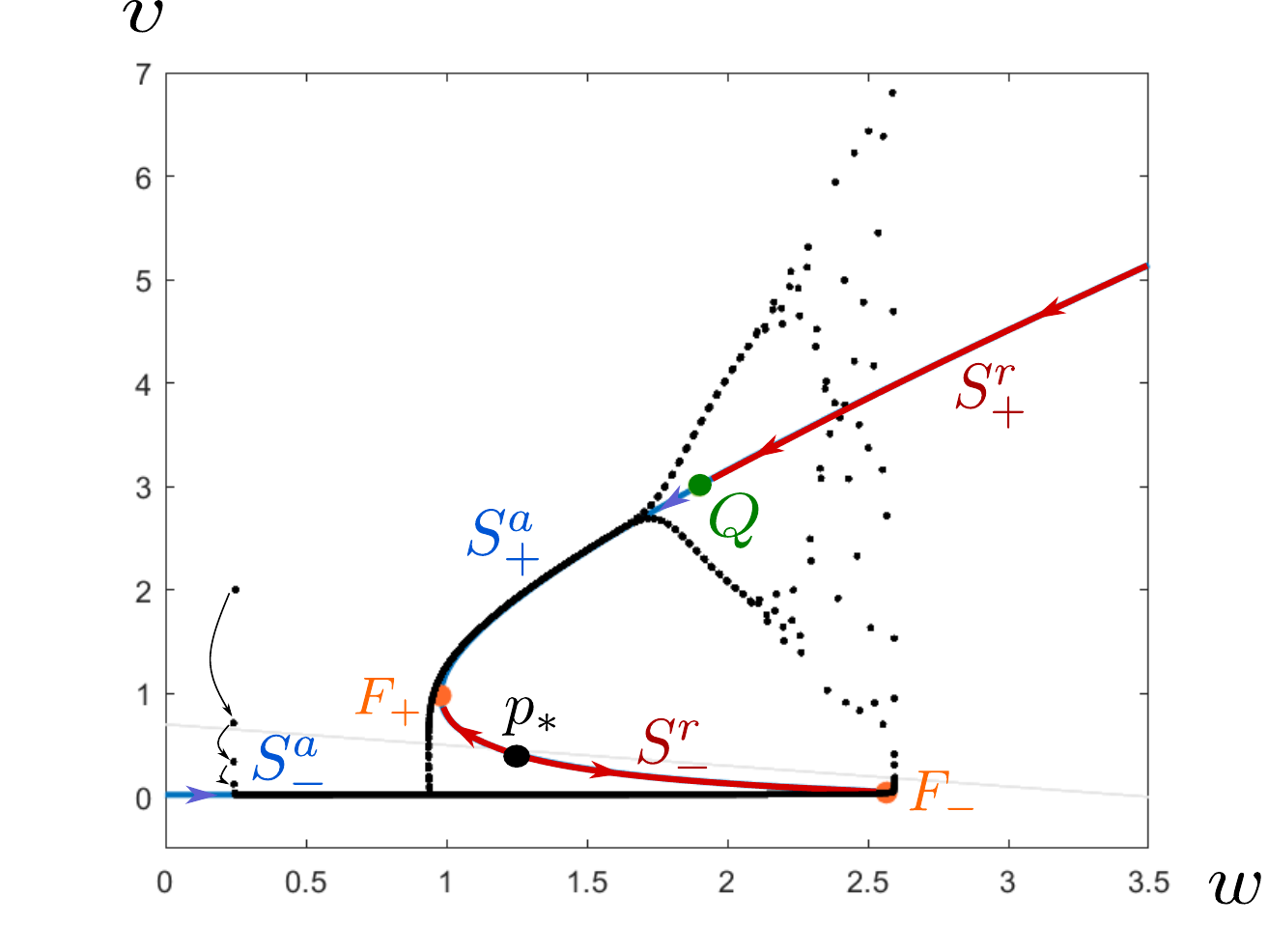}}
		\caption{Global dynamics of the fast-slow Chialvo map \eqref{eq:Chialvo} for $0 < \eps \ll 1$. Singular structure for $\eps = 0$ including critical manifolds $S_{\pm}^a$, $S_{\pm}^r$, the fold, flip and equilibrium points $F_{\pm}$, $Q$ and $p_\ast$ respectively are shown overlaid with numerically computed iterates of \eqref{eq:Chialvo} with initial condition $(w,v) = (1/4,2)$, $\eps = 10^{-3}$, $\tilde a = 1$, $\tilde b = 5$ and varying values of $\tilde c, k$. The `slow nullcline' $\{g(w,v,0) = 0\}$ is also shown in shaded grey. Dynamics away from the non-normally hyperbolic points $F_\pm$ and $Q$ is described by the theory in Sections \ref{sec:a_coordinate-independent_framework_for_fast-slow_maps} and \ref{sec:slow_manifold_theorems}. We identify four distinct behaviours depending on $(\tilde c, k)$, which are described in the text. (a) Case I: Excitability, with $(\tilde c,k) = (7,0.07)$. (b) Case II: Relaxation, with $(\tilde c,k) = (3.5,0.07)$. (c) Case III: Non-chaotic bursting, with $(\tilde c,k) = (3.5,0.035)$. (d) Case IV: Chaotic bursting, with $(\tilde c,k) = (3.5,0.02)$.}
		\label{fig:Chialvo_global_dynamics}
	\end{figure}

	\subsubsection*{Case I: Excitability}
	
	Excitable dynamics occurs for parameters such that $p_\ast \in S_+^a$ with $w_\ast = \varphi_0(v_\ast) < w_-$ so that it lies to the left of the lower fold point $F_-$, and $w_{flip} > w_-$ so that the flip point $Q$ lies to the right of $F_-$. Figure \ref{fig:Chialvo_global_dynamics}(a) shows the singular dynamics for a particular choice of parameter values satisfying these conditions, overlaid with iterates of the map \eqref{eq:Chialvo} for $\eps = 10^{-3}$.
	
	Starting from the point $(w,v) = (1/4,2)$, subsequent iterates are exponentially attracted to the attracting slow manifold $S^a_{-,\eps}$ along stable fibers according to Theorem \ref{thm:foliations}. After reaching a neighbourhood of $S_{-,\eps}^a$, iterates are well approximated by the map \eqref{eq:Chialvo_slow} (up to an error of $O(\eps^2)$), and move to the right for an $O(\eps^{-1})$ but finite number of iterates before reaching a neighbourhood of $F_-$. Unfortunately the dynamics near the fold point near $F_-$ are not covered by the normally hyperbolic theory in Sections \ref{sec:a_coordinate-independent_framework_for_fast-slow_maps}-\ref{sec:slow_manifold_theorems}. We therefore omit this part in our description; see however Remark \ref{rem:fold} below. After leaving a neighbourhood of $F_-$, iterates are exponentially attracted along stable fibers to $S^a_{+,\eps}$, before moving to the left as governed by the map \eqref{eq:Chialvo_slow}. After $O(\eps^{-1})$ but finitely many iterations, iterates reach (and remain thereafter in) a neighbourhood of the asymptotically stable fixed point $p_{\ast,\eps} \in S_{+,\eps}^a$ which perturbs from $p_\ast \in S_+^a$ as described by Proposition \ref{prop:Chialvo_fixed_points}.

	\begin{remark}
		\label{rem:fold}
		In particular cases, fold bifurcations in maps arising after discretization of ODE systems with a (continuous-time) fold bifurcation have been treated in detail using geometric techniques in e.g.~\cite{Nipp2013,Nipp2009}. In general, however, a detailed geometric description of the dynamics near the discrete fold bifurcation remains a topic of future work. 
	\end{remark}

	\subsubsection*{Case II: Relaxation}
	
	Relaxation-type dynamics occurs for parameter values such that $p_\ast \in S^r_-$ and $w_{flip} > w_-$, for which flip point $Q$ lies to the right of the lower fold point $F_-$. Figure \ref{fig:Chialvo_global_dynamics}(b) shows the singular dynamics for particular choice of parameter values, overlaid with iterates of the map \eqref{eq:Chialvo} for $\eps = 10^{-3}$.
	
	As before, we take an initial condition $(w,v) = (1/4,2)$. The geometric description of the dynamics is analogous to the excitable case I, up to the point where iterates reach a neighbourhood of $S^a_{+,\eps}$. In the relaxation case, iterates continues to track iterates of \eqref{eq:Chialvo_slow} until they reach a neighbourhood of the upper fold $F^+$. As with the lower fold $F_-$, the normally hyperbolic theory of Sections \ref{sec:a_coordinate-independent_framework_for_fast-slow_maps}-\ref{sec:slow_manifold_theorems} breaks down near $F_+$, and we omit a proper treatment of the dynamics here. Normally hyperbolic theory applies once more after iterates leave a neighbourhood of $F_+$, after which iterates are exponentially attracted along stable fibers to $S_{-,\eps}$. This provides the basic mechanism for relaxation-type dynamics observed in Figure \ref{fig:Chialvo_global_dynamics}(b).

	\subsubsection*{Case III: Non-chaotic bursting}
	
	Here we consider parameter values such that $p_\ast \in S^r_-$, $w_{flip} \in (w_+,w_-)$ and $|w_{flip} - w_-| < \delta$ is sufficiently small that the flip points $Q^2_{\pm}$ in the second iterate layer map $v \mapsto \bar v^2$ lie to the right of $F_-$ as in Figure \ref{fig:Chialvo_singular}. Figure \ref{fig:Chialvo_global_dynamics}(c) shows the singular dynamics for the same choice of parameter values as Figure \ref{fig:Chialvo_singular}, overlaid with iterates of the map \eqref{eq:Chialvo} for $\eps = 10^{-3}$.
	
	Starting from $(w,v) = (1/4,2)$, iterates approach $S_{-,\eps}^a$ and subsequently follow iterates of the map \eqref{eq:Chialvo_slow} up to a neighbourhood of $F_-$, similarly to cases I and II. After leaving a neighbourhood of $F_-$, iterates are exponentially attracted to an attracting slow manifold $S_{+,\eps}^{a,2}$ which perturbs
	from (a suitable submanifold of) the normally hyperbolic and attracting manifold of period-2 points $S^2$ in Figure \ref{fig:Chialvo_singular}; the existence and properties of the slow manifold $S_{+,\eps}^{a,2}$ and its stable manifold follow after applying the results in Section \ref{sec:slow_manifold_theorems} to the second iterate map $(w,v)^\textnormal{T} \mapsto (\bar w^2, \bar v^2)^\textnormal{T}$. Iterates move leftward along $S_{+,\eps}^{a,2}$ in a `flip-like' orientation-reversing manner about the repelling slow manifold $S_{+,\eps}^r$. They subsequently traverse a neighbourhood of the flip point $Q$, before being attracted to $S_{+,\eps}^a$ after a delay effect which is described in detail in \cite{Baesens1991}. Finally, iterates track iterates of the map \eqref{eq:Chialvo_slow} along $S_{+,\eps}^a$ towards to fold $F_+$, and return to a neighbourhood of $S_{-,\eps}^a$ after leaving a neighbourhood of $F_+$. Assuming a simple jump-type dynamics at $F_\pm$ (as is indeed observed numerically), the geometric sequence repeats, providing a mechanism for non-chaotic bursting.
	
	\begin{remark}
		Similar non-chaotic dynamics is expected as long as $w_{flip}$ and $w_-$ are sufficiently close, i.e.~$|w_{flip} - w_-| < \delta$ for $\delta > 0$ sufficiently small.
	\end{remark}

	\subsubsection*{Case IV: Chaotic bursting}
	
	Lastly, we consider parameter values such that $p_\ast \in S^r_-$, $w_{flip} \in (w_+,w_-)$ and $|w_{flip} - w_-| > \delta$, where the latter condition implies that the layer map is chaotic for $w-$values close to $w_-$. Figure \ref{fig:Chialvo_global_dynamics}(d) shows the singular dynamics for particular choice of parameter values satisfying these conditions, overlaid with iterates of the map \eqref{eq:Chialvo} for $\eps = 10^{-3}$.
	
	Starting again from $(w,v) = (1/4,2)$, iterates are initially attracted to $S_{-,\eps}^a$ and subsequently track iterates of the map \eqref{eq:Chialvo_slow} towards the fold $F_-$, similarly to cases I, II and III. After leaving a neighbourhood of the fold, complicated dynamics corresponding to a slow (leftward) drift through part of a period-doubling cascade occurs. The details of this transition are beyond the scope of this work, however the reader is referred to \cite{Trujillo2021} for details on the period-doubling cascade in the corresponding layer problem, and to \cite{Baesens1991} for details on the slow drift through a period-doubling cascade in a particular `normal form'. Once iterates are close enough to $w_{flip}$, the dynamics becomes non-chaotic and orientation-reversing as in case III, and eventually approach a neighbourhood of $S_{+,\eps}^a$ after a delay effect (again described by \cite{Baesens1991}). The subsequent tracking of iterates along $S_{+,\eps}^a$ and the observed return to $S_{-,\eps}^a$ is similar to that described for cases II and III. Assuming a simple-jump like behaviour near the folds $F_\pm$, this geometric sequence repeats itself and provides a mechanism for chaotic bursting behaviour.

	\subsection{Euler-discretized fast-slow ODEs in non-standard form}
	\label{sub:Euler_discretization}
	
	We consider the existence and properties of invariant manifolds in Euler-discretizations of singularly perturbed ODEs in the general non-standard form
	\begin{equation}
		\label{eq:general_ode}
		z' = N(z) f(z) + \eps G(z, \eps) , \qquad z \in \mathbb R^n , \qquad 0 < \eps \ll 1 ,
	\end{equation}
	where the right-hand-side is $C^r-$smooth in $(z,\eps)$ and assumed to have a normally hyperbolic critical manifold for $\eps = 0$. More precisely, we impose the following assumptions:
	
	\begin{assumption}
		\label{ass:ODEs_1}
		\textup{(Restriction to fast-slow ODEs)} The ODE \eqref{eq:general_ode} satisfies the geometric definition of a $C^r-$smooth singularly perturbed ODE in e.g.~\cite{Fenichel1979,Wechselberger2019}, i.e.~the level set
		\[
		C_{ode} := \left\{ z \in \mathbb R^n : N(z) f(z) = O_n \right\}
		\]
		contains a $(k<n)-$dimensional regularly embedded critical manifold
		\[
		S_{ode} := \left\{ z \in \mathbb R^n : f(z) = O_{n-k} \right\} .
		\]
		We assume for simplicity that $S \subseteq C \subset \mathbb R^n$ is connected.
	\end{assumption}
	
	\begin{assumption}
		\label{ass:ODEs_2}
		\textup{(Regularity of the matrix $N(z)$)} We assume that the $n \times (n-k)$ matrix $N(z)$ has full column rank for all $z \in S_{ode}$, and that equilibria of $N(z)f(z)$ for all $z \notin S_{ode}$, if they exist, are isolated.
	\end{assumption}
	
	Assumptions \ref{ass:ODEs_1}-\ref{ass:ODEs_2} are the defining assumptions for fast-slow ODEs in \cite{Wechselberger2019}, which are directly analogous to the defining Assumptions \ref{ass:1}-\ref{ass:factorisation} for fast-slow maps in Section \ref{sec:a_coordinate-independent_framework_for_fast-slow_maps}.
	
	In the following we are concerned with the dynamics of maps induced by Euler discretization of general fast-slow ODEs \eqref{eq:general_ode} under Assumptions \ref{ass:ODEs_1}-\ref{ass:ODEs_2}. Specifically, Euler discretization in time yields the map
	\begin{equation}
		\label{eq:discretized_map_Euler}
		z \mapsto \bar z = z + h N(z) f(z) + \eps h G(z,\eps) ,
	\end{equation}
	where $h > 0$ is the associated step-size. It is straightforward to verify that the map \eqref{eq:discretized_map_Euler} has critical manifold $S = \{z \in \mathbb R^n : f(z) = O_{n-k}\}$ and satisfies Assumptions \ref{ass:1}-\ref{ass:factorisation}; this follows directly from Assumptions \ref{ass:ODEs_1}-\ref{ass:ODEs_2} on the ODE \eqref{eq:general_ode}. Hence, the formalism of Section \ref{sec:a_coordinate-independent_framework_for_fast-slow_maps} applies.
	
	\begin{remark}
		All results obtained for the map \eqref{eq:discretized_map_Euler} extend and are similar in nature to their pre-existing counterparts in the context of Euler-discretized fast-slow ODEs in standard form
		\begin{equation}
			\label{eq:stnd_ode}
			\begin{split}
				x' &= \eps \tilde g(x,y,\eps) , \\
				y' &= \tilde f(x,y,\eps) ,
			\end{split}
			\qquad (x,y) \in \mathbb R^k \times \mathbb R^{n-k} , 
			\qquad 0 < \eps \ll 1 ,
		\end{equation}
		since this system can be written in the general form \eqref{eq:general_ode} by setting $z = (x,y)^\textnormal{T}$ and
		\[
		N(z) = 
		\begin{pmatrix}
			O_{k,n-k} \\
			I_{n-k}
		\end{pmatrix}
		, \qquad 
		f(z) = \tilde f(x,y,0) , 
		\qquad
		G(z,\eps) = 
		\begin{pmatrix}
			\tilde g(x,y,\eps) \\
			\tilde f_{rem}(x,y,\eps)
		\end{pmatrix}
		,
		\]
		where $\tilde f_{rem}(x,y,\eps) = \tilde f(x,y,\eps) - \tilde f(x,y,0)$. See \cite{Nipp1995,Nipp1996} and the references therein for more on existing results for discretizations of fast-slow ODEs in standard form \eqref{eq:stnd_ode}.
	\end{remark}

	\subsubsection*{Layer map}
	
	The layer map is given by
	\begin{equation}
		\label{eq:layer_discretization_Euler}
		z \mapsto \bar z = z + h N(z) f(z) .
	\end{equation}
	By Proposition \ref{prop:EVs} and Corollary \ref{cor:nh}, the $n-k$ non-trivial multipliers $\mu_j(z)$ along the critical manifold $S$ are given by
	\[
	\mu_j(z) = 1 + \lambda_j(z) , \qquad i = 1, \ldots , n - k ,
	\]
	where $\lambda_j(z)$ are eigenvalues of the matrix
	\[
	h Df(z) N(z) , \qquad z \in S .
	\]
	For the ODE, non-trivial eigenvalues $\lambda_{j,ode}(z)$ along $S_{ode}$ are given by eigenvalues to the matrix $Df(z)N(z)$. It follows that
	\begin{equation}
		\label{eq:multipliers_Euler}
		\mu_j(z) = 1 + \lambda_j(z) = 1 + h \lambda_{j,ode}(z) .
	\end{equation}
	We obtain the following stability conditions for $S$, depending on the real part of the ODE eigenvalues along $S_{ode}$ and the step-size $h$ associated to the discretization.
	
	\begin{lemma}
		\label{lem:nh_Euler}
		Consider the layer map \eqref{eq:layer_discretization_Euler} with step-size $h > 0$. The critical manifold $S$ is normally hyperbolic at $z \in S$ if
		\[
		2 \textup{Re}\, \lambda_{j,ode}(z) \neq - h | \lambda_{j,ode}(z) |^2 ,
		\]
		for all $j = 1, \ldots , n-k$. A normally hyperbolic point $z \in S$ is attracting if
		\[
		2 \textup{Re}\, \lambda_{j,ode}(z) < - h | \lambda_{j,ode}(z) |^2 ,
		\]
		for all $j = 1, \ldots , n-k$, repelling if
		\[
		2 \textup{Re}\, \lambda_{j,ode}(z) > - h | \lambda_{j,ode}(z) |^2 ,
		\]
		for all $j = 1, \ldots , n-k$, and saddle-type otherwise.
	\end{lemma}
	
	\begin{proof}
		This follows directly from equation \eqref{eq:multipliers_Euler} and Corollary \ref{cor:nh}.
	\end{proof}
	
	It follows that normal hyperbolicity breaks down at fold-type singularities $z \in S$ where $\lambda_{j,ode}(z) = 0$ for some $j = 1, \ldots , n-k$, independently of $h$. On the other hand, Hopf-type non-normally hyperbolic singularities in the ODE have eigenvalues such that $\textup{Im}\,\lambda_{j,ode}(z) \neq 0$, and their location and may be shifted by $O(h)$ distances in the discretization. The discretization also introduces additional singularities that are not present in the ODE, with multipliers that cross through the left hand side of the unit circle in $\{ z \in \mathbb C : \textup{Im}\, z \leq 0 \}$, however the location of the corresponding singularities depends on $h$, and tend to infinity in the limit $h \to 0^+$.

	\subsubsection*{Reduced map}
	
	We now consider the reduced map for the Euler discretization \eqref{eq:discretized_map_Euler}. If we assume that $S_{ode}$ is normally hyperbolic, it follows from Lemma \ref{lem:nh_Euler} that for sufficiently small step-sizes $h>0$, compact submanifolds of $S$ are also normally hyperbolic since any non-hyperbolic points on $S$ induced by the discretization are pushed to infinity as $h \to 0^+$. Therefore, the oblique projection operator in \eqref{eq:proj_def} and \eqref{eq:proj} is well-defined on compact submanifolds of $S$ for sufficiently small $h > 0$. In fact, the projection is invariant with respect to the discretization in the sense that
	\begin{equation}
		\label{eq:proj_Euler}
		\Pi_{\mathcal N}^{S} = I_n - h N h^{-1} (Df N)^{-1} Df |_S =  I_n - N (Df N)^{-1} Df |_S = \Pi_{\mathcal N,ode}^{S} ,
	\end{equation}
	where $\Pi_{\mathcal N,ode}^{S}$ denotes the projection operator used to define the reduced problem for the ODE \eqref{eq:general_ode}, see \cite[Def.~3.8]{Wechselberger2019}. It follows by Proposition \ref{prop:reduced_map} that the reduced map is
	\[
	z \mapsto \bar z = z + \eps h \Pi_{\mathcal N}^{S} G(z,0) |_S = z + \eps h \Pi_{\mathcal N,ode}^{S} G(z,0)|_S ,
	\]
	with $\Pi_{\mathcal N}^{S} = \Pi_{\mathcal N,ode}^{S}$ given by \eqref{eq:proj_Euler}. 
	As noted already in Remark \ref{rem:reduced_map_discretizations}, in this case it is possible to rescale the step-size via $h = \eps^{-1} \tilde h$ in order to obtain a well-defined reduced map on $S$ in the limit $\eps \to 0$. Explicitly, we obtain
	\[
	z \mapsto \bar z = z + \eps h \Pi_{\mathcal N}^{S} G(z,0) |_S = z + \tilde h \Pi_{\mathcal N,ode}^{S} G(z,0)|_S ,
	\]
	which is $\eps -$independent when formulated in terms of the rescaled step-size parameter $\tilde h = \eps h$.

	\subsubsection*{Dynamics for $0 < \eps \ll 1$}
	
	The results of Section \ref{sec:slow_manifold_theorems} apply to the Euler discretization \eqref{eq:discretized_map_Euler} under conditions of normally hyperbolicity. Moreover, Lemma \ref{lem:nh_Euler} shows that the normally hyperbolicity of the critical manifold $S$ after discretization is entirely determined by the non-trivial eigenvalues along the critical manifold of the ODE \eqref{eq:general_ode} and the step-size $h > 0$. This allows for the statements on the existence and properties of slow manifolds and corresponding stable/unstable manifolds in the discretized system, purely in terms of ODE properties and the step-size parameter $h$.
	
	It will be helpful for the statement of our results to introduce additional notation for spectral bounds in the ODE. Given a $k-$dimensional normally hyperbolic critical manifold $S_{ode}$, denote the $n_a \leq n-k$ non-trivial eigenvalues such that $\textrm{Re}\, \lambda_{j,ode}(z) < 0$ by $\lambda_{a,j,ode}(z)$, and the $n_r = n-k-n_a$ non-trivial eigenvalues such that $\textrm{Re}\, \lambda_{j,ode}(z) > 0$ by $\lambda_{r,j,ode}(z)$. We define the corresponding spectral bounds by
	\begin{equation}
		\label{eq:ode_spectral_bounds}
		\nu_{A,ode} := \sup_{z \in S_{ode},j=1,\ldots,n_a} | \textrm{Re}\, \lambda_{a,j,ode}(z) | , \qquad 
		\nu_{R,ode} := \inf_{p \in S_{ode},j=1,\ldots,n_r} | \textrm{Re}\, \lambda_{r,j,ode}(z) | .
	\end{equation}
	
	Our main results for the Euler discretization \eqref{eq:discretized_map_Euler} are as follows.
	
	\begin{thm}
		\label{thm:Euler}
		\textup{(Existence and characterisation of slow and stable/unstable manifolds in Euler discretized ODEs)} Assume that the fast-slow ODE system \eqref{eq:general_ode} satisfies Assumptions \ref{ass:ODEs_1}-\ref{ass:ODEs_2}, with compact, connected and normally hyperbolic critical manifold $S_{ode}$. Then for the Euler discretization \eqref{eq:discretized_map_Euler} with $\eps=0$ we have $S = S_{ode}$ and coincident transverse linear fiber bundles $\mathcal N = \mathcal N_{ode}$. 
		Moreover, there exists an $\eps_0 > 0$ and a $h_0 > 0$ such that for all $\eps \in (0,\eps_0)$ and $h \in (0,h_0)$ the following assertions are true: 
		\begin{enumerate}
			\item[(i)] $S$ and $S_{ode}$ perturb to nearby $C^r-$smooth locally invariant slow manifolds $S_{\eps,h}$ and $S_{ode,\eps}$ respectively, such that
			\begin{equation}
				\label{eq:distance}
				\textup{dist}(S_{\eps,h} , S_{ode,\eps} ) = O(\eps^2 h) .
			\end{equation}
			For each fixed $h$ the manifold $S_{\eps,h}$ has all the properties described in Theorem \ref{thm:slow_manifolds}, and the manifold $S_{ode,\eps}$ is a Fenichel slow manifold.
			\item[(ii)] $S_{\eps,h}$ and $S_{ode,\eps}$ have local graph representations
			\[
			S_{\eps,h} = \left\{ (x, \varphi_{\eps,h}(x) ) : x \in \mathcal K \right\} , \qquad 
			S_{ode,\eps} = \left\{ (x, \varphi_{ode,\eps}(x) ) : x \in \mathcal K \right\} ,
			\]
			where
			\begin{equation}
				\label{eq:varphi_h_ode}
				\begin{split}
					\varphi_{\eps,h}(x) &= \varphi_0(x) - \eps (D_yf)^{-1} (Df N)^{-1} Df G |_S + O(\eps^2 h) , \\
					\varphi_{ode,\eps}(x) &= \varphi_0(x) - \eps (D_yf)^{-1} (Df N)^{-1} Df G |_{S_{ode}} + O(\eps^2) .
				\end{split}
			\end{equation}
			\item[(iii)] $W^{s/u}_{loc}(S)$ and $W^{s/u}_{loc}(S_{ode})$ perturb to $C^r-$smooth stable/unstable manifolds $W^{s/u}_{loc}(S_{\eps,h})$ and $W^{s/u}_{loc}(S_{\eps, ode})$ respectively, such that
			\begin{equation}
				\label{eq:distance_stable_manifolds}
				\textup{dist}(W^{s/u}_{loc}(S_{\eps,h}) , W^{s/u}_{loc}(S_{ode,\eps}) ) = O(\eps^2 h) .
			\end{equation}
			For each fixed $h$, the manifolds $W^{s/u}_{loc}(S_{\eps,h})$ have all the properties described in Theorem \ref{thm:stable_manifolds}, while the manifolds $W^{s/u}_{loc}(S_{ode,\eps})$ are described by Fenichel theory \cite{Fenichel1979,Jones1995,Kuehn2015,Wiggins2013}.
			\item[(iv)] The manifolds $W^{s/u}_{loc}(S_{\eps,h})$ and $W^{s/u}_{loc}(S_{\eps, ode})$ admit foliations by stable/unstable fibers $w^{s/u}_{loc}(S_{\eps,h})$ and $w^{s/u}_{loc}(S_{\eps, ode})$ respectively. For each fixed $h$, the foliations of $W^{s/u}_{loc}(S_{\eps,h})$ are described by Theorem \ref{thm:foliations}, while the foliations of $W^{s/u}_{loc}(S_{ode,\eps})$ are described by Fenichel theory. 
			Contraction along stable fibers $w^{s}_{loc}(S_{\eps,h})$ occurs with rate $\chi_A(h)$ fixed and bounded in the interval
			\begin{equation}
				\label{eq:Euler_contraction}
				1 - \nu_{A,ode} h + c_A h^2 < \chi_A(h) < 1 ,
			\end{equation}
			for a constant $c_A$, and repulsion along unstable fibers $w^{u}_{loc}(S_{\eps,h})$ occurs with rate $\chi_R(h)$ fixed and satisfying
			\begin{equation}
				\label{eq:Euler_expansion}
				\chi_R(h) > 1 + \nu_{R,ode} h + c_R h^2 > 1 ,
			\end{equation}
			for a constant $c_R$.
		\end{enumerate}
	\end{thm}
	
	\begin{proof}
		The equality $S = S_{ode}$ is immediate from the form of the layer map \eqref{eq:layer_discretization_Euler}. Coincidence of the transverse linear fiber bundles $\mathcal N = \cup_{z \in S} \mathcal N_z$ and $\mathcal N_{ode} = \cup_{z \in S_{ode}} \mathcal N_{ode,z}$ follows from the pointwise equality
		\[
		\mathcal N_{z,ode} = \text{span} \{N^i(z) : 1,\ldots,n-k \} = \text{span} \{h N^i(z) : 1,\ldots,n-k \} = \mathcal N_z ,
		\]
		where $N^i(z)$ are the columns of $N(z)$, which holds for all $h>0$ and $z \in S_{ode} = S$.
		
		\
		
		Except for the statement pertaining to the distance between slow manifolds in \eqref{eq:distance}, statement (i) follows directly from Fenichel theory \cite{Fenichel1979,Jones1995,Kuehn2015,Wiggins2013} and Theorem \ref{thm:slow_manifolds}. The distance in \eqref{eq:distance} follows from statement (ii) after extending the local estimate
		\[
		| \varphi_{\eps,h}(x) - \varphi_{ode,\eps}(x) | = O(\eps^2 h) 
		\]
		obtained from the equations in \eqref{eq:varphi_h_ode} via compactness and a partition of unity. Thus it suffices to prove the local graph representations in \eqref{eq:varphi_h_ode}. The expansion for $\varphi_{ode,\eps}(x)$ follows directly from \cite[eqns~(3.33), (3.36)]{Wechselberger2019}, and the expansion for $\varphi_{\eps,h}(x)$ up to and including $O(\eps)$ follows directly from Theorem \ref{thm:graph_slow_manifolds}. Carrying out the matching argument in the proof of Proposition \ref{prop:reduced_map} for the map \eqref{eq:discretized_map_Euler} and expanding in powers of $\eps h$ instead of $\eps$ leads to the higher order correction $O(\eps^2 h)$ (a factor of $h$ cancels out since $h$ is a common factor of $\bar z - z$).
		
		\
		
		Except for the statement pertaining to the distance between slow manifolds in \eqref{eq:distance_stable_manifolds}, statement (iii) follows directly from Fenichel theory \cite{Fenichel1979,Jones1995,Kuehn2015,Wiggins2013} and Theorem \ref{thm:stable_manifolds}. The estimate \eqref{eq:distance_stable_manifolds} follows from the following facts:
		\begin{itemize}
			\item $W_{loc}^{s/u}(S)$ and $W_{loc}^{s/u}(S_{ode})$ are linearly tangent along $S = S_{ode}$ since $\mathcal N = \mathcal N_{ode}$;
			\item By assertion (i), $\text{dist}(S_{\eps,h},S_{ode,\eps}) = O(\eps^2 h)$; 
			\item $S_{ode,\eps} \subseteq W_{loc}^{s/u}(S_{ode,\eps})$ and $S_{\eps,h} \subseteq W_{loc}^{s/u}(S_{\eps,h})$, since both slow manifolds are contained within the intersection of their corresponding perturbed stable/unstable manifolds by Theorem \ref{thm:stable_manifolds}.
		\end{itemize}
		
		\
		
		It remains to prove assertion (iv). Statements pertaining to foliations of $W^{s/u}_{loc}(S_{ode,\eps})$ follows directly from Fenichel theory \cite{Fenichel1979,Jones1995,Kuehn2015,Wiggins2013}. Those pertaining to foliations of $W^{s/u}_{loc}(S_{\eps,h})$ for fixed $h$ follow from Theorem \ref{thm:foliations}. To see that the bounds on the contraction rate $\chi_A(h)$ are given by \eqref{eq:Euler_contraction}, notice that by \eqref{eq:multipliers_Euler} we have
		\[
		|\mu_{a,j}(z)| = | 1 + h \lambda_{a,j,ode}(z) | = 1 + h \textrm{Re}\, \lambda_{a,j,ode}(z) + O(h^2) ,
		\]
		as $h \to 0^+$, for all $z \in S = S_{ode}$. Combining this with \eqref{eq:ode_spectral_bounds} and $\nu_A = \sup_{z \in S , j = 1, \ldots , n_a} |\mu_{a,j}(z)|$ yields the bounds on the contraction rate in \eqref{eq:Euler_contraction} via Theorem \ref{thm:foliations} (iii). Similar arguments yield the bounds for the expansion rate $\chi_R(h)$ in \eqref{eq:Euler_expansion} via Theorem \ref{thm:foliations} (iv).
	\end{proof}
	
	\begin{remark}
		The latter part of assertion (iii) as well as assertion (iv) in Theorem \ref{thm:Euler} hold only for fixed $h>0$ and asymptotically small $\eps \sim 0$. They do not account for asymptotic dependence in the case that both $\eps \to 0$ and $h \to 0$, in which case one has an instance of the general class of double singular limit problems described in \cite{Kuehn2022}. 
	\end{remark}
	

	%
	%
	%
	%
	%
	%
	%

	\subsection{Fast-slow Poincar\'e maps}
	\label{sub:Poincare_maps}
	
	Poincar\'e maps associated to parameter-dependent systems of ODEs with hyperbolic limit cycles naturally give rise to fast-slow Poincar\'e maps if the parameter is allowed to evolve slowly in time. In this section we apply the formalism and results of Sections \ref{sec:a_coordinate-independent_framework_for_fast-slow_maps}-\ref{sec:slow_manifold_theorems} in order to study the dynamics of such systems. We begin by characterising the dynamics of the Poincar\'e map, before considering the implications in the corresponding (higher dimensional) ODE.
	

	%

	
	\
	
	Consider the ODE system
	\begin{equation}
		\label{eq:stnd_fast_slow_limit_cycles}
		\begin{split}
			z' &= \tilde f(z,\alpha,\eps) , \\
			\alpha' &= \eps \tilde g(z,\alpha,\eps) ,
		\end{split}
	\end{equation}
	where $(z,\alpha) \in \mathbb R^n \times \mathbb R$, $0 < \eps \ll 1$ and the right-hand-side 
	is $C^{r \geq 1}-$smooth in $(z,\eps)$. 
	System \eqref{eq:stnd_fast_slow_limit_cycles} is in the standard form for fast-slow systems, with $n$ fast variables $z$ and a single slow variable $\alpha$. For the systems considered here, however, we do not require the existence of a critical manifold $\{(z,\alpha) \in \mathbb R^{n+1} : f(z,\alpha)=0\}$. Thus \eqref{eq:stnd_fast_slow_limit_cycles} is not necessarily slow-fast in the geometric sense of \cite{Fenichel1979}, see also \cite[Definition 3.2]{Wechselberger2019}. Nevertheless, we shall refer to the limiting system
	\begin{equation}
		\label{eq:stnd_limit_cycles_layer}
		\begin{split}
			z' &= \tilde f(z,\alpha,0) , \\
			\alpha' &= 0 ,
		\end{split}
	\end{equation}
	as $\eps \to 0$ as the \textit{layer problem} for system \eqref{eq:stnd_fast_slow_limit_cycles}. 
	Our key assumption is the following.
	
	\begin{assumption}
		\label{ass:limit_cycles}
		\textup{(Existence of a hyperbolic limit cycle)} 
		The $\alpha-$family of ODEs defined by the layer problem \eqref{eq:stnd_limit_cycles_layer} has a hyperbolic limit cycle $\Gamma_{\alpha_\ast}$ of period $T_{\alpha_\ast}$ passing through the point $z = z_\ast$ for $\alpha = \alpha_\ast$.
	\end{assumption}
	
	Since $\Gamma_{\alpha_\ast}$ is hyperbolic, it persists for nearby $\alpha$ in a sufficiently small open interval $V_\alpha \ni \alpha_\ast$ in $\mathbb R$, i.e.~the layer problem \eqref{eq:stnd_limit_cycles_layer} has a hyperbolic limit cycle $\Gamma_{\alpha}$ of period $T_\alpha$ for each fixed $\alpha \in V_\alpha$. It follows that
	\[
	\mathcal M := \left\{ \Gamma_\alpha : \alpha \in V_\alpha \right\} \subseteq \mathbb R^n \times V_\alpha
	\]
	defines a $2-$dimensional manifold of limit cycles topologically equivalent to a cylinder in $\mathbb R^n \times V_\alpha$.
	
\begin{figure}[t!]
		\centering
		\includegraphics[scale=0.13]{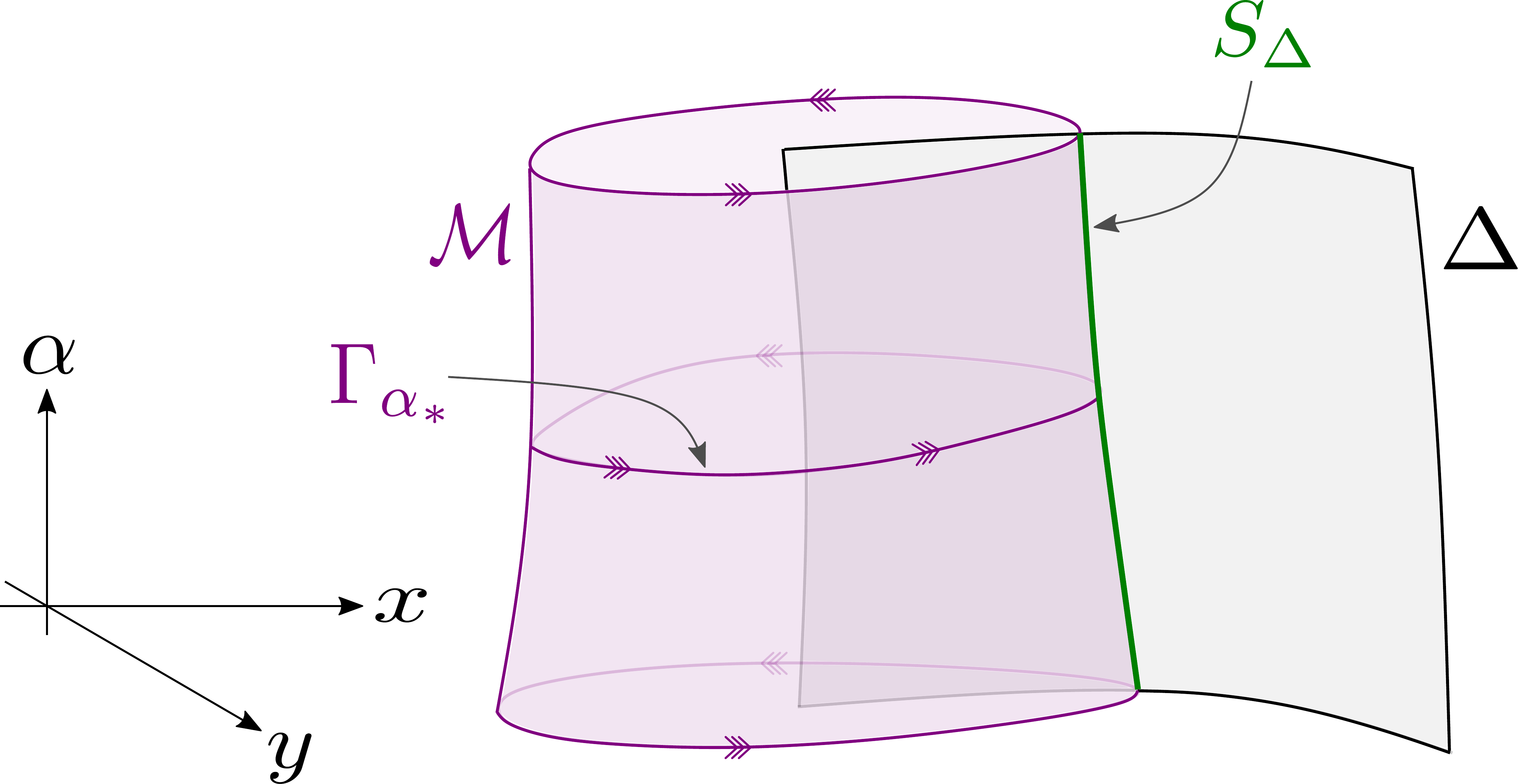}
		\caption{Setup for $\eps = 0$. The layer problem \eqref{eq:stnd_limit_cycles_layer} has a manifold of hyperbolic limit cycles $\mathcal M$, shown here in shaded purple, which contains the limit cycle $\Gamma_{\alpha_\ast}$ and extends over the open $\alpha-$interval $V_\alpha \ni \alpha_\ast$. A Poincar\'e map $P_\Delta : \Delta \to \Delta$ on the cross-section $\Delta$ shown in shaded green is induced by the flow of the layer problem \eqref{eq:stnd_limit_cycles_layer}. Hyperbolicity of the limit cycles $\Gamma_\alpha$ implies that the 1-dimensional submanifold $S_\Delta = \mathcal M \cap \Delta$ (bold green) is a normally hyperbolic critical manifold for $P_\Delta$, see Lemma \ref{lem:Poincare_map}.}
		\label{fig:Poincare_setup}
	\end{figure}

	In the following we describe the Poincar\'e map $\mathcal P_\Delta : \Delta \to \Delta$ induced by the flow of the perturbed system \eqref{eq:stnd_fast_slow_limit_cycles} on a fixed codimension$-1$ cross-section $\Delta$ which is assumed to satisfy the usual non-degeneracy conditions with respect to $\mathcal M$, i.e.~
	\begin{itemize}
		\item $S_\Delta := \Delta \cap \mathcal M$ defines a $1-$dimensional submanifold of $\mathbb R^n \times V_\alpha$;
		\item $T_p \mathcal M \oplus T_p \Delta = T_p \mathbb R^{n+1}$ at each $p \in S_\Delta$;
	\end{itemize}
	see Figure \ref{fig:Poincare_setup}. 
	Without loss of generality, coordinates $z = (x,y)^\textnormal{T} \in \mathbb R^{n-1} \times \mathbb R$ and $\tilde f = (\tilde f_1, \tilde f_2)^\textnormal{T}$ can be chosen such that system \eqref{eq:stnd_fast_slow_limit_cycles} can be rewritten as
	\begin{equation}
		\label{eq:stnd_fast_slow_limit_cycles_2}
		\begin{split}
			x' &= \tilde f_1(x,y,\alpha,\eps) , \\
			y' &= \tilde f_2(x,y,\alpha,\eps) , \\
			\alpha' &= \eps \tilde g(x,y,\alpha,\eps) ,
		\end{split}
	\end{equation}
	and the section $\Delta$ can be written as a graph
	\begin{equation}
		\label{eq:Delta}
		\Delta = \left\{ (x, Y(x,\alpha), \alpha) : x \in V_x, \alpha \in V_\alpha \right\} ,
	\end{equation}
	where $V_x$ is a sufficiently small neighbourhood of the point $x_\ast$ defined by $z_\ast = (x_\ast,y_\ast)$, and $Y : \mathbb R^{n-1} \times \mathbb R \to \mathbb R$ is $C^r-$smooth function. It follows that $S_\Delta$ has graph representation
	\begin{equation}
		\label{eq:S_Delta}
		S_\Delta = \left\{ (\varphi_0(\alpha), Y(\varphi_0(\alpha)), \alpha) : \alpha \in V_\alpha \right\} ,
	\end{equation}
	for a $C^r-$smooth function $\varphi_0 : \mathbb R \to \mathbb R^{n-1}$.
	
	\begin{remark}
		In the following we shall often consider $S_\Delta$ only in $(x,\alpha)-$coordinates on $\Delta$, writing simply $S_\Delta = \{(\varphi_0(\alpha),\alpha) : \alpha \in V_\alpha \}$, while still identifying it with the intersection $\Delta \cap \mathcal M$. The existence of the parameterization in \eqref{eq:S_Delta} prevents this slight abuse of notation from causing confusion.
	\end{remark}
	
	We have the following result.
	
	
	\begin{lemma}
		\label{lem:Poincare_map}
		The Poincar\'e map $\mathcal P_\Delta$ is a fast-slow map satisfying Assumptions \ref{ass:1}-\ref{ass:factorisation} with normally hyperbolic critical manifold $S_\Delta$. 
	\end{lemma}
	
	\begin{proof}
		The structure of the singular set and its spectrum is invariant with respect to the choice of cross-section $\Delta$, as well as the choice of local coordinates on $\Delta$; see e.g.~\cite[Lemma 1.3]{Kuznetsov2013}. Thus 
		it suffices to show the result for a particular choice of $\Delta$. The derivation of the Poincar\'e map is similar to derivation in the context 3-dimensional systems in \cite{Stiefenhofer1998}. 
		%
		
		We work with system \eqref{eq:stnd_fast_slow_limit_cycles_2} and assume without loss of generality that
		\[
		z_\ast = (x_\ast,0)^\textnormal{T} , \qquad \tilde f_2(x_\ast,0,\alpha_\ast,0) \neq 0 ,
		\]
		i.e.~that the limit cycle $\Gamma_{\alpha_\ast}$ in the layer problem \eqref{eq:stnd_limit_cycles_layer} intersects the $n-$dimensional hyperplane 
		\[
		\Sigma := \left\{ (x,0,\alpha) : x \in V_x, \alpha \in V_\alpha \right\} ,
		\] 
		transversally at $(x_\ast,0,\alpha_\ast)$, for a suitable (sufficiently small) neighbourhood $x_\ast \ni V_x \subset \mathbb R^{n-1}$. 
		The neighbourhoods $V_x$ and $V_\alpha$ can be chosen such that $\Sigma$ is a valid choice of cross-section, i.e.~such that
		the choice $\Delta = \Sigma$ satisfies the non-degeneracy conditions for a cross-section given in bullet points prior to the statement of Lemma \ref{lem:Poincare_map}. Note that with this choice, $Y(x,\alpha) \equiv 0$ on $V_x \cap V_\alpha$. 
		
		System \eqref{eq:stnd_fast_slow_limit_cycles_2} induces a flow
		\[
		\Phi = (\Phi_1, \Phi_2, \Phi_3) : \mathbb R^{n-1} \times \mathbb R \times \mathbb R \times [0,\eps_0) \times \mathbb R \to \mathbb R^{n-1} \times \mathbb R \times \mathbb R
		\]
		which by Assumption \ref{ass:limit_cycles} satisfies
		\[
		\Phi_2 (x_\ast, 0, \alpha_\ast, 0, T_{\alpha_\ast}) = 0 \qquad \text{and} \qquad
		\Phi_2 ' (x_\ast, 0, \alpha_\ast, 0, T_{\alpha_\ast}) \neq 0 .
		\]
		Therefore, the implicit function theorem implies the existence of a smooth and (locally) invertible function $\tilde t : V_x \times V_y \times V_\alpha \times [0,\eps_0) \to \mathbb R$ such that $\tilde t(x_\ast, 0, \alpha_\ast, 0) = T_{\alpha_\ast}$ and
		\[
		\Phi_2(x,y,\alpha,\eps,\tilde t(x,y,\alpha,\eps)) \equiv 0 
		\]
		on $V_x \times V_y \times V_\alpha \times [0,\eps_0)$, where $V_y$ is a neighbourhood about $y_\ast = 0$ in $\mathbb R$. In terms of earlier notation, a limit cycle $\Gamma_\alpha \subset \mathcal M$ such that $\Gamma_\alpha \cap S_\Sigma = \{(\varphi_0(\alpha) , 0, \alpha)\}$ has period $T_\alpha = \tilde t(\varphi_0(\alpha), 0 , \alpha, 0)$.
		
		The preceding argument guarantees the existence of a well-defined Poincar\'e map $\mathcal P_\Sigma : \Sigma \to \Sigma$ 
		such that $\mathcal P_\Sigma(S_\Sigma) = S_\Sigma$. Explicitly, 
		\begin{equation}
			\label{eq:Poincare_map}
			\mathcal P_\Sigma : 
			\begin{pmatrix}
				x \\
				\alpha
			\end{pmatrix}
			\mapsto
			\begin{pmatrix}
				\bar x \\
				\bar \alpha
			\end{pmatrix}
			=
			\begin{pmatrix}
				x \\
				\alpha
			\end{pmatrix}
			+
			\begin{pmatrix}
				\widehat f(x,\alpha,\eps) \\
				\eps g(x,\alpha,\eps)
			\end{pmatrix}
			,
		\end{equation}
		with
		\begin{equation}
			\label{eq:Poincare_map_2}
			\begin{split}
				\widehat f(x,\alpha,\eps) &:= \Phi_1(x,0,\alpha,\eps,\tilde t(x,0,\alpha,\eps)) - x , \\
				g(x,\alpha,\eps) &:= 
				D_\eps \Phi_3(x,0,\alpha,0,\tilde t(x,0,\alpha,0)) + O(\eps) ,
			\end{split}
		\end{equation}
		where the latter expression follows after expansion in $\eps$. In particular, the function $g(x,\alpha,\eps)$ can be expressed in terms of the original function $\tilde g(z,\alpha,\eps) = \tilde g(x,y,\alpha,\eps)$ via the integral formula
		\[
		g(x,\alpha,\eps) = \int_0^{\tilde t(x,0,\alpha,0)} \tilde g \left( \Phi_1(x,0,\alpha,0,s), 0, \alpha, 0 \right) ds + O(\eps) .
		\]
		
		
		\
		
		We need to verify Assumptions \ref{ass:1}-\ref{ass:factorisation} for the map \eqref{eq:Poincare_map}, as well as the normal hyperbolicity of $S_\Sigma$. Notice first that \eqref{eq:Poincare_map} can be rewritten in the general form \eqref{eq:gen_maps} after defining $\widehat f(x,\alpha,\eps) =: \widehat f(x,\alpha,0) + f_{rem}(x,\alpha,\eps)$ and
		\[
		N(x,\alpha) := 
		\begin{pmatrix}
			I_{n-1}  \\
			O_{1,n-1}
		\end{pmatrix}
		, 
		\qquad
		f(x,\alpha) := \widehat f(x,\alpha,0) ,
		\qquad
		\tilde G(x,\alpha,\eps) := 
		\begin{pmatrix}
			f_{rem}(x,\alpha,\eps) \\
			g(x,\alpha,\eps)
		\end{pmatrix}
		.
		\]
		We have that
		\[
		f(x_\ast,\alpha_\ast) = \Phi_1(x_\ast,0,\alpha_\ast,0,T) - x_\ast = 
		\Phi_1(x_\ast,0,\alpha_\ast,0,0) - x_\ast = x_\ast - x_\ast = 0,
		\]
		and, by Assumption \ref{ass:limit_cycles}, that the $n-1$ characteristic multipliers $\mu_j$ of the matrix $I_{n-1} + D_x f(x_\ast,\alpha_\ast)$ satisfy $|\mu_j| \neq 1$. It follows that the $n-1$ eigenvalues $\lambda_j$ of the matrix $D_x f(x_\ast,\alpha_\ast)$ satisfy $\text{Re}\ \lambda_j \neq 0$, i.e.~it is regular, so that the implicit function theorem implies the existence of a $C^r-$smooth function $\tilde \varphi_0 : V_\alpha \to V_x$ such that $\tilde \varphi_0(\alpha_\ast) = x_\ast$ and $f(\varphi_0(\alpha),\alpha) \equiv O_{n-1}$ on $V_\alpha$. The fact that $\tilde \varphi_0 = \varphi_0$, i.e.~that $S_\Sigma = \{(\varphi_0(\alpha),\alpha) : \alpha \in V_\alpha\} = \{(x,\alpha) : f(x,\alpha) = O_{n-1}\}$ defines the critical manifold of \eqref{eq:Poincare_map}, follows from the local uniqueness of $\varphi$ and $\tilde \varphi_0$ and their coincidence at $(x_\ast,\alpha_\ast)$.
		Thus, \eqref{eq:Poincare_map} is a fast-slow map satisfying Assumption \ref{ass:1}. Assumption \ref{ass:factorisation} is immediate since we consider a neighbourhood with $(x,\alpha) \in V_x \times V_\alpha$ within which the only zeroes of $N(x,\alpha) f(x,\alpha)$ are those identified above for $f(x,\alpha) = O_{n-1}$ (uniqueness follows from the implicit function theorem).
		
		\
		
		It remains to verify that $S_\Sigma$ is normally hyperbolic. By Proposition \ref{prop:EVs}, non-trivial multipliers along $S_\Sigma$ are in $1-1$ correspondence with eigenvalues of the $(n-1) \times (n-1)$ square matrix
		\begin{equation}
			\label{eq:monodromy}
			I_{n-1} + N(\varphi_0(\alpha),\alpha) Df(\varphi_0(\alpha),\alpha) = I_{n-1} + D_x f(\varphi_0(\alpha),\alpha) .
		\end{equation}
		For each fixed $\alpha \in V_\alpha$ the $n-1$ multipliers $\mu_j(\alpha)$ of the matrix \eqref{eq:monodromy} are $1-1$ correspondence with the characteristic multipliers associated a hyperbolic limit cycle $\Gamma_\alpha \subset \mathcal M$, and therefore satisfy $|\mu_j(\alpha)| \neq 1$ for all $\alpha \in V_\alpha$. Thus $S_\Sigma$ is normally hyperbolic.
	\end{proof}

	
	
	%

	Lemma \ref{lem:Poincare_map} implies that the results of Section \ref{sec:slow_manifold_theorems} can be applied directly to the Poincar\'e map $\mathcal P_\Delta$. The derivation of the Poincar\'e map on $\Delta = \Sigma$ generalises straightforwardly to the case where $\Delta$ is given by the graph \eqref{eq:Delta}; one must simply restrict to $y = Y(x,\alpha)$ instead in the defining expressions \eqref{eq:Poincare_map_2}. Explicitly,
	\begin{equation}
		\label{eq:Poincare_map_3}
		\mathcal P_\Delta : 
		\begin{pmatrix}
			x \\
			\alpha
		\end{pmatrix}
		\mapsto
		\begin{pmatrix}
			\bar x \\
			\bar \alpha
		\end{pmatrix}
		=
		\begin{pmatrix}
			x \\
			\alpha
		\end{pmatrix}
		+
		\begin{pmatrix}
			\widehat f(x,\alpha,\eps) \\
			\eps g(x,\alpha,\eps)
		\end{pmatrix}
		,
	\end{equation}
	with
	\begin{equation}
		\label{eq:Poincare_map_4}
		\begin{split}
			\widehat f(x,\alpha,\eps) &:= \Phi_1(x,Y(x,\alpha),\alpha,\eps,\tilde t(x,Y(x,\alpha),\alpha,\eps)) - x , \\
			g(x,\alpha,\eps) &:= 
			D_\eps \Phi_3(x,Y(x,\alpha),\alpha,0,\tilde t(x,Y(x,\alpha),\alpha,0)) + O(\eps) ,
		\end{split}
	\end{equation}
	where in particular we have
	\begin{equation}
		\label{eq:averaged_eqn}
		g(x,\alpha,\eps) = \int_0^{\tilde t(x,Y(x,\alpha),\alpha,0)} \tilde g \left( \Phi_1(x,Y(x,\alpha),\alpha,0,s), Y(x,\alpha), \alpha, 0 \right) ds + O(\eps) .
	\end{equation}
	
	\
	
	Our main results on the persistence of the critical manifold $S_\Delta$, its local stable/unstable manifolds $W^{s/u}_{loc}(S_\Delta)$ and their foliations for $0 < \eps \ll 1$, are stated for the map $\mathcal P_\Delta$ defined in equations \eqref{eq:Poincare_map_3}-\eqref{eq:Poincare_map_4}.
	

	\begin{thm}
		\label{thm:Poincare_map}
		Consider the Poincar\'e map $\mathcal P_\Delta$ defined by equations \eqref{eq:Poincare_map_3}-\eqref{eq:Poincare_map_4}. There exists and $\eps_0 > 0$ such that following properties hold for all $\eps \in (0,\eps_0)$:
		\begin{enumerate}
			\item[(i)] The critical manifold $S_\Delta$ persists as a locally invariant slow manifold $S_{\Delta,\eps}$ described by Theorems \ref{thm:slow_manifolds}-\ref{thm:graph_slow_manifolds}. In particular,
			\[
			S_{\Delta,\eps} = \{(\varphi_\eps(\alpha), \alpha) : \alpha \in V_\alpha \} ,
			\]
			where
			\[
			\varphi_\eps(\alpha) = \varphi_0(\alpha) - \eps \left( (D_xf) f_{rem} + (D_\alpha g) \right) \big|_{S_\Delta} + O(\eps^2) ,
			\]
			and slow iteration along $S_\eps$ is governed by
			\begin{equation}
				\label{eq:slow_Poincare}
				\mathcal P_\Delta|_{S_{\Delta,\eps}} : 
				\begin{pmatrix}
					x \\
					\alpha
				\end{pmatrix}
				\mapsto 
				\begin{pmatrix}
					\bar x \\
					\bar \alpha
				\end{pmatrix}
				=
				\begin{pmatrix}
					\varphi_0(\alpha) \\
					\alpha
				\end{pmatrix}
				+ \eps
				\begin{pmatrix}
					- (D_xf)^{-1} (D_\alpha f) g \\
					g
				\end{pmatrix}
				\bigg|_{S_\Delta} + O(\eps^2) ,
			\end{equation}
			as long as iterates remain in $V_x \times V_\alpha$.
			
			\item[(ii)] The stable/unstable manifolds $W_{loc}^{s/u}(S_\Delta)$ persist as locally invariant manifolds $W^{s/u}_{loc}(S_{\Delta,\eps})$ described by Theorem \ref{thm:stable_manifolds}.
			
			\item[(iii)] Foliations of the stable/unstable manifolds $W_{loc}^{s/u}(S_\Delta)$ persist as locally invariant foliations of $W_{loc}^{s/u}(S_{\Delta,\eps})$ described by Theorem \ref{thm:foliations}. The spectral bounds $\nu_A$ and $\nu_R$ determining the contraction and repulsion rates in Theorem \ref{thm:foliations} (iii) and (iv) respectively are given by
			\[
			\nu_A = \sup_{\alpha \in V_\alpha, j=1,\ldots,n-1} \{ |\mu_j(\alpha)| : |\mu_j(\alpha)| < 1 \} , \ \ 
			\nu_R = \sup_{\alpha \in V_\alpha, j=1,\ldots,n-1} \{ |\mu_j(\alpha)| : |\mu_j(\alpha)| > 1 \} ,
			\]
			where $\mu_j(\alpha)$, $j=1,\ldots,n-1$ denote the characteristic multipliers associated to the limit cycle $\Gamma_\alpha \subset \mathcal M$ in the layer problem \eqref{eq:stnd_limit_cycles_layer}, which are in $1-1$ correspondence with 
			multipliers of the matrix $I_{n-1} + D_xf{(\varphi_0(\alpha),\alpha)}$.
		\end{enumerate}
	\end{thm}
	
	\begin{proof}
		It follows by Lemma \ref{lem:Poincare_map} that the results of Sections \ref{sec:a_coordinate-independent_framework_for_fast-slow_maps}-\ref{sec:slow_manifold_theorems} can be applied directly to the map $P_\Delta$.
		
		The existence and form of the parameterization of the slow manifold $S_{\Delta,\eps}$ is assertions (i) follows directly from Theorems \ref{thm:slow_manifolds}-\ref{thm:graph_slow_manifolds}, and the form of the map describing slow iteration along $S_{\Delta,\eps}$ follows from Proposition \ref{prop:reduced_map}.
		
		Assertions (ii) and (iii) follow directly from Theorem \ref{thm:stable_manifolds} and Theorem \ref{thm:foliations} respectively.
	\end{proof}
	
	
	\begin{figure}[t!]
		\centering
		\includegraphics[scale=0.12]{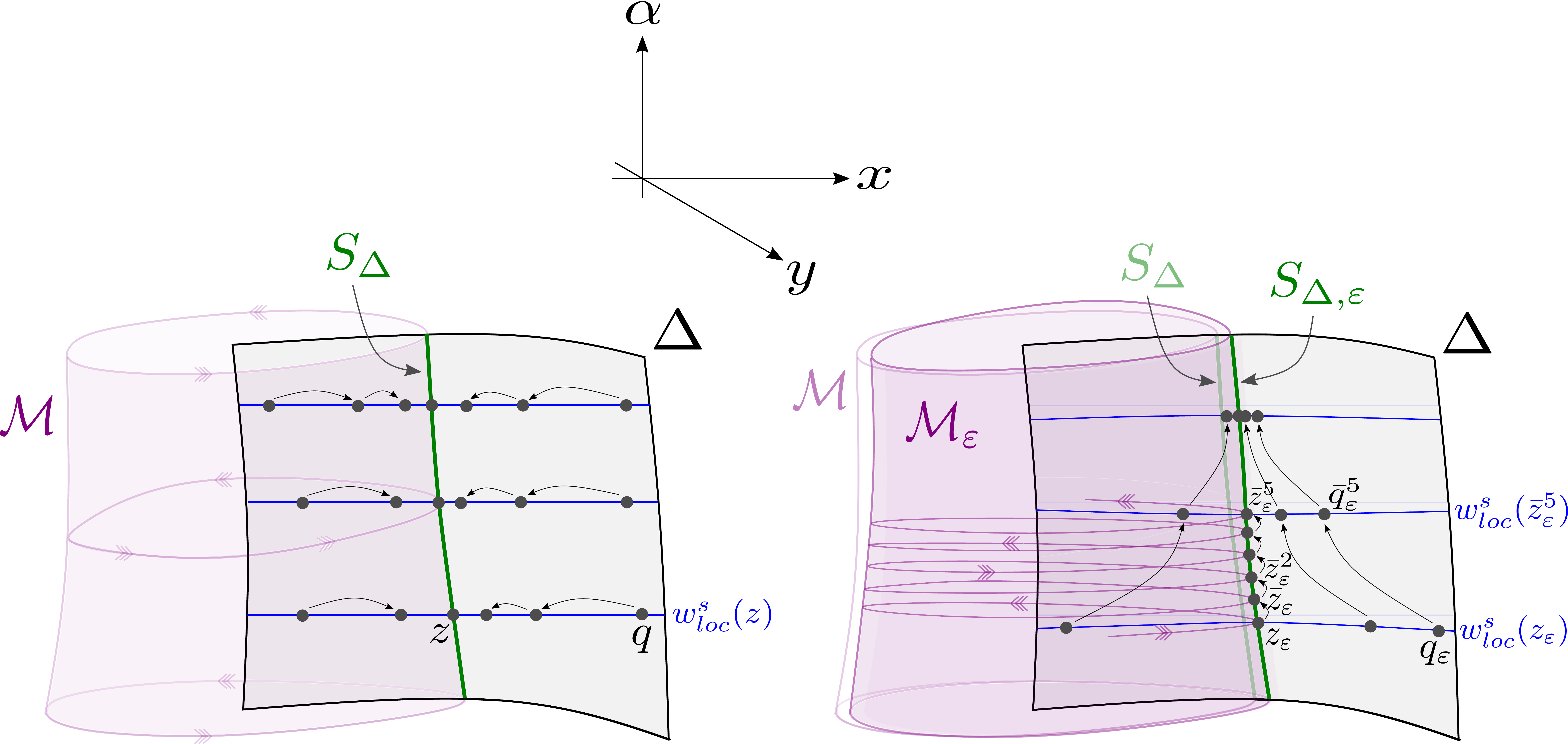}
		\caption{Perturbation of the singular geometry and dynamics for $\eps = 0$ (left) for $0 \leq \eps \ll 1$ (right). We sketch the 3-dimensional case of a normally hyperbolic and attracting limit cycle manifold $\mathcal M$. 
			Theorem \ref{thm:Poincare_map} describes the persistence properties for the Poincar\'e map $P_{\Delta}$ on an arbitrary cross-section $\Delta$. 
			Three representative stable fibers on $\Delta$ for $\eps = 0$ are shown in blue on the left, which perturb to invariant stable fibers for $0 < \eps \ll 1$ (also in blue) on the right, e.g.~$w_{loc}^s(z)$ perturbs to $w_{loc}^s(z_\eps)$. 
			A number of example iterates are shown on $\Delta$ for both $\eps = 0$ and $0 < \eps \ll 1$. Note that the arrows connecting points on e.g.~$w_{loc}^s(z_\eps)$ and $w_{loc}^s({\bar z_\eps}^5)$ only show the direction of iteration (they do not indicate a single iteration). The persistence of the manifold of limit cycles $\mathcal M$ as a nearby locally invariant `slow manifold' $\mathcal M_\eps$ as described by Theorem \ref{thm:Poincare_map} is also shown ($\mathcal M_\eps$ is shown in a darker shade or purple).}
		\label{fig:Poincare_map}
	\end{figure}
	
	The situation is sketched in Figure \ref{fig:Poincare_map}.
	
	\
	

	Our final result for this section characterises the persistence of the limit cycle manifold $\mathcal M$ in the ODE system \eqref{eq:stnd_fast_slow_limit_cycles} 
	for $0 < \eps \ll 1$. Our findings are consistent with previous results due to Anosova \cite{Anosova1999,Anosova2002}. We also provide a criterion for the existence of hyperbolic limit cycles in the perturbed system \eqref{eq:stnd_fast_slow_limit_cycles}.
	
	\begin{thm}
		\label{thm:limit_cycle_manifolds}
		Consider the fast-slow ODE system \eqref{eq:stnd_fast_slow_limit_cycles} under Assumption \ref{ass:limit_cycles}. Then there exists and $\eps_0 > 0$ such that for all $\eps \in (0,\eps_0)$, 
			the limit cycle manifold $\mathcal M$ perturbs to a locally invariant manifold $\mathcal M_\eps$ which is $O(\eps)-$close, $C^r-$smooth and diffeomorphic to $\mathcal M$.
			
			The manifold $\mathcal M_\eps$ contains a hyperbolic limit cycle of the ODE \eqref{eq:stnd_fast_slow_limit_cycles} intersecting the point with local coordinates $(x_{lc}, y_{lc}, \alpha_{lc}) \in \mathcal M_\eps$ if and only if
			\begin{equation}
				\label{eq:limit_cycle_conds}
				g(\varphi_0(\alpha),\alpha,0) = 0, \qquad \text{and} \qquad
				D_\alpha g(\varphi_0(\alpha),\alpha,0) \neq 0 ,
			\end{equation}
			where $g$ is defined in the $(x,y,\alpha)-$coordinates of the Poincar\'e map \eqref{eq:Poincare_map_4} via the integral formula \eqref{eq:averaged_eqn}.
			%
			%
			
	\end{thm}
	
	\begin{proof}
		The limit cycle manifold $\mathcal M$ can be parameterised by a one-parameter family of submanifolds $S_\sigma$ defined via the flow of the layer problem \eqref{eq:stnd_limit_cycles_layer}. Specifically, define
		\[
		S_\sigma := \Phi_{\sigma T_\alpha}(S)
		= \{ \Phi(\varphi_0(\alpha),Y(\varphi_0(\alpha),\alpha),\alpha,0,\sigma T_\alpha) : \alpha \in V_\alpha \}  , \quad \sigma \in [0,1) ,
		\]
		where we set $S_0 := S = \mathcal M \cap \Delta$, such that $\mathcal M = \cup_{\sigma \in [0,1)} S_\sigma$. Since it is always possible to associate a non-degenerate cross-section $\Delta_\sigma$ such that $S_\sigma = \Delta_\sigma \cap \mathcal M$, a family of Poincar\'e maps $\{\mathcal P_{\Delta_{\sigma}} : \sigma \in [0,1) \}$ can be obtained by constructing each $\mathcal P_{\Delta_\sigma} : \Delta_\sigma \to \Delta_\sigma$ analogously to the map \eqref{eq:Poincare_map} described in Theorem \ref{thm:Poincare_map}. In particular, $S_\sigma$ is a normally hyperbolic critical manifold for $\mathcal P_\sigma$, for each $\sigma \in [0,1)$.
		
		By Theorem \ref{thm:Poincare_map} (i), each $S_\sigma$ perturbs to a nearby slow manifold $S_{\sigma,\eps} \subset \Delta_\sigma$. Thus for any fixed choice of $\sigma$, we may define the manifold $\mathcal M_\eps$ by
		\[
		\mathcal M_\eps := \Phi\left( S_{\sigma,\eps}, \eps, [-\tilde T, \tilde T] \right) \cap \{(x,y,\alpha) : \alpha \in V_\alpha \}
		\]
		where $\tilde T > 0$ is fixed larger than the maximum of $T_{\alpha,\eps}$ on $V_\alpha \times [0,\eps_0)$. 
		Since the flow map $\Phi$ is $C^r-$smooth, $\mathcal M_\eps$ is exactly as smooth as $S_{\sigma,\eps}$, which by Theorem \ref{thm:Poincare_map} is exactly as smooth as $S_\sigma$, which is exactly as smooth as $\mathcal M$. The fact that $\mathcal M_\eps$ is $O(\eps)-$close and diffeomorphic to $\mathcal M$ follows from the corresponding results for $S_{\sigma,\eps}$ and $S_\sigma$ implied by Theorem \ref{lem:Poincare_map} (i) and regular perturbation theory.
		
		\
		
		It remains to prove the statement pertaining to limit cycles in $\mathcal M_\eps$. We work in local coordinates $(x,y,\alpha) \in V_x \times V_y \times V_\alpha$ for which we have $x_{lc} = \varphi_\eps(\alpha_{lc})$ and $y_{lc} = Y(\varphi_0(\alpha_{lc}),\alpha_{lc})$ since we may assume that $(x_{lc},y_{lc},\alpha_{lc}) \in S_{\Delta,\eps} = \Delta \cap \mathcal M_\eps$. Limit cycles of \eqref{eq:stnd_fast_slow_limit_cycles} in $\mathcal M_\eps$ are in 1$-$1 correspondence with fixed points of the restricted Poincar\'e map on $P_\Delta|_{S_{\Delta,\eps}}$ given by \eqref{eq:slow_Poincare}. Since the matrix $D_xf$ is locally regular under normally hyperbolic conditions, fixed points and their hyperbolicity are completely determined by the restricted, 1-dimensional map
		\begin{equation}
			\label{eq:Poincare_slow_reduced}
			\alpha \mapsto \bar \alpha = \alpha + \eps g(\varphi_\eps(\alpha),\alpha,\eps)
			= \alpha + \eps g(\varphi_0(\alpha),\alpha,0) + O(\eps^2) .
		\end{equation}
		It follows from the implicit function theorem that the map \eqref{eq:Poincare_slow_reduced} has a hyperbolic fixed point for all $0 < \eps \ll 1$ sufficiently small if and only if the conditions in \eqref{eq:limit_cycle_conds} are satisfied. This concludes the proof.
	\end{proof}
	
	Figure \ref{fig:Poincare_map} shows both $\mathcal M$ and it's perturbation $\mathcal M_\eps$, which intersects $\Delta$ along the slow manifold $S_{\Delta,\eps}$ of the Poincar\'e map $P_\Delta$.
	
	\begin{remark}
		\label{rem:averaging}
		Equation \eqref{eq:averaged_eqn} implies an integral formulation of the conditions in \eqref{eq:limit_cycle_conds} in terms of the expression
		\begin{equation}
			\label{eq:averaging_eqn_2}
			g(\varphi_0(\alpha),\alpha,0) = \int_0^{T_\alpha} \tilde g \left( \Phi_1(\varphi_0(\alpha),Y(\varphi_0(\alpha),\alpha),\alpha,0,s), Y(\varphi_0(\alpha),\alpha), \alpha, 0 \right) ds ,
		\end{equation}
		where $T_\alpha = \tilde t(\varphi_0(\alpha),Y(\varphi_0(\alpha),\alpha),\alpha,0)$ is the period associated to the limit cycle $\Gamma_\alpha \subset \mathcal M$ of the layer problem \eqref{eq:stnd_limit_cycles_layer} which passes through the point $(\varphi_0(\alpha), Y(\varphi_0(\alpha),\alpha), \alpha) \in S_\Delta$. In classical averaging theory, the integral \eqref{eq:averaging_eqn_2} defines the so-called \textit{averaged equation} \cite{Guckenheimer1983,Sanders2007}.
	\end{remark}

	\begin{remark}
		\label{rem:isochrons}
		Although the contraction and expansion along invariant foliations under the Poincar\'e map $P_{\Delta,\eps}$ are described in detail by Theorem \ref{thm:Poincare_map}, the persistence of stable/unstable manifolds $W^{u/s}_{loc}(\mathcal M)$ and their corresponding foliations in the ODE \eqref{eq:stnd_fast_slow_limit_cycles} is not described in Theorem \ref{thm:limit_cycle_manifolds}. A detailed proof of persistence is left for future work. 
	\end{remark}

	\section{Outlook}
	\label{sec:outlook}
	
	Discrete multi-scale dynamical systems induced by maps arise in a wide variety applied and theoretical settings. However, to the best of our knowledge, a complete geometric theory for their analysis analogous to the established GSPT for multi-scale continuous-time systems does not yet exist. The aim of the present manuscript is to provide such a theory, referred to herein as \textit{DGSPT}, for the class of fast-slow maps defined by \eqref{eq:nonstandard_maps} or, more precisely, general maps \eqref{eq:gen_maps} under Assumptions \ref{ass:1} and \ref{ass:factorisation}. 
	
	\
	
	The singular theory is presented in Section \ref{sec:a_coordinate-independent_framework_for_fast-slow_maps}, where the layer map is introduced and used to define a discrete notion of normal hyperbolicity in terms of a spectral condition on the non-trivial multipliers along the critical manifold $S$, recall Definition \ref{def:nh}. By Proposition \ref{prop:EVs}, the problem of checking for normal hyperbolicity of $S$ reduces to a direct evaluation of the spectrum associated to the matrix $I_{n-k} + Df N|_S$, which is given purely in terms of $N$ and $f$. Under normally hyperbolic assumptions, we introduced the reduced map \eqref{eq:red_map} in Definition \ref{def:red_map}. Due to the fact that there is no discrete analogue of the time rescaling $\tau = \eps t$ in the continuous-time setting, the reduced map depends necessarily on $\eps$ in the leading order. By Proposition \ref{prop:reduced_map}, however, it nevertheless approximates slow iteration along perturbed slow manifolds $S_\eps$ up to an accuracy of $O(\eps^2)$, which is typically all that is required in applications. An expression for the $m$'th iterate map on $S_\eps$ was also derived using the asymptotic self-similarity of the map near $S$, see again Proposition \ref{prop:red_jump_map}. Interestingly, this map can to leading order in $\eps$ be related to the Euler discretization of the corresponding continuous-time reduced problem. It is also possible that the reduced $m$'th iterate map is well-defined on $S$ under a suitable formulation of the dual limit $(\eps,m) \to (0,\infty)$, recall Remark \ref{rem:mth_map}, however the details remain for future work.
	
	The persistence of key dynamic and geometric features for $0 < \eps \ll 1$ under normally hyperbolic conditions is described in detail by the Fenichel-like theorems provided in Section \ref{sec:slow_manifold_theorems}. The persistence of a compact, normally hyperbolic critical manifold as a nearby locally invariant slow manifold $S_\eps$ is described by Theorems \ref{thm:slow_manifolds} and \ref{thm:graph_slow_manifolds}, with the latter providing an explicit parameterization for $S_\eps$ up to $O(\eps^2)$ in local coordinates. The persistence of local stable/unstable manifolds $W_{loc}^{s/u}(S)$ and their corresponding asymptotic rate foliations are described by Theorems \ref{thm:stable_manifolds} and \ref{thm:foliations} respectively. Theorem \ref{thm:foliations} in particular gives quantitative estimates for the contraction and expansion rates close to $S_\eps$ in terms of the scalar quantities $\nu_A$ and $\nu_R$ defined in \eqref{eq:spectral_bounds_original}, which define the (annular) spectral gap about the unit circle associated to the matrix $I_{n-k} + Df N|_S$. While the invariant manifold theorems in Section \ref{sec:slow_manifold_theorems} are in some regards less general than their counterparts in \cite{Nipp2013} (e.g.~we require $C^1-$smoothness as opposed to Lipschitz continuity), the `Fenichel-like' formulation of the results in Section \ref{sec:slow_manifold_theorems} are expected to be advantageous in many applications due the fact that they (i) apply directly to fast-slow maps \eqref{eq:nonstandard_maps} without the need for prior transformations into special coordinates, and (ii) apply directly near compact submanifolds of $S$ in the locally invariant setting.
	
	
	The utility of DGSPT was demonstrated in Section \ref{sec:examples} in the context of three non-trivial applications. In Section \ref{sub:Chialvo_model} we used DGSPT in order to provide a geometric analysis of the Chialvo map-based model for neuronal bursting in a fast-slow parameter regime, thereby validating and extending previous work in \cite{Chialvo1995,Jing2006,Trujillo2021}. The DGSPT formalism provided a systematic approach to the identification of distinct singular structure corresponding to four types of neuronal dynamics: excitability, relaxation, regular (non-chaotic) bursting and chaotic bursting, recall Figure \ref{fig:Chialvo_global_dynamics}. In Section \ref{sub:Euler_discretization} we used DGSPT to 
	prove results on the relationship between the geometry and dynamics of fast-slow ODEs in the general (non-standard) form \eqref{eq:general_ode} and its corresponding Euler discretization, see Theorem \ref{thm:Euler}. These results extend pre-existing results on the discretization of fast-slow ODEs in standard form \eqref{eq:stnd_form_ode} in e.g.~\cite{Nipp1995,Nipp1996}. Finally in Section \ref{sub:Poincare_maps}, we showed that fast-slow ODE systems with a single slow variable (which often arise in standard form \eqref{eq:stnd_form_ode} after allowing a parameter to evolve slowly in time) give rise to fast-slow Poincar\'e maps with a normally hyperbolic critical manifold if the layer problem has a hyperbolic limit cycle, see Lemma \ref{lem:Poincare_map}. This critical manifold, which lies within the intersection of a limit cycle manifold $\mathcal M$ of the (continuous-time) layer problem and a transverse section $\Delta$, perturbs to a nearby slow manifold for the Poincar\'e map for $0 < \eps \ll 1$. The geometry and dynamics of the Poincar\'e map are described in detail in Theorem \ref{thm:Poincare_map}, which follows after a direct application of the invariant manifold theorems in Section \ref{sec:slow_manifold_theorems}. Properties of the Poincar\'e map were then applied in order to prove Theorem \ref{thm:limit_cycle_manifolds}, which describes the persistence of the limit cycle manifold $\mathcal M$ as a nearby locally invariant manifold $\mathcal M_\eps$ in the ODE and, additionally, provides a criteria for the existence and hyperbolicity of limit cycles within $\mathcal M_\eps$.
	
	\
	
	The present work aims to provides a mathematical framework for the analysis of a large class of multi-scale discrete dynamical systems, thereby paving the way for a wealth of future work. In particular, the question of how to supplement DGSPT with geometric methods for studying the non-normally hyperbolic regime frequently arose in this work. The singular theory of Section \ref{sec:a_coordinate-independent_framework_for_fast-slow_maps} provides the means for classifying non-normally hyperbolic singularities on $S$ in terms of `$\eps = 0$ conditions' on the layer map \eqref{eq:layer_map} and the function $\mathcal G(z) : = \Pi_{\mathcal N}^S G(z,0)|_S$ which determines the dynamics of the reduced map \eqref{eq:red_map}. A detailed treatment of the codimension-1 singularities, i.e.~the fast-slow fold, flip/period-doubling and Neimark-Sacker/torus singularities, in the low dimensional setting using DGSPT in combination with adaptations of the geometric blow-up technique constitutes very recent, ongoing and future work.

	\section{Acknowledgements}
	
	SJ and CK acknowledge funding from the SFB/TRR 109 Discretization and Geometry in Dynamics grant. SJ would like to thank Prof.~Martin Wechselberger for many helpful discussions pertaining to geometric singular perturbation theory ``beyond the standard form" in general. CK thanks the VolkswagenStiftung for support via a Lichtenberg Professorship.

	\bibliographystyle{siam}
	\bibliography{dgspt}

	\unmarkedfntext{The corresponding author S.~Jelbart can be contacted via jelbart@ma.tum.de}
	
	\unmarkedfntext{C.~Kuehn can be contacted via ckuehn@ma.tum.de}

\end{document}